\documentclass[11pt,letterpaper]{amsart}
\usepackage[left=2.2cm,tmargin=2.5cm,bmargin=2.5cm,right=2.2cm]{geometry}

\usepackage{graphicx} 
\usepackage{tikz}
\usepackage{amssymb}
\usepackage{amsmath}
\usepackage{amsthm}
\usepackage{tikz-cd}
\usepackage{array}
\usepackage{mathtools}
\usepackage{adjustbox}
\usepackage{float}
\usetikzlibrary {positioning}
\usepackage[raggedrightboxes]{ragged2e}
\usepackage{biblatex}
\newcolumntype{P}[1]{>{\centering\arraybackslash}p{#1}}

\counterwithin{equation}{section}
\counterwithin{figure}{section}
\newtheorem{theorem}{Theorem}[section]

\newtheorem{lemma}[theorem]{Lemma}
\newtheorem{proposition}[theorem]{Proposition}
\newtheorem{corollary}[theorem]{Corollary}

\newcommand{\Z}{\mathbb{Z}}
\newcommand{\Q}{\mathbb{Q}}
\newcommand{\R}{\mathbb{R}}
\newcommand{\res}[2]{\left(\frac{#1}{#2}\right)}
\newcommand{\C}{\mathcal{C}}

\newcommand{\A}{\mathcal{A}}
\newcommand{\B}{\mathcal{B}}

\title{The Local-Global Conjecture is False for Generalized Circle Packings}
\author{Hanqi Shi, Wenyuan Shi, Ian Whitehead, \\ Ham Williams-Tracy, Jeffrey Zhirui Zhang}

\begin{document}

\begin{abstract}
Haag, Kertzer, Rickards, and Stange disprove the Local-Global Conjecture for Apollonian circle packings. We extend their disproof to four more types of integral circle packing: the octahedral, cubic, square, and triangular packings. In each case, we find quadratic invariants which imply quadratic reciprocity obstructions to the conjecture in certain packings. We utilize an explicit parametrization of circles tangent to a fixed circle in each packing type, and a quadratic reciprocity argument. Even in the packings where we do not find quadratic obstructions, the curvatures exhibit a predictable reciprocity structure. This leads to partial obstructions on integers appearing as curvatures in subsets of the packing.  
\end{abstract}

\maketitle

\section{Introduction}
The Local-Global Conjecture for Apollonian circle packings was disproven in stunning 2024 work of Haag, Kertzer, Rickards and Stange \cite{HKRS}. Their work uses classical tools, but a wholly original strategy. It opens up a new field of inquiry: reciprocity obstructions to integer points in group orbits, like the Brauer-Manin obstruction to integer points on varieties. 

Apollonian packings are one species within a broader taxonomy of fractal circle packings with number theoretic structure. In this article, we extend the methods of Haag, Kertzer, Rickards, and Stange to four other circle packing types. We find that their strategy generalizes naturally to these cases. We disprove the analogue of the Local-Global Conjecture for all four packing types.  This suggests that there is much more to the story of these reciprocity obstructions.

\begin{figure}[h]
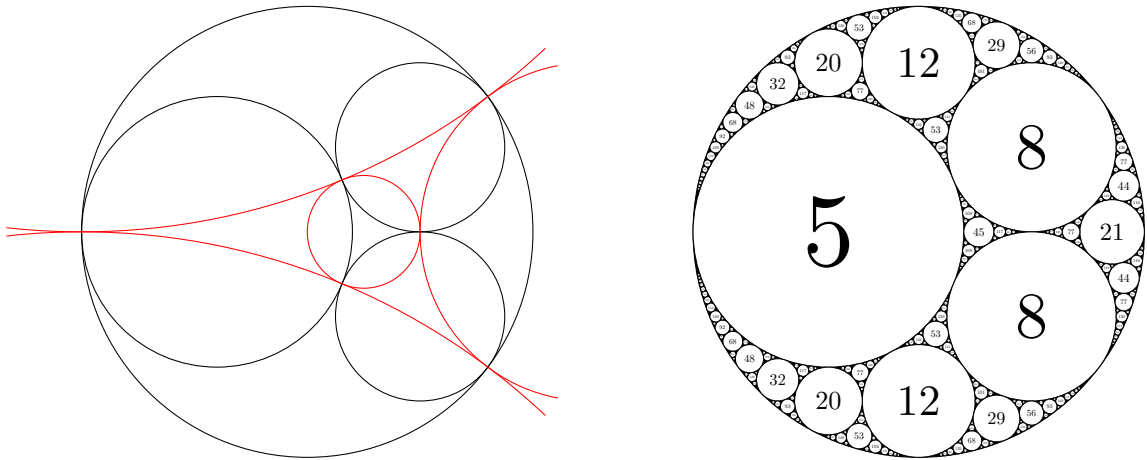
 
\begin{center}
\vspace{-.5in}
\include{apollonianfigure}
\end{center}
\vspace{-.75in}
\caption{Base and Dual Circle Configurations; Apollonian Packing}
\label{ApollonianFig}
\end{figure}

To construct an Apollonian packing, we begin with a base configuration of four circles all tangent to one another. It is possible to draw a dual configuration of four more circles which cover the interstitial regions in the base configuration. Any intersections between base circles and dual circles are orthogonal at shared points of tangency. The Apollonian group is generated by the reflections through the dual circles. The orbit of the four base circles under the Apollonian group is an infinite collection of circles, pairwise tangent or disjoint, called an Apollonian packing. 

The number theory of Apollonian packings begins with the Descartes Circle Theorem. If four circles are mutually tangent, then their curvatures $a$, $b$, $c$, $d$ satisfy the quadratic equation
\begin{equation} \label{Descartes}
(a+b+c+d)^2=2(a^2+b^2+c^2+d^2)
\end{equation}
The curvature of a circle is the reciprocal of the radius. We consider a circle as having a negative curvature if it contains the others; we can also consider a straight line as a circle of curvature $0$. The Apollonian group acts on quadruples of curvatures. A reflection in the group preserves three curvatures, say $a$, $b$, $c$ of circles which intersect the dual circle orthogonally, and transforms the fourth curvature $d$ to $d'=2a+2b+2c-d$, the other solution of the quadratic equation \eqref{Descartes}. Based on this interpretation of the group action, we can see that if we start with four integer curvatures, we will obtain an integral packing where every curvature is an integer. An integral packing is called primitive if the greatest common divisor of the curvatures is $1$. It is natural to ask which integers appear as curvatures in a given primitive integral packing.

The Local-Global Conjecture, first posed in \cite{GLMWY} and further developed in \cite{FuchsSanden, Fuchs}, was an attempt to answer this question. There are modular restrictions on the curvatures that can appear in a given Apollonian packing: for example, in the packing shown in Figure \ref{ApollonianFig}, no curvatures congruent to 2 or 3 mod 4 and no curvatures congruent to 1 mod 3 appear. The specific restrictions depend on the packing, but they can be expressed as congruences mod 24 for all Apollonian packings.  The Local-Global Conjecture predicts that every sufficiently large integer subject to the modular restrictions will appear in a given primitive integral packing. Bourgain and Kontorovich proved that asymptotically $100\%$ of these integers appear, but this still leaves open the possibility of a zero-density but infinite set of excluded integers \cite{BourgainKontorovich}.

Haag, Kertzer, Rickards, and Stange find such a set. We will describe their argument using the packing $\A$ in Figure \ref{ApollonianFig} as an example. This is an especially simple case, but illustrative of the general approach. They construct a quadratic invariant $\chi_2(\A)$ for the packing. In this case, the invariant takes the value of the Kronecker symbol $\res{b}{a}$ for any pair of tangent circles with coprime curvatures $a$ and $b$. Observe that in Figure 1, $\res{b}{a}$ appears to be $-1$ for all such pairs. If we can prove that the invariant is well-defined, then it follows that no perfect square integers appear as curvatures, despite the fact that infinitely many square numbers satisfy the modular restrictions. The set of square integers is a quadratic obstruction to the Local-Global Conjecture.

The proof that $\chi_2(\A)=\res{b}{a}$ takes the same value for all tangent, coprime pairs $a,b$ has two main steps. The first step is to show that for a fixed circle of curvature $a$, the Kronecker symbol $\res{b}{a}$ takes the same value for any tangent circle of coprime curvature $b$. This uses a theorem of Sarnak: that all tangent curvatures $b$ are parametrized by a shifted binary quadratic form $Q(x,y)-a$ for $x, y \in \Z$ coprime \cite{Sarnak}. For example, in Figure \ref{ApollonianFig}, the curvatures of circles tangent to the circle of curvature $5$ are parametrized by $2x^2+2xy+13y^2-5$. Moreover, the discriminant of $Q(x,y)$ is divisible by $a$. Therefore, $Q(x,y)$ factors as a constant times a perfect square in $\Z/a\Z$, and all of its values coprime to $a$ are equal up to squares in $(\Z/a\Z)^*$. This establishes that the quadratic invariant is well-defined for $a$ fixed and $b$ variable.

The second step of varying $a$ uses quadratic reciprocity. If we have tangent circles of coprime curvatures $a$ and $b$ in the packing above, then $\res{b}{a} = \res{a}{b}$. This simple form of quadratic reciprocity is valid because all circles have curvatures congruent to $0$ or $1$ mod $4$. For other packings, the definition of the quadratic invariant is more complicated, but the modular restrictions always conspire to make a reciprocity argument work. The quadratic invariant is the same for the circle of curvature $b$ as it is for the circle of curvature $a$. For any two circles in a given primitive, integral packing, Haag, Kertzer, Rickards, and Stange show that there is a path of pairwise tangent, coprime circles between them. This completes their proof that $\chi_2(\A)$ is well-defined across the packing.

A quadratic invariant of this type exists for every primitive integral Apollonian packing. Sometimes the invariant does not lead to a quadratic obstruction--e.g. if $\chi_2(\A)=1$ then there is no obstruction to squares appearing as curvatures. But the packings with an obstruction are counterexamples to the Local-Global Conjecture. Additionally, Haag, Kertzer, Rickards and Stange find quartic invariants, which are well-defined for packings satisfying some congruence conditions mod 24 but not others. These lead to further obstructions and counterexamples to the conjecture. 

We now introduce the circle packings studied in this article. Instead of starting with four mutually tangent circles, we start with a base configuration of circles whose tangency relations are described using a planar graph. For each packing type, we define an Apollonian group generated by reflections across a dual circle configuration. The orbit of the base configuration under the Apollonian group is an infinite fractal collection of circles, pairwise tangent or disjoint. The particular planar graphs we use are the graphs of the octahedron, the cube, the square tiling of $\R^2$, and the triangular tiling of $\R^2$. We call the four resulting circle packings the octahedral, cubic, square, and triangular packings. The first two are examples of polyhedral packings, defined by Kontorovich and Nakamura \cite{KontorovichNakamura}. The latter two belong to a family of packings constructed from tilings of the plane, defined by Rehwinkel, Yang, Yang, and the third author \cite{RWYY}.

All four circle packing types have interesting number theoretic properties. They all  include integral packings, where all the curvatures are integers, subject to modular restrictions depending on the packing. They are all superintegral, crystallographic packings in the terminology of \cite{KontorovichNakamura}. Furthermore, they are all familial packings in the terminology of \cite{FuchsStangeZhang}. In particular, each Apollonian group is conjugate to a subgroup of $\mathrm{SL}(2, \mathcal{O}_K)$ for $K=\Q(\sqrt{-d})$ an imaginary quadratic field ($d=1$ for the square packing, $d=2$ for the cubic and octahedral packings, and $d=3$ for the triangular packing), and each Apollonian group contains a congruence subgroup of $\mathrm{SL}(2, \Z)$. Under these hypotheses, Fuchs, Stange and Zhang prove the asymptotic Local-Global Conjecture \cite{FuchsStangeZhang}, generalizing \cite{BourgainKontorovich}. So asymptotically $100\%$ of integers satisfying the modular restrictions appear as curvatures in our packings.

Nevertheless, we prove
\begin{theorem}
The Local-Global Conjecture is false for the octahedral, cubic, square, and triangular packings.
\end{theorem}
Counterexamples to the conjecture are shown in Figure \ref{CounterexampleFig}. 

Sections 3-6 treat the four packing types. In each section, we begin by determining the possible modular restrictions for packings of that type. This uses \cite[Theorem 8.1]{FuchsStangeZhang} to determine the modulus, and further elementary but intricate computations to determine the precise congruences. Then we generalize the first step of \cite{HKRS} by parametrizing all circles tangent to a fixed circle in a packing. We use a particular set of circles, a generalization of the Ford circles in the Apollonian packing case, to concretely describe the orbit of a congruence subgroup. In some cases depending on the modular restrictions, we find a quadratic invariant which is well-defined at a given circle. Next, we generalize the second step of \cite{HKRS}, carrying out reciprocity calculations and checking that there is a path of tangent, coprime circles between any two circles in a packing. Again in some cases, depending on the modular restrictions, we find a quadratic invariant which is well-defined across the full packing. We use the quadratic invariant to find obstructions to the Local-Global Conjecture. Even in cases where a quadratic invariant cannot be defined across the full packing, we find partial obstructions within subsets of circles.

\begin{figure}[h] 
\include{intro_graphics}
\caption{Counterexamples to the Local-Global Conjecture: Octahedral Packing Containing No Curvatures $n^2$, $2n^2$; Cubic Packing Containing No Curvatures $n^2$, $2n^2$; Square Packing Containing No Curvatures $n^2$; Triangular Packing Containing No Curvatures $3n^2$.}
\label{CounterexampleFig}
\end{figure}

There are some important complexities that arise in our work but not in \cite{HKRS}. The first is in the parametrization of circles tangent to a fixed circle. Sarnak's insight is that the circles tangent to a fixed circle in an Apollonian packing are an orbit of a congruence subgroup of $\mathrm{SL}(2, \Z)$. This allows for the parametrization by a shifted binary quadratic form. In each case, we have an analogous congruence subgroup, but in the octahedral, square, and triangular cases the circles tangent to a fixed circle are a union of two orbits. We then obtain a parametrization by two shifted binary quadratic forms. The value of $\res{b}{a}$ for fixed $a$ and $b$ tangent and coprime to $a$ may be the same for all $b$, or it may take one value for all $b$ in one orbit and the opposite value for all $b$ in the other orbit. This prevents us from having a well-defined quadratic invariant under some modular conditions, but in these cases there is still a simple and predictable pattern of Kronecker symbols in the packing.

The second complexity is in the quadratic reciprocity law. We observed above that in \cite{HKRS}, the modular restrictions in Apollonian packings seem to conspire with quadratic reciprocity to make an invariant constructed from Kronecker symbols well-defined. For cubic, square, and triangular packings, some possible modular restrictions conspire with quadratic reciprocity, but some modular restrictions conspire against it. In particular, in the cases when we have two tangent circles of curvature $3$ mod $4$, and no further modular information to distinguish these circles, any definition of a quadratic invariant using Kronecker symbols and modular congruences will not be consistent from one circle to the other. But again in these cases, there is still a simple and predictable pattern of Kronecker symbols in the packing. 

Because of these two complexities, we do not find a well-defined quadratic invariant for every packing we study--only those with certain congruence restrictions. This is similar to the situation with quartic invariants in \cite{HKRS}. But even in the cases without a well-defined quadratic invariant across the packing, we observe interesting reciprocity phenomena. We often express these as partial quadratic obstructions within subsets of the packing. In several cases, we two-color or three-color our packings such that tangent circles have different colors; then a partial obstruction can be expressed as a restriction on curvatures in one color. See Figures \ref{OctPartialFig} and \ref{CubePartialFig} for examples. Even in the most complicated cases, the Kronecker symbols within the packing exhibit a rigidity; if the value of one Kronecker symbol between adjacent, coprime curvatures is known, that determines the values of all the other nonzero Kronecker symbols in the packing. See Figures  \ref{SquarePartialFigure}, \ref{TrianglePartialFigure}, and \ref{TrianglePartialEvenFigure} for examples of this. Although these partial obstructions do not cause the packings to violate the Local-Global Conjecture, they are still natural motivating examples for a broader theory of reciprocity obstructions in group orbits. 

In cases where we do not find quadratic obstructions to the Local-Global Conjecture, we still cannot rigorously rule them out. The scope of possible obstructions has not yet been precisely defined. But in Section \ref{Data}, we collect data to suggest that the obstructions we find may be the only ones. We computed curvatures for several examples of each packing type, including examples satisfying each possible collection of modular restrictions. Some of these examples exhibit the quadratic obstructions we prove in Theorems \ref{OctObstruction}, \ref{CubeObstruction}, \ref{SquareObstruction}, and \ref{TriangleObstruction}. The rest appear to satisfy the Local-Global Conjecture. We did not find any evidence of cubic, quartic, or other higher-order obstructions, but our search was not comprehensive enough to conjecture that these do not exist.

Throughout this project, we were continually impressed by the flexibility of the argument in \cite{HKRS}. Even in our most complicated cases, where it is not possible to define a quadratic invariant or disprove the Local-Global Conjecture, carrying out this argument reveals a rich reciprocity structure among the curvatures. We had not expected to find such a structure in every case. Although this article only analyzes individual examples, we suspect that there is a general framework which would unify all of these examples.

There are also many more examples which should be analyzed. The hexagonal packing is dual to the triangular packing we consider here. We did enough computation to conjecture that the Local-Global Conjecture is false for the hexagonal packing, but we did not fully develop the theory of this packing type. A natural class of packings with number theoretic properties are the superintegral polyhedral packings defined in \cite{KontorovichNakamura}. A forthcoming article of Allcock, Devlin, Felikson, Kontorovich and the third author classifies all superintegral polyhedral packings: they are obtained only from the tetrahedron, cube, octahedron, cubeoctahedron, rhombic dodecahedron, and certain polyhedra constructed from these. We investigated the cubeoctahedral packing but did not reach a conclusion. The issue is that the modular restrictions in this packing type lead to many pairs of tangent circles whose curvatures are both divisible by powers of 2, making reciprocity calculations more difficult. A theory of local $2$-adic or $p$-adic reciprocity obstructions may be necessary to resolve cases like this one. 

It would also be illuminating to examine reciprocity obstructions in higher dimensional packings. Kontorovich proves the Local-Global Conjecture for integral Soddy sphere packings in \cite{Kontorovich}; however it is plausible that reciprocity obstructions might exist in other higher-dimensional sphere packing types. Additionally, we propose investigating reciprocity phenomena in the Apollonian superpacking defined in \cite{GLMWYSuper}, and its generalizations. Because the superpacking contains all primitive integral Apollonian packings, there should not be a well-defined quadratic invariant across the full superpacking. Still, understanding the propagation of Kronecker symbols across tangent circles in the superpacking may give a bigger picture of quadratic invariants. The symmetry group of the superpacking is $\mathrm{SL}(2, \Z[i])$; perhaps lifting to metaplectic covers of this group would help to explain reciprocity phenomena generally. Finally, it is an important open problem to find more examples of reciprocity obstructions in other integral group orbits which have nothing to do with circle packings, to see how much further the theory can extend. 

\subsection{Acknowledgements} We thank Edna Jones, Alex Kontorovich, James Rickards, Kate Stange, Spencer Whitehead, and David Yang for helpful conversations relating to this project. We thank Swarthmore College for funding the summer research project that led to this article, and we thank the Swarthmore Mathematics and Statistics Department for providing us with opportunities and guidance in the field.

\section{Background and Notation}
To any configuration of pairwise tangent or disjoint circles, we can associate a planar tangency graph, which has a vertex for each circle and an edge between vertices if and only if the associated circles are tangent. Our constructions begin with a base configuration of circles $\B$, and a dual configuration $\hat{\B}$. We assume that the circles in $\B$ are pairwise disjoint or tangent, and that the same holds for $\hat{\B}$. We also assume that their tangency graphs are a pair of dual planar graphs, so that vertices in the tangency graph of $\hat{\B}$ correspond to faces in the graph of $\B$. Finally, we assume that a circle in $\B$ intersects a circle in $\hat{\B}$ orthogonally if they correspond to an incident vertex/face pair, and that these circles are disjoint otherwise. These conditions imply that any intersections between base circles and dual circles are orthogonal at shared points of tangency.

The Koebe-Andreev-Thurston theorem in geometry implies that any for any dual pair of three-dimensional spherical polyhedra, their graphs can be realized as tangency graphs for a pair of dual circle configurations. Moreover, this realization is unique up to M\"obius transformations. Kontorovich and Nakamura use this theorem to define \emph{polyhedral packings} \cite{KontorovichNakamura}. The polyhedral Apollonian group is generated by the reflections through the dual circles. The orbit of the base circles under the Apollonian group is an infinite collection of circles, pairwise tangent or disjoint, called a polyhedral packing. From the Koebe-Andreev-Thurston theorem, we can deduce that all packings for a given polyhedron are M\"obius-equivalent. For a quadruple of mutually tangent circles, the tangency graph is the graph of a tetrahedron, and thus the Apollonian packing is the tetrahedral packing. 

The article \cite{RWYY} introduces a further generalization: \emph{packings constructed from tilings of the plane}. In this case, we begin with an infinite base configuration of circles whose tangency graph is the graph of a cellular decomposition of $\R^2$; the dual configuration corresponds to the dual cellular decomposition. Again, the Apollonian group is generated by reflections through the dual circles, and the orbit of the base circles under this group is a fractal circle packing. The Koebe-Andreev Thurston theorem is not known for arbitrary cellular decompositions of the plane, but it is known for triangulations \cite{BeardonStephenson} and for the square tiling \cite{Schramm}. Thus we know that all triangular packings are M\"obius-equivalent, and the same for square packings.

\begin{figure}[h]
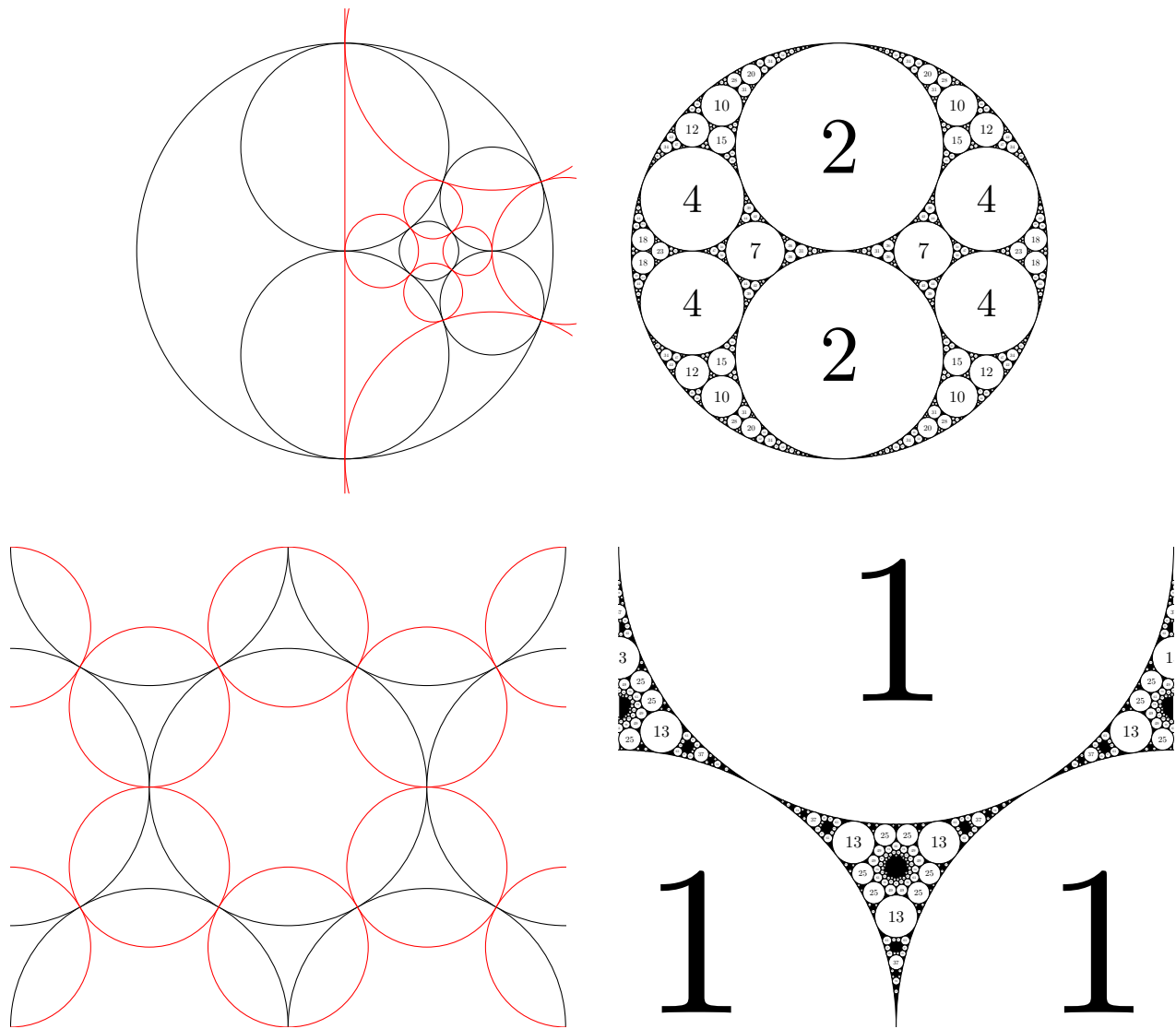
 
\begin{center}
\vspace{-.25in}
\include{octpackingfigure}
\vspace{-.25in}
\include{tripackingfigure}
\end{center}
\caption{Base and Dual Octahedral Circle Configurations; Octahedral Packing; Base and Dual Triangular Circle Configurations; Triangular Packing}
\label{PolyhedralFig}
\end{figure}

The article \cite{ABBRTVWWY} works out some of the algebraic theory of polyhedral circle packings, generalizing Equation \eqref{Descartes}. The following Proposition summarizes the results of \cite[Section 2]{ABBRTVWWY}:
\begin{proposition} \label{GeneralDescartes}
    Given a polyhedron $\Pi$ with $n$ vertices, let $\B$ be a $\Pi$-polyhedral configuration of circles, and let $\hat{\B}$ be the dual configuration. Then the $n$-tuple of curvatures of $\B$ satisfies $n-4$ independent linear equations, and an additional independent homogeneous quadratic equation. Moreover, the Apollonian group acts on the set of solutions to these equations. The generator for a given dual circle fixes the curvatures orthogonal to that dual circle, and transforms the remaining curvatures \cite{ABBRTVWWY}.
\end{proposition}
The particular linear relations and quadratic forms for the octahedral and cubic packings were already worked out in the articles \cite{GuettlerMallows}, \cite{Stange} respectively. In the case of packings constructed from tilings, a general form of Proposition \ref{GeneralDescartes} is not yet proven, but particular linear and quadratic relations for the triangular and square packings we study are shown in \cite{RWYY}. Each of the following four sections of this article begins with a statement of the linear and quadratic relations for the packing type under consideration.

A packing is called \emph{integral} if the curvatures are all integers. Such a packing can be rescaled so that the greatest common divisor of all the curvatures is 1; in this case, the packing is called \emph{primitive}. The superpacking is defined as the orbit of the base circle configuration $\B$ under the supergroup generated by reflections through both $\B$ and $\hat{\B}$; a packing is called \emph{superintegral} if its superpacking is integral. Superintegrality is equivalent to arithmeticity of the supergroup--see \cite{KontorovichNakamura}. Only a small subset of polyhedra and tilings admit integral packings, and an even smaller subset admit superintegral packings. These are natural classes in which to seek quadratic invariants. All the packing types in this article admit superintegral packings. 

All integral packings we study have modular restrictions on the integers that can appear as curvatures. For each packing type, there is an integer $M$ such that the modular restrictions can all be expressed as congruences mod $M$. This is an important consequence of a spectral gap property for the symmetry groups of the packings. Fuchs, Stange, and Zhang prove a general result \cite[Theorem 8.1]{FuchsStangeZhang} which allows $M$ to be computed explicitly in our cases. The ingredients in the calculation are the level of a congruence subgroup of $\mathrm{SL}(2, \Z)$ contained in the symmetry group, and information about the image of the adjoint representation of the symmetry group. For the Apollonian packing, $M=24$. We check that for the octahedral, cubic, square, and triangular packings, the $M$ values are $8$, $4$, $8$, and $12$ respectively. And we go on to determine the possible congruences mod $M$ for each packing type, following \cite{HKRS} in the Apollonian case. See Propositions \ref{OctModular}, \ref{CubeMod}, \ref{SquareMod}, and \ref{TriangleMod}. The \emph{Generalized Local-Global Conjecture} \cite[Conjecture 1.5]{FuchsStangeZhang}, which we disprove, is that every sufficiently large integer subject to the modular restrictions will appear as a curvature in a given packing. 

For certain calculations in this article, it is helpful to introduce an \emph{inversive coordinate system} for oriented generalized circles (i.e. circles and lines) in the complex plane. An oriented generalized circle is described by four coordinates, $(\tilde{b}, b, h_1, h_2) \in \R^4$. Here $b$ is the curvature, $\tilde{b}$ is the cocurvature or curvature after reflecting across the unit circle, and $(h_1, h_2)$ are the coordinates of the center multiplied by the curvature. In the special case of a line, $b$ is zero and $(h_1, h_2)$ are the coordinates of a unit normal vector to the line. The orientation of the circle is determined by the sign of $b$; for lines, it is determined by the direction of $(h_1, h_2)$. All oriented generalized circles satisfy the quadratic equation
\begin{equation}
h_1^2+h_2^2-b \tilde{b} = 1
\end{equation}
The quadratic form on the left side gives rise to a symmetric bilinear form, which we view as an inner product between circles. If two circles intersect, then their inner product is the dot product of their oriented unit normal vectors; if they are disjoint, then their inner product has absolute value greater than one. This inner product can also be used to compute the reflection of one circle across another. A final key property is that M\"obius transformations act linearly on generalized circles in this coordinate system. For further details, see \cite{LagariasMallowsWilks}.

Each packing type we study contains a primitive integral \emph{strip packing} bounded by two parallel lines. We can take the bottom line to be the real axis, and we call the collection of all circles tangent to this axis the \emph{generalized Ford circles}. In the case of the Apollonian packing, the Ford circles have a long history of study \cite{Ford}. They are parametrized by rational numbers $\frac{x}{y}$ in reduced form; each circle is tangent to the real axis at $\frac{x}{y}$ and has curvature $y^2$. The Ford circles for the Apollonian packing give a natural visualization of the action of $\mathrm{SL}(2, \Z)$ on $\Q$. In all cases, the Ford circles are one or more orbits of the Apollonian subgroup generated by the reflections through dual circles orthogonal to the real axis. This subgroup is conjugate to a congruence subgroup of $\mathrm{SL}(2, \Z)$, of level 2 for the octahedral, cubic, and square packings, and level 6 for the triangular packing. Working with the generalized Ford circles helps us concretely parametrize the orbits of these subgroups. Since all packings of a given type are M\"obius equivalent, it suffices to find the parametrization for the Ford circles, and then use linearity in inversive coordinates to extend to other packings. We use the same technique to find chains of circles with coprime curvatures in each packing type: the argument begins by parametrizing circles or chains of circles tangent to the top and bottom line in the strip packing. 

The \emph{Kronecker symbol} is defined as follows. For odd primes $p$, it is
\begin{equation}
\res{a}{p}=\begin{cases} 1 & a \equiv x^2 \not\equiv 0 \bmod p \text{ for some } x\in \Z \\ -1 & a \not\equiv x^2 \bmod p \text{ for all } x \in \Z \\ 0 & a \equiv 0 \bmod p \end{cases}
\end{equation}
It is further defined by 
\begin{equation}
\res{a}{2}=\begin{cases} 1 & a \equiv 1, 7 \bmod 8 \\ -1 & a \equiv 3,5 \bmod 8 \\ 0 & a \equiv 0 \bmod 2 \end{cases} \qquad \qquad \res{a}{-1} = \begin{cases} 1 & a \geq 0 \\ -1 & a < 0 \end{cases}
\end{equation}
Finally, $\res{a}{b}$ is extended to all nonzero integers $b$ by multiplicativity in $b$. Let $a_0$ and $b_0$ denote the odd parts or greatest odd divisors of $a$ and $b$. When at least one of $a$, $b$ is positive, we have the quadratic reciprocity law
\begin{equation}
\res{a}{b} = \res{b}{a} (-1)^{(a_0-1)(b_0-1)/4}
\end{equation}

\section{Octahedral Packing}
The octahedral packing is built from a configuration of six circles with an octahedral tangency graph as shown.

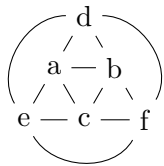
\begin{center}
\begin{tikzpicture}[scale=.4]
    \node (C) at (0,0)    {c};
    \node (A) at (-1,1.732)    {a};
    \node (B) at (1,1.732)    {b};
    \node (D) at (0,3.464)    {d};
    \node (E) at (2,0)    {f};
    \node (F) at (-2,0)    {e};

    \draw (A) -- (B);
    \draw (A) -- (C);
    \draw (B) -- (C);
    \draw (A) -- (D);
    \draw (D) -- (B);
    \draw (A) -- (F);
    \draw (F) -- (C);
    \draw (B) -- (E);
    \draw (C) -- (E);
    \draw (D) to[out=0,in=60] (E);
    \draw (E) to[out=240,in=300] (F);
    \draw (D) to[out=180,in=120] (F);
\end{tikzpicture}
\end{center}

Any sextuple in an octahedral packing satisfies the following linear and quadratic equations, 
\begin{align}
    &a + f = b + e = c + d \label{OctLinear}\\
    &\left(\frac{a+f}{2}\right)^2 - (a+f)(a+b+c) + a^2 + b^2 + c^2 = 0 \label{OctQuadratic}
\end{align}
first given in \cite{GuettlerMallows}. Note that in an integral packing, $a + f$ must be even by Equation \eqref{OctQuadratic}. Define $2w = a + f$. If we assume that the sextuple is oriented in the standard way, with at most one exterior circle of negative curvature, then $w$ is a positive integer.

From these equations, we can write $w$ in terms of $a, b, c$.
\begin{align}
    w = a + b + c \pm \sqrt{2(ab + ac + bc)}.
\end{align}
Then $ab + ac + bc$ must be twice a square. Let $ab + ac + bc = 2m^2$ where $m$ is an integer. Without loss of generality, we can make $m$ non-negative.

We now establish modular restrictions on curvatures in an octahedral packing, following the strategy of \cite[Section 3.1]{HKRS}. Zhang proves that the modulus for congruence restrictions in this packing is 8 \cite[Theorem 1.1]{Zhang}. We will prove the following proposition:

\begin{proposition} \label{OctModular}
    Any primitive integral octahedral sextuple of curvatures is congruent to one of the following modulo 8:
    \begin{center} \, \hfill
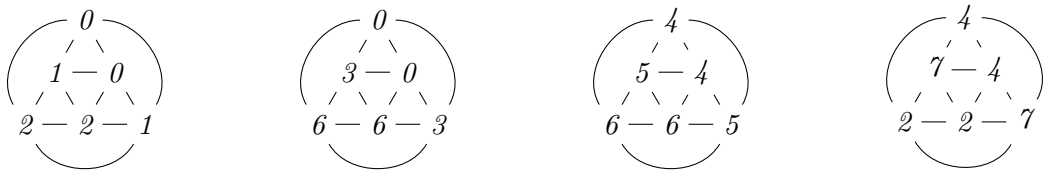
\begin{tikzpicture}[scale=.4]
    \node (C) at (0,0)    {2};
    \node (A) at (-1,1.732)    {1};
    \node (B) at (1,1.732)    {0};
    \node (D) at (0,3.464)    {0};
    \node (E) at (2,0)    {1};
    \node (F) at (-2,0)    {2};

    \draw (A) -- (B);
    \draw (A) -- (C);
    \draw (B) -- (C);
    \draw (A) -- (D);
    \draw (D) -- (B);
    \draw (A) -- (F);
    \draw (F) -- (C);
    \draw (B) -- (E);
    \draw (C) -- (E);
    \draw (D) to[out=0,in=60] (E);
    \draw (E) to[out=240,in=300] (F);
    \draw (D) to[out=180,in=120] (F);
\end{tikzpicture} \hfill
\begin{tikzpicture}[scale=.4]
    \node (C) at (0,0)    {6};
    \node (A) at (-1,1.732)    {3};
    \node (B) at (1,1.732)    {0};
    \node (D) at (0,3.464)    {0};
    \node (E) at (2,0)    {3};
    \node (F) at (-2,0)    {6};

    \draw (A) -- (B);
    \draw (A) -- (C);
    \draw (B) -- (C);
    \draw (A) -- (D);
    \draw (D) -- (B);
    \draw (A) -- (F);
    \draw (F) -- (C);
    \draw (B) -- (E);
    \draw (C) -- (E);
    \draw (D) to[out=0,in=60] (E);
    \draw (E) to[out=240,in=300] (F);
    \draw (D) to[out=180,in=120] (F);
\end{tikzpicture} \hfill
\begin{tikzpicture}[scale=.4]
    \node (C) at (0,0)    {6};
    \node (A) at (-1,1.732)    {5};
    \node (B) at (1,1.732)    {4};
    \node (D) at (0,3.464)    {4};
    \node (E) at (2,0)    {5};
    \node (F) at (-2,0)    {6};

    \draw (A) -- (B);
    \draw (A) -- (C);
    \draw (B) -- (C);
    \draw (A) -- (D);
    \draw (D) -- (B);
    \draw (A) -- (F);
    \draw (F) -- (C);
    \draw (B) -- (E);
    \draw (C) -- (E);
    \draw (D) to[out=0,in=60] (E);
    \draw (E) to[out=240,in=300] (F);
    \draw (D) to[out=180,in=120] (F);
\end{tikzpicture} \hfill
\begin{tikzpicture}[scale=.4]
    \node (C) at (0,0)    {2};
    \node (A) at (-1,1.732)    {7};
    \node (B) at (1,1.732)    {4};
    \node (D) at (0,3.464)    {4};
    \node (E) at (2,0)    {7};
    \node (F) at (-2,0)    {2};

    \draw (A) -- (B);
    \draw (A) -- (C);
    \draw (B) -- (C);
    \draw (A) -- (D);
    \draw (D) -- (B);
    \draw (A) -- (F);
    \draw (F) -- (C);
    \draw (B) -- (E);
    \draw (C) -- (E);
    \draw (D) to[out=0,in=60] (E);
    \draw (E) to[out=240,in=300] (F);
    \draw (D) to[out=180,in=120] (F);
\end{tikzpicture} \hfill \,
\end{center}
\end{proposition}

The proof is via a sequence of lemmas. 

\begin{lemma} \label{OctModLemma1}
    Let $(a, b, c, d, e, f)$ be a primitive integral sextuple. Then $w$ must be odd. 
\end{lemma}

\begin{proof}
    Rewriting Equation \eqref{OctQuadratic}, we get 
    \begin{align*}
        w^2 - 2w(a+b+c) + a^2+b^2+c^2 = 0.
    \end{align*}
    If $w$ is even, $w^2 - 2w(a+b+c) \equiv 0 \bmod 4$, so necessarily $a^2+b^2+c^2 \equiv 0 \bmod 4$. A square of any integer is congruent to either $0 \bmod 4$ when the integer is even or $1 \bmod 4$ when the integer is odd. Therefore, $a, b, c$ must all be even, which implies that $d, e, f$  are even as well. Then the sextuple is not primitive. 
\end{proof}

Because $w$ is odd, $2w \equiv 2 \; \text{or} \; 6 \bmod 8$.

\begin{lemma} \label{OctModLemma2}
    Any primitive integral sextuple contains exactly two odd curvatures and four even curvatures.
\end{lemma}
\begin{proof}
    Each diagonal pair $(a,f), (b,e), (c,d)$ have the same parity, so there are an even number of odd curvatures in the sextuple. We will show that there is one diagonal pair of odd curvatures and two diagonal pairs of even curvatures. If all the curvatures are even, then the packing is not primitive. Next, suppose there are exactly four odd curvatures. Choosing any three tangent circles, we must have two odd curvatures and one even curvature. Then, $a+b+c$ is even and $\sqrt{2(ab+ac+bc)}$ is also even, so $w$ is even, which leads to a contradiction. Finally, if all the curvatures are odd, then $ab + ac + bc$ is odd and cannot be $2m^2$. 
\end{proof}

Within a sextuple, we have one diagonal pair of odd curvatures. For the other two diagonal pairs of even curvatures, the residues modulo 8 must as follows:
\begin{itemize}
    \item If $2w \equiv 2 \bmod 8: (0, 2)$ or $(4, 6)$;
    \item if $2w \equiv 6 \bmod 8: (0, 6)$ or $(2, 4)$.
\end{itemize}

\begin{lemma} \label{OctModLemma3}
    Let $(a, b, c, d, e, f)$ be a primitive integral sextuple. Then the residues of any pair of tangent circles modulo 8 cannot be $(2,6)$.
\end{lemma}
\begin{proof}
    Suppose the opposite, and choose $b, c$ as the tangent circles. Then $a$ is a circle tangent to both $b$ and $c$, and $a$ is necessarily odd. Then
    \begin{align*}
        ab+ac+bc \equiv a(b+c) + bc \equiv 0+4 \equiv 4 \bmod 8.
    \end{align*}
    This means that $ab+ac+bc$ has exactly 2 factors of 2, then $ab+ac+bc$ cannot be $2m^2$ for some integer $m$.
\end{proof}

It follows from Lemma \ref{OctModLemma3} that the set of residues modulo 8 of the four even curvatures in a sextuple must be:
\begin{itemize}
    \item If $2w \equiv 2 \bmod 8: (0,0,2,2)$ or $(4,4,6,6)$;
    \item if $2w \equiv 6 \bmod 8: (0,0,6,6)$ or $(2,2,4,4)$.
\end{itemize}

\begin{lemma} \label{OctModLemma4}
    Let $(a, b, c, d, e ,f)$ be a primitive integral sextuple. If $2w \equiv 2 \bmod 8$, then the odd curvatures are both $1$ or both $5$ mod $8$. If $2w \equiv 6 \bmod 8$, the the odd curvatures are both $3$ or both $7$ mod $8$.
\end{lemma}
\begin{proof}
    Choose $a$ as odd and $b, c$ such that $b \not\equiv c \bmod 8$. Because the residues modulo 8 of the even diagonals are the same, $b + c \equiv b + e \equiv 2w \bmod 8$. Then since $bc \equiv 0 \bmod 8$, $ab+ac+bc\equiv a(b+c) \equiv 2aw \bmod 8$ and $v_2(ab+ac+bc) = 1$. Thus $ab + ac + bc = 2m^2$ with $m$ odd, which implies $ab + ac + bc \equiv 2aw \equiv 2 \bmod 8$. When $2w \equiv 2 \bmod 8$, $a \equiv 1 \; \text{or} \; 5 \bmod 8$; when $2w \equiv 6 \bmod 8$, $a \equiv 3 \; \text{or} \; 7 \bmod 8$. Because $a+f \equiv 2w \bmod 8$, in all cases $a \equiv f \bmod 8$. 
\end{proof}

The following lemma is the analogue of Proposition 3.1 in \cite{HKRS}. 

\begin{lemma} \label{OctModLemma5}
    Let $(a, b, c, d, e ,f)$ be a primitive integral sextuple. The sum of the curvatures of any two tangent circles is not congruent to $5,7 \bmod 8$.
\end{lemma}
\begin{proof}
    Suppose that $a, b, c$ are curvatures of mutually tangent circles with $a$ odd and $b, c$ even. We know $ab+ac+bc = 2m^2$ where $m$ is a positive integer. Then
    \begin{align*}
        (a+b)(a+c) = a^2 + 2m^2.
    \end{align*}
    Consider a prime number $p > 2$ such that $p \nmid a, m$. Then if $p \mid a^2 + 2m^2$, $-2$ is a square modulo $p$, so $p \equiv 1 \; \text{or} \; 3 \bmod 8$. Contrapositively, if an odd prime $p \equiv 5 \; \text{or} \; 7 \bmod 8$ is factor of $a^2 + 2m^2$, then $p \mid a, \, m$ and $v_{p}(a^2 + 2m^2)$ is even. 
    
    Assume $a+b \equiv 5 \; \text{or} \; 7 \bmod 8$. Then $a+b$ must have a prime factor $p \equiv 5 \; \text{or} \; 7 \bmod 8$ with $v_{p}(a+b)$ odd. Because $p$ is also a factor of $a^2 + 2m^2$, $p\mid a, \, m$ and $v_{p}(a^2 + 2m^2)$ is even. So $v_{p}(a+c)$ is odd. Then $p$ is a common factor of $a+b$, $a+c$, and $a$, so the packing is not primitive. 
\end{proof}

Note that this proof assumes $a+b$ is a positive integer, a consequence of how we orient our circles.

Given all the restrictions proven in the lemmas above, any primitive integral sextuple must be congruent to $(1, 0, 2, 0, 2, 1)$, $(3, 0, 6, 0, 6, 3)$, $(5, 4, 6, 4, 6, 5)$, or $(7, 2, 4, 2, 4, 7)$ modulo 8, proving Proposition \ref{OctModular}. 

Now we need to generalize this modular restriction to the entire packing. 

\begin{proposition} \label{OctGroupMod}
    In a primitive integral octahedral packing, every sextuple of circles belongs to the same type modulo 8. 
\end{proposition}
\begin{proof}
    It suffices to show that the octahedral Apollonian group fixes each congruence type. Each generator of the Apollonian group fixes the curvatures of three tangent circles arranged in a triangle. We can see from the cases in Proposition \ref{OctModular} that these three curvatures uniquely determine the sextuple mod 8, so the Apollonian group cannot map between different sextuples mod 8.
\end{proof}

In particular, any primitive integral octahedral packing has all curvatures in exactly three congruence classes mod 8. Section \ref{Data} below contains examples of all four packing types mod 8. 

The next step is to parametrize all circles tangent to a fixed circle in an octahedral packing. It is shown in \cite{Zhang} that the stabilizer of a fixed circle in the octahedral Apollonian group is conjugate to a congruence subgroup of $\mathrm{SL}(2, \Z)$. Our approach is similar; we first parametrize the octahedral Ford circles, or circles tangent to the real axis in the octahedral strip packing. We then can apply a M\"{o}bius transformation to parametrize the circles tangent to any fixed circle. Another description of the orbits of this congruence subgroup is given by the Rosen continued fractions; see \cite{MR65632, MR3609844}.

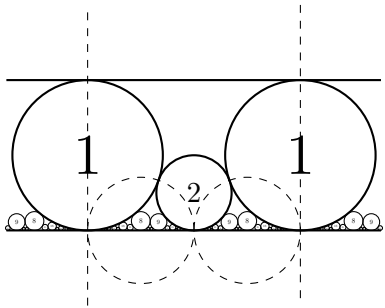
\begin{figure}[h] \label{OctFordFig}
\begin{center}
\begin{tikzpicture}
\draw[black, dashed] (-0.0,-2.0) -- (0.0,2.0);
\draw[black, dashed] (2.828425,2.0) -- (2.828425,-2.0);
\draw[black, dashed] (0.7071085623775818,-1.0) circle (0.7071085623775818);
\draw[black, dashed] (2.1213256871327455,-1.0) circle (0.7071085623775818);
\draw[black, thick] (0.0,0.0) circle (1.0)node[scale=2.0]{1};
\draw[black, thin] (-0.70710675,-0.875) circle (0.125)node[scale=0.25]{8};
\draw[black, thin] (0.70710675,-0.875) circle (0.125)node[scale=0.25]{8};
\draw[black, thin] (3.535533875,-0.875) circle (0.125)node[scale=0.25]{8};
\draw[black, thin] (2.3570226111111108,-0.9444444444444444) circle (0.05555555555555555)node[scale=0.1111111111111111]{18};
\draw[black, thin] (1.21218306122449,-0.9795918367346939) circle (0.02040816326530612)node[scale=0.04081632653061224]{49};
\draw[black, thin] (1.7677669375,-0.96875) circle (0.03125)node[scale=0.0625]{32};
\draw[black, thin] (3.771236111111111,-0.8888888888888888) circle (0.1111111111111111)node[scale=0.2222222222222222]{9};
\draw[black, thin] (-0.942809,-0.8888888888888888) circle (0.1111111111111111)node[scale=0.2222222222222222]{9};
\draw[black, thin] (0.942809,-0.8888888888888888) circle (0.1111111111111111)node[scale=0.2222222222222222]{9};
\draw[black, thin] (2.26274168,-0.96) circle (0.04)node[scale=0.08]{25};
\draw[black, thin] (1.1785113055555554,-0.9861111111111112) circle (0.013888888888888888)node[scale=0.027777777777777776]{72};
\draw[black, thin] (1.61624406122449,-0.9795918367346939) circle (0.02040816326530612)node[scale=0.04081632653061224]{49};
\draw[black, thick] (3.9,-1.0) -- (-1.08,-1.0);
\draw[black, thick] (3.9,1.0) -- (-1.08,1.0);
\draw[black, thin] (-0.4714045,-0.9444444444444444) circle (0.05555555555555555)node[scale=0.1111111111111111]{18};
\draw[black, thin] (0.4714045,-0.9444444444444444) circle (0.05555555555555555)node[scale=0.1111111111111111]{18};
\draw[black, thin] (3.299831666666667,-0.9444444444444444) circle (0.05555555555555555)node[scale=0.1111111111111111]{18};
\draw[black, thin] (2.47487371875,-0.96875) circle (0.03125)node[scale=0.0625]{32};
\draw[black, thin] (1.9798989800000002,-0.98) circle (0.02)node[scale=0.04]{50};
\draw[black, thin] (2.020305081632653,-0.9795918367346939) circle (0.02040816326530612)node[scale=0.04081632653061224]{49};
\draw[black, thin] (1.6499158194444443,-0.9861111111111112) circle (0.013888888888888888)node[scale=0.027777777777777776]{72};
\draw[black, thin] (-0.56568544,-0.96) circle (0.04)node[scale=0.08]{25};
\draw[black, thin] (0.56568544,-0.96) circle (0.04)node[scale=0.08]{25};
\draw[black, thin] (3.39411256,-0.96) circle (0.04)node[scale=0.08]{25};
\draw[black, thin] (2.4243661020408163,-0.9795918367346939) circle (0.02040816326530612)node[scale=0.04081632653061224]{49};
\draw[black, thin] (-0.353553375,-0.96875) circle (0.03125)node[scale=0.0625]{32};
\draw[black, thin] (0.353553375,-0.96875) circle (0.03125)node[scale=0.0625]{32};
\draw[black, thin] (3.1819805,-0.96875) circle (0.03125)node[scale=0.0625]{32};
\draw[black, thin] (2.54558442,-0.98) circle (0.02)node[scale=0.04]{50};
\draw[black, thin] (1.885618111111111,-0.8888888888888888) circle (0.1111111111111111)node[scale=0.2222222222222222]{9};
\draw[black, thick] (1.4142135,-0.5) circle (0.5)node[scale=1.0]{2};
\draw[black, thin] (2.121320375,-0.875) circle (0.125)node[scale=0.25]{8};
\draw[black, thin] (3.88908728125,-0.96875) circle (0.03125)node[scale=0.0625]{32};
\draw[black, thin] (-0.28284272,-0.98) circle (0.02)node[scale=0.04]{50};
\draw[black, thin] (0.28284272,-0.98) circle (0.02)node[scale=0.04]{50};
\draw[black, thin] (3.11126984,-0.98) circle (0.02)node[scale=0.04]{50};
\draw[black, thin] (1.13137084,-0.96) circle (0.04)node[scale=0.08]{25};
\draw[black, thin] (-0.4040610204081633,-0.9795918367346939) circle (0.02040816326530612)node[scale=0.04081632653061224]{49};
\draw[black, thin] (0.4040610204081633,-0.9795918367346939) circle (0.02040816326530612)node[scale=0.04081632653061224]{49};
\draw[black, thin] (2.592724861111111,-0.9861111111111112) circle (0.013888888888888888)node[scale=0.027777777777777776]{72};
\draw[black, thin] (3.232488142857143,-0.9795918367346939) circle (0.02040816326530612)node[scale=0.04081632653061224]{49};
\draw[black, thin] (-1.06066015625,-0.96875) circle (0.03125)node[scale=0.0625]{32};
\draw[black, thin] (1.06066015625,-0.96875) circle (0.03125)node[scale=0.0625]{32};
\draw[black, thin] (3.636549163265306,-0.9795918367346939) circle (0.02040816326530612)node[scale=0.04081632653061224]{49};
\draw[black, thick] (2.828427,0.0) circle (1.0)node[scale=2.0]{1};
\draw[black, thin] (3.6769552599999997,-0.98) circle (0.02)node[scale=0.04]{50};
\draw[black, thin] (-0.23570226388888887,-0.9861111111111112) circle (0.013888888888888888)node[scale=0.027777777777777776]{72};
\draw[black, thin] (0.23570226388888887,-0.9861111111111112) circle (0.013888888888888888)node[scale=0.027777777777777776]{72};
\draw[black, thin] (3.0641293888888885,-0.9861111111111112) circle (0.013888888888888888)node[scale=0.027777777777777776]{72};
\draw[black, thin] (-0.8081220408163265,-0.9795918367346939) circle (0.02040816326530612)node[scale=0.04081632653061224]{49};
\draw[black, thin] (0.8081220408163265,-0.9795918367346939) circle (0.02040816326530612)node[scale=0.04081632653061224]{49};
\draw[black, thin] (-0.84852814,-0.98) circle (0.02)node[scale=0.04]{50};
\draw[black, thin] (0.84852814,-0.98) circle (0.02)node[scale=0.04]{50};
\draw[black, thin] (1.69705628,-0.96) circle (0.04)node[scale=0.08]{25};
\end{tikzpicture}
\caption{Octahedral Ford Circles and Duals}
\end{center}
\end{figure}

The octahedral strip packing has base circles of curvatures $(1, 0, 2, 0, 2, 1)$. We place it so that the bottom line is the real axis, and a circle of curvature $1$ is centered on the imaginary axis. There are four base circles tangent to the real line and four dual circles orthogonal to the real line. The generalized Ford circles are all circles in the packing tangent to the real line. The base, dual, and Ford circles are shown in Figure \ref{OctFordFig}.

We will show that each Ford circle intersects the real axis at a point $\frac{x}{y} \sqrt{2}$ for $\frac{x}{y} \in \Q$, and that the Ford circles are in bijection with these points. More specifically, we find the following parametrization. 

\begin{proposition} \label{OctFord}
    The generalized Ford circles are parametrized by $x, y \in \Z$ with $y \geq 0$ and $\gcd(x,y)=1$. Their inversive coordinates are as follows:
    \begin{align*}
    & c_\alpha(x,y)=(4x^2,2y^2,2\sqrt{2}xy,1)\hspace{3pt} \text{ if } x \text{ is odd} \\
    & c_\beta(x,y)=(2x^2,y^2, \sqrt{2}xy,1) \hspace{13pt} \text{ if } x \text{ is even} 
    \end{align*}
    Each circle is tangent to the real axis is at $\frac{x}{y} \sqrt{2}$. 
\end{proposition}

We will prove this proposition in full detail. The analogous propositions in later sections have similar proofs, which we will omit.

\begin{proof} Let $S_{\alpha} = \lbrace c_{\alpha}(x,y) \mid  \gcd(x,y)=1, \, x \text{ odd} \rbrace$ and $S_{\beta} = \lbrace c_{\beta}(x,y) \mid \gcd(x,y)=1, \, x \text{ even} \rbrace$. We will show that these sets have the following properties. 
\begin{enumerate}
    \item The base circles tangent to $\hat{\R}$ are in $S_{\alpha} \cup S_{\beta}$. 
    \item The sets $S_{\alpha}$ and $S_{\beta}$ are preserved by reflection through the dual circles orthogonal to $\hat{\R}$.
    \item The only circles in $S_{\alpha} \cup S_{\beta}$ that intersect one of the dual circles orthogonal to $\hat{\R}$ are the base circles tangent to $\hat{\R}$.
\end{enumerate}

The reason that these properties imply the proposition is as follows. All Ford circles are obtained from the base circles tangent to $\hat{\R}$ under the group generated by reflections through the dual circles orthogonal to $\hat{\R}$. Thus properties (1) and (2) immediately imply that every generalized Ford circle is in $S_{\alpha} \cup S_{\beta}$. To show the opposite inclusion, note that the dual circles orthogonal to $\hat{\R}$ generate a geometrically finite hyperbolic reflection group acting on the upper half-plane, whose fundamental domain is the region outside the circles. Given any circle in $S_{\alpha} \cup S_{\beta}$, we can map it to a new circle that intersects the fundamental domain. Since this new circle is still tangent to $\hat{\R}$, it must intersect some dual circle orthogonal to $\hat{\R}$. Then by properties (2) and (3), the new circle must be one of the base circles, so the original circle must be a generalized Ford circle. 

To check Property 1, note that $c_{\alpha}(1,0)=(4, 0, 0, 1)$, $c_{\beta}(0,1) = (0, 1, 0, 1)$, $c_{\beta}(2,1) = (8, 1, 2\sqrt{2}, 1)$, and $c_{\alpha}(1,1)=(4, 2, 2\sqrt{2}, 1)$. These are the four base circles tangent to $\hat{\R}$.

The dual circles orthogonal to $\hat{\R}$ are $d_1=(0,0,-1,0)$, $d_2=(4\sqrt{2},0,1,0)$, $d_3=(0,\sqrt{2},1,0)$, and $d_4=(4\sqrt{2},\sqrt{2},3,0)$. To check Property (2), we compute the reflections of each $\alpha$-circle and $\beta$-circle across each dual circle, in the following table:
\begin{center}
\begin{tabular}{|c|c|c|} 
\hline
\, & $c_{\alpha}(x,y)$ & $c_{\beta}(x,y)$ \\
 \hline
 $d_1$ & $c_{\alpha}(-x,y)$ & $c_{\beta}(-x,y)$ \\ 
 \hline
 $d_2$ & $c_{\alpha}(2y-x,y)$ & $c_{\beta}(2y-x,y)$ \\ 
 \hline
 $d_3$ & $c_{\alpha}(x,2x-y)$ & $c_{\beta}(x,2x-y)$ \\
 \hline
 $d_4$ & $c_{\alpha}(4y-3x,3y-2x)$ & $c_{\beta}(4y-3x,3y-2x)$\\
 \hline
\end{tabular}
\end{center}
We see that the sets $S_{\alpha}$, $S_{\beta}$ are preserved by each reflection.

To check Property (3), we compute the inner product of each $\alpha$-circle and $\beta$-circle with each dual circle, in the following table:

\begin{center}
\begin{tabular}{|c|c|c|} 
\hline
\, & $c_{\alpha}(x,y)$ & $c_{\beta}(x,y)$ \\
 \hline
 $d_1$ & $-2\sqrt{2} xy$ & $-\sqrt{2} xy$ \\ 
 \hline
 $d_2$ & $2\sqrt{2} (x-2y)y$ & $\sqrt{2} (x-2y)y$ \\ 
 \hline
 $d_3$ & $2\sqrt{2} x(y-x)$ & $\sqrt{2} x(y-x)$ \\
 \hline
 $d_4$ & $2\sqrt{2}(y-x)(x-2y)$ & $\sqrt{2}(y-x)(x-2y)$\\
 \hline
\end{tabular}
\end{center}

We observe that in each case, the inner product is an integer times an irrational number greater than 1. Therefore, the only way for the inner product to have absolute value $\leq 1$ is if it is 0. This occurs only when $(x, y)=(0,1)$, $(1, 0)$, $(1, 1)$, or $(2,1)$. All of these circles are in the base configuration. 
\end{proof}

Next we parametrize the curvatures of all circles tangent to a fixed circle in an arbitrary octahedral packing.  

\begin{proposition} \label{OctParametrization}
Suppose an octahedral packing contains a sextuple with curvatures $(a,b,c,d,e,f)$. Then all circles tangent to the circle of curvature $a$ are parametrized by $x, y \in \Z$ with $y \geq 0$, $\gcd(x,y)=1$. Their curvatures are as follows:
\begin{align*}
& Q_\alpha(x,y)-a=\left(a+c\right)x^2 +  \left(-2a-2b-c+d\right)xy + 2(a+b)y^2 -a\hspace{16pt} \text{ if } x \text{ is odd} \\
& Q_\beta(x,y)-a=\left(\frac{a+c}{2}\right)x^2 +  \left(\frac{-2a-2b-c+d}{2}\right)xy + (a+b)y^2 -a \hspace{3pt} \text{ if } x \text{ is even} 
\end{align*}
\end{proposition}

\begin{proof}
We can apply a M\"{o}bius transformation which maps the four circles $(0,0,0,-1)$, $(0, 1, 0, 1)$, $(4, 2, 2\sqrt{2}, 1)$, $(4, 0, 0, 1)$ to circles of curvature $a$, $b$, $c$, $d$ respectively. This transformation is linear in the inversive coordinate system, so when applied to $c_{\alpha}(x,y)$ and $c_{\beta}(x,y)$, it gives the stated formulas.
\end{proof}

Using Equations \eqref{OctLinear} and \eqref{OctQuadratic}, we find that the discriminant of the quadratic form $Q_{\alpha}(x,y)$ is $-8a^2$ and the discriminant of $Q_{\beta}(x,y)$ is $-2a^2$.

We now have the information necessary to define a quadratic invariants for octahedral packings. We find an invariant $\chi_2$ for packings of type $(0, 1, 2)$ or $(2, 4, 7)$ modulo 8. We follow the approach of \cite{HKRS}, first defining $\chi_2$ at a circle, and then showing that it is the same for all circles in the packing. 

Fix a circle $\mathcal{C}$ of nonzero curvature $a$ in a primitive integral sextuple $(a, b, c, d, e, f)$. The quadratic form $Q_{\alpha}(x,y)$ represents the values $a+c$, $a+d$, and $2(a+b)$. Since $\gcd(a, b, c, d)=1$, for each odd prime $p$ dividing $a$, $Q_{\alpha}(x,y)$ represents a value not divisible by $p$. Therefore by the Chinese remainder theorem, one can find $x$, $y$ such that $Q_{\alpha}(x,y)$ has no prime factors in common with $a$ except possibly $2$. Moreover, it is possible to choose $x, y$ with $y \geq 0$, $\gcd(x,y)=1$, and $x$ odd. Similarly, one can find $x, y$ with $y \geq 0$, $\gcd(x,y)=1$, and $x$ even such that $Q_{\beta}(x,y)$ has no prime factors in common with $a$ except possibly $2$. We can also avoid common factors of 2. If $a$ is even, then by Proposition \ref{OctModular}, one of $b$ and $c$ is even and one is odd. Then either all the $\alpha$-circles have even curvature and all the $\beta$-circles have odd curvature, or vice versa. Thus either $Q_{\alpha}$ or $Q_{\beta}$ represents odd integers. The upshot is that it is always possible to choose a value $\rho$ coprime to $a$ which is represented by $Q_{\alpha}(x,y)$ for $y \geq 0$, $\gcd(x,y)=1$, and $x$ odd, or by $Q_{\beta}(x,y)$ for $y \geq 0$, $\gcd(x,y)=1$, and $x$ even. 

Next, following \cite[Propositions 4.1, 4.2]{HKRS}, we study the Kronecker symbols $\res{\rho}{a}$ where $\rho$ is an integer represented by $Q_\alpha$ or $Q_\beta$.

\begin{lemma} \label{QuadFormKronecker}
Suppose that $Q(x,y)=Ax^2+Bxy+Cy^2$ is a quadratic form which represents some integer coprime to $a$. Suppose further that $2^k a\mid B^2-4AC$, where $k=4$ if $a$ has exactly one factor of $2$, $k=2$ if $a$ has exactly 3 factors of $2$, and $k=0$ otherwise. Then $\res{\rho}{a}$ has the same value for all integers $\rho$ represented by $Q(x,y)$ with $\gcd(\rho, a)=1$. 
\end{lemma}
\begin{proof}
By replacing $Q(x,y)$ with an $\mathrm{SL}(2, \Z)$-equivalent form if necessary, we may assume that $\gcd(A, a)=1$. If $a$ is odd, then $2A$ is invertible mod $a$ and
$$Q(x,y) = Ax^2+Bxy+Cy^2 \equiv A\left(x-\frac{B}{2A}y\right)^2 \bmod a$$
Thus $\res{\rho}{a}=\res{A}{a}$ for every integer $\rho$ represented by $Q(x,y)$ with $\gcd(\rho, a)=1$. 

If $a$ is even, write $a=2^\ell a'$ with $a'$ odd. If $\ell$ is even then $\res{\rho}{a}=\res{\rho}{a'}$ and the above argument holds. If $\ell$ is odd then $\res{\rho}{a}=\res{\rho}{2}\res{\rho}{a'}$, and it suffices to show that $\res{\rho}{2}$ is constant. In this case, since $2^ka\mid B^2-4AC$, we have $32\mid B^2-4AC$. Thus $B$ is even and we may write $B=2B'$ and $8\mid B'^2-AC$. Then 
$$Q(x,y) = Ax^2+Bxy+Cy^2 \equiv A\left(x-\frac{B'}{A}y\right)^2 \bmod 8$$
Thus $\res{\rho}{2}=\res{A}{2}$ for every odd integer $\rho$ represented by $Q(x,y)$.
\end{proof}
The example of $x^2+2xy-3y^2$, which has discriminant 16 and represents integers $\rho$ such that $\res{\rho}{2}=1$ and $\res{\rho}{2}=-1$, shows that the choice of $k$ in the lemma is as small as possible.

\begin{proposition} \label{OctNode}
Suppose $(a,b,c,d,e,f)$ that is a primitive integral octahedral sextuple. Let $\rho_1$, $\rho_2$ be two integers coprime to $a$ which are represented by either $Q_{\alpha}$ or $Q_{\beta}$. Let $a'=a/2$ if $a \equiv 2 \bmod 4$, and $a'=a$ otherwise. Then if $a \not\equiv 3,5 \bmod 8$,
\begin{equation*}
\res{\rho_1}{a'} = \res{\rho_2}{a'}
\end{equation*}
If $a \equiv 3, 5 \bmod 8$, then the equality holds if and only if $\rho_1$, $\rho_2$ are represented by the same quadratic form; if $\rho_1$ is represented by $Q_\alpha$ and $\rho_2$ is represented by $Q_\beta$, then 
\begin{equation*}
\res{\rho_1}{a} = -\res{\rho_2}{a}
\end{equation*}
\end{proposition}
\begin{proof}
If $\rho_1$, $\rho_2$ are represented by the same quadratic form, then the result follows from Lemma \ref{QuadFormKronecker}. The discriminants of $Q_\alpha$ and $Q_\beta$, $-8a^2$ and $-2a^2$ respectively, are always divisible by $a'$ and by a sufficient power of $2$. If $a+c$ is odd, then we may replace the quadratic form $Q_{\beta}(x,y)$ with $Q_{\beta}(2x,y)$ to represent only integer values, and the discriminant becomes $8a^2$.

Finally, suppose that $\rho_1$ is represented by $Q_{\alpha}$ and $\rho_2$ is represented by $Q_{\beta}$. This implies that $a$ is odd; if $a$ is even, then either $Q_{\alpha}$ or $Q_{\beta}$ represents only even integers. Thus $a'=a$. If $\rho_2=Q_{\beta}(x,y)$, $2\rho_2=Q_{\alpha}(x,y)$ and $2\rho_2$ is still coprime to $a$. Then we have
\begin{equation*}
\res{\rho_1}{a}=\res{2\rho_2}{a} = \res{2}{a}\res{\rho_2}{a}
\end{equation*}
which completes the proof.
\end{proof}

For a circle $\C$ of curvature $a$ in a primitive integral octahedral packing of type $(0,1,2)$ or $(2, 4, 7)$ mod 8, we set
\begin{equation} \label{OctInvariant}
\chi_2(\mathcal{C})=\left\lbrace \begin{array}{cc} \res{\rho}{a} & a \not \equiv 2, 4 \bmod 8\\  \res{\rho}{a/2} & a \equiv 2 \bmod 8 \\  \res{-\rho}{a} & a \equiv 4 \bmod 8 \end{array} \right.
\end{equation}
where $\rho$ is any integer coprime to $a$ which is represented by $Q_{\alpha}(x,y)$ for $y \geq 0$, $\gcd(x,y)=1$, and $x$ odd, or by $Q_{\beta}(x,y)$ for $y \geq 0$, $\gcd(x,y)=1$, and $x$ even. Proposition \ref{OctNode} ensures that this expression is well-defined.

Now we study the propagation of $\chi_2$ from circle to circle.

\begin{proposition} \label{OctEdge}
    Suppose that $\mathcal{C}_a$, $\mathcal{C}_b$ are tangent circles with coprime curvatures $a$, $b$ in an integral octahedral packing of type $(0,1,2)$ or $(2,4,7)$ mod $8$. Then $\chi_2(\mathcal{C}_a)=\chi_2(\mathcal{C}_b)$.
\end{proposition}

\begin{proof} We will use $\rho=a+b$ to compute $\chi_2(\mathcal{C}_a)$ and $\chi_2(\mathcal{C}_b)$. Assume without loss of generality that $a$ is odd. Then $a=a'$ and since $a \equiv \pm 1 \bmod 8$, $\res{a}{2}=1$ and $\res{a}{b}=\res{a}{b'}$. If $a \equiv 1 \bmod 8$ or $b \equiv 2 \bmod 8$, then quadratic reciprocity for the Kronecker symbol yields $\res{b}{a}=\res{a}{b}$. Therefore
\begin{equation*}
\chi_2(\mathcal{C}_a) = \res{a+b}{a} = \res{b}{a} = \res{a}{b} = \res{a}{b'} =  \res{a+b}{b'} = \chi_2(\mathcal{C}_b)
\end{equation*}
The only remaining case is $a \equiv 7 \bmod 8$, $b \equiv 4 \bmod 8$. Then quadratic reciprocity for the Kronecker symbol yields $\res{b}{a}=\res{-a}{b}$. Therefore 
\begin{equation*}
\chi_2(\mathcal{C}_a) = \res{a+b}{a} = \res{b}{a} = \res{-a}{b} = \res{-a}{b'} =  \res{-a-b}{b'} = \chi_2(\mathcal{C}_b)
\end{equation*}
\end{proof}

Now we show that $\chi_2(\mathcal{C})$ is identical for all circles $\mathcal{C}$ in the packing. Given two circles, we will find a path between them where each step is between tangent circles of coprime curvatures, allowing us to apply Proposition \ref{OctEdge}. 

\begin{lemma} \label{OctSimultaneous}
All octahedral Ford circles tangent to both the real axis and the circle $(0,1,0,1)$ have coordinates $(4,2n^2,2\sqrt{2}n,1)$ for $n \in \Z$. If an octahedral packing contains a sextuple $(a,b,c,d,e,f)$, then all circles tangent to both $\mathcal{C}_a$ and $\mathcal{C}_b$ have curvatures $c+(-2a-2b-c+d)n+2(a+b)n^2$ for $n\in \Z$.  
\end{lemma}
\begin{proof}
    Based on the parametrization in Proposition \ref{OctFord}, the only circles which have inner product $-1$ with $(0,1,0,1)$ are $c_{\alpha}(1, n) = (4,2n^2,2\sqrt{2}n,1)$. Applying a M\"{o}bius transformation as in the proof of Proposition \ref{OctParametrization} gives curvatures $Q_{\alpha}(1, n)-a=c+(-2a-2b-c+d)n+2(a+b)n^2$.
\end{proof} 

\begin{lemma} \label{OctInsert1}
    For any two tangent circles $\mathcal{C}_a$ and $\mathcal{C}_b$ where $a$ and $b$ are even in a primitive octahedral circle packing, there exists a circle $\mathcal{C}'$ tangent to both $\mathcal{C}_a$ and $\mathcal{C}_b$ whose curvature is coprime to both $a$ and $b$.
\end{lemma}
\begin{proof} We shall show that $Q_{\alpha}(1, n)-a$ takes on a value coprime to both $a$ and $b$. By the Chinese remainder theorem, it suffices to show that for all primes $p$ dividing $a$ or $b$, $Q_{\alpha}(1, n)-a \not\equiv 0 \bmod p$. Whenever $a$ and $b$ are both even, $c$ is odd and $Q_{\alpha}(1, n)-a \not\equiv 0 \bmod 2$.

Suppose for some $p>2$ dividing $a$ or $b$, $Q_{\alpha}(1, n)-a \equiv 0 \bmod p$. Then, since $Q_{\alpha}(1, n)-a$ represents $c$, $d$, and $4a+4b+2c+d$, these are all divisible by $p$. Hence $p\mid a,b,c,d$. But this is impossible in a primitive packing. 
\end{proof}

\begin{lemma} \label{OctInsert2}
    For any circles $\mathcal{C}_a$ and $\mathcal{C}_b$ in a primitive octahedral circle packing where $a$ is odd and $b$ is even, there exists a circle $\mathcal{C}'$ tangent to both $\mathcal{C}_a$ and $\mathcal{C}_b$ whose curvature is even and coprime to $a$.
\end{lemma}
\begin{proof} Note that $Q_{\alpha}(1, n)-a$ will always be even in this case. Again, it suffices to show that for all primes $p$ dividing $a$, $Q_{\alpha}(1, n)-a \not\equiv 0 \bmod p$. It cannot be that $p=2$ since $a$ is odd, and for any $p>2$, the argument of Lemma \ref{OctInsert1} may be used. 
\end{proof}

\begin{corollary} \label{OctPath}
Let $\mathcal{C},\mathcal{C}' \in \mathcal{A}$ be two circles of nonzero curvature in a primitive octahedral circle packing. Then there exists a path of circles $\mathcal{C}_1,\mathcal{C}_2,...,\mathcal{C}_k$ such that \begin{enumerate}
    \item $\mathcal{C}_1=\mathcal{C}$ and $\mathcal{C}_k=\mathcal{C}'$;
    \item $\mathcal{C}_i$ is tangent to $\mathcal{C}_{i+1}$ for all $1 \leq i \leq k-1$;
    \item The curvatures of $\mathcal{C}_i$ and $\mathcal{C}_{i+1}$ are non-zero and coprime for all $1 \leq i \leq k-1$.
\end{enumerate}
\end{corollary}
\begin{proof} Clearly, there exists a path of circles that satisfies the first two requirements. Then, for any consecutive circles $\mathcal{C}_i$ and $\mathcal{C}_{i+1}$ in the path whose curvatures are not coprime, we can insert one circle between them with Lemma \ref{OctInsert1} if their curvatures are both even, or two circles between them with Lemma \ref{OctInsert2} and Lemma \ref{OctInsert1} if one of their curvatures is odd and the other is even, such that the path satisfies all three requirements. 
\end{proof}

Proposition \ref{OctEdge} and Corollary \ref{OctPath} imply the following result. 
\begin{proposition} \label{OctChi}
The value of $\chi_2$ is constant across all circles in a fixed primitive octahedral circle packing $\mathcal{A}$ of type $(0,1,2)$ or $(2,4,7)$ mod $8$.
\end{proposition}

Based on this Proposition, we will refer to $\chi_2(\mathcal{A})$ as a quadratic invariant of the packing. This is used to prove the main result of this section.

\begin{theorem} \label{OctObstruction}
In primitive integral octahedral packing $\mathcal{A}$ of type $(0,1,2)$ or $(2,4,7)$ mod $8$ with $\chi_2(\mathcal{A})=-1$, no integers of the form $n^2$ or $2n^2$ appear as curvatures. In particular, the Local-Global Conjecture is false for these packings.
\end{theorem}
\begin{proof}
Suppose that some circle in $\mathcal{A}$ has curvature $n^2$. Choose a tangent circle of curvature $a$ with $\gcd(a, n^2)=1$. Then 
\begin{equation*}
\chi_2(\mathcal{A}) = \res{n^2+a}{a'} = \res{n^2}{a'} =1
\end{equation*}
Similarly, if a circle in $\mathcal{A}$ has curvature $2n^2$, we choose a tangent circle of curvature $a$ with $\gcd(a, 2n^2)=1$. Then 
\begin{equation*}
\chi_2(\mathcal{A}) = \res{2n^2+a}{a} = \res{2n^2}{a} = \res{2}{a} \res{n^2}{a} = 1
\end{equation*}
using the fact that $a \equiv \pm 1 \bmod 8$. In either case, we contradict $\chi_2(\mathcal{A})=-1$. 
\end{proof}

To conclude, we describe the situation in packings of type $(0,3,6)$ or $(4,5,6)$ mod $8$. In type $(0,3,6)$, for two tangent circles $\mathcal{C}_1$, $\mathcal{C}_2$ of coprime curvatures $a$, $b$, we define
\begin{equation*} 
\chi_2(\mathcal{C}_1, \mathcal{C}_2)=\left\lbrace \begin{array}{cc} \res{a+b}{a} & a \not \equiv 0, 6 \bmod 8\\  \res{a+b}{a/2} & a \equiv 6 \bmod 8 \\  \res{-(a+b)}{a} & a \equiv 0 \bmod 8 \end{array} \right.
\end{equation*}
In type $(4,5,6)$ we define $\chi_2(\mathcal{C}_1, \mathcal{C}_2) = \res{a+b}{a}$. 

A reciprocity argument proves the analogue of Proposition \ref{OctEdge} 
\begin{proposition} \label{OctPartialEdge}
    Suppose that $\mathcal{C}_1$, $\mathcal{C}_2$ are tangent circles with coprime curvatures $a$, $b$ in an integral octahedral packing of type $(0,3,6)$ or $(4,5,6)$ mod $8$. Then $\chi_2(\mathcal{C}_1, \mathcal{C}_2)=\chi_2(\mathcal{C}_2, \mathcal{C}_1)$.
\end{proposition}
Proposition \ref{OctNode} implies that $\chi_2(\C_1, C_2)$ does not take the same value for all $\C_2$ tangent to $\C_1$ with coprime curvature. In fact, $\chi_2(\C_1, \C_2)=\chi_2(\C_1, \C_3)$ if $\C_2, \C_3$ have curvatures represented by the same quadratic form; $\chi_2(\C_1, \C_2)=-\chi_2(\C_1, \C_3)$ otherwise. Note that we only have the situation when $\chi_2(\C_1, \C_2)=-\chi_2(\C_1, \C_3)$ when $a$ is odd. 

Now we give the analogue of Proposition \ref{OctChi}, describing a $\chi_2$ invariant for subsets of circles in a packing. We can color the vertices of an octahedron with three colors so that diagonal vertices are colored identically and adjacent vertices are colored differently. This gives rise to the unique $3$-coloring of the octahedral packing such that no two tangent circles are the same color. Because of Proposition \ref{OctModular}, the odd circles in the packing will all be the same color. Furthermore, if we fix a circle of one color and parametrize all circles tangent to it with Proposition \ref{OctParametrization}, the $\alpha$-circles will all be the same color and the $\beta$-circles will all be the same color. 

\begin{figure}[h] \label{OctPartialFig}
    \centering
    \include{octpartial.tex}
    \vspace{-.25in}
    \caption{Octahedral Packing with Partial $\chi_2$ Invariant: No Blue Curvatures $n^2$ and No Red Curvatures $2n^2$}
\end{figure}

\begin{proposition} \label{OctPartialChi}
Suppose that a primitive integral octahedral packing $\mathcal{A}$ of type $(0,3,6)$ or $(4,5,6)$ mod $8$ is $3$-colored so that the odd circles are yellow and the even circles are red and blue. Then $\chi_2(\mathcal{C}_1, \mathcal{C}_2)$ is constant for all pairs of tangent, coprime circles colored yellow and red. It takes the opposite, constant value for all pairs of tangent, coprime circles colored yellow and blue. 
\end{proposition}
Note that $\chi_2(\mathcal{C}_1, \mathcal{C}_2)$ is undefined for a red-blue pair. This Proposition follows from Proposition \ref{OctPartialEdge}, Proposition \ref{OctNode}, and Corollary \ref{OctPath}. 

We obtain partial quadratic obstructions. 
\begin{proposition} \label{OctPartialObstruction}
    Suppose that a primitive integral octahedral packing $\mathcal{A}$ of type $(0,3,6)$ or $(4,5,6)$ mod $8$ is $3$-colored so that the odd circles are yellow and the even circles are red and blue. Suppose that $\chi_2(\mathcal{C}_1, \mathcal{C}_2)=1$ for yellow-red pairs and $\chi_2(\mathcal{C}_1, \mathcal{C}_2)=-1$ for yellow-blue pairs. Then no integers of the form $n^2$ appear as curvatures of blue circles and no integers of the form $2n^2$ appear as curvatures of red circles. 
\end{proposition}

There are also no yellow circles of curvature $n^2$, but this is already ruled out by modular restrictions. However, we do not find a quadratic obstruction across the entire packing. And computational evidence suggests that such obstructions do not exist. 

\section{Cubic Packing}
The cubic packing is built from a configuration of eight circles with a cubic tangency graph as shown.
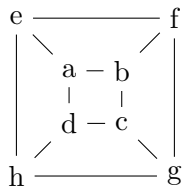
\begin{center}
\begin{tikzpicture}[scale=.35]
    \node (D) at (-1,-1)    {d};
    \node (C) at (1,-1)    {c};
    \node (B) at (1,1)    {b};
    \node (A) at (-1,1)    {a};
    \node (H) at (-3,-3)    {h};
    \node (G) at (3,-3)    {g};
    \node (F) at (3,3)    {f};
    \node (E) at (-3,3)    {e};
    
    \draw (A) -- (B);
    \draw (B) -- (C);
    \draw (C) -- (D);
    \draw (D) -- (A);
    \draw (A) -- (E);
    \draw (B) -- (F);
    \draw (C) -- (G);
    \draw (D) -- (H);
    \draw (E) -- (F);
    \draw (F) -- (G);
    \draw (G) -- (H);
    \draw (H) -- (E);
\end{tikzpicture}
\end{center}
The eight curvatures in a cubic circle configuration are subject to the following linear and quadratic relations
\begin{align}
    &a+g=b+h=c+e=d+f \label{CubeBodyDiagonals} \\
    &a+c=b+d  \label{CubeFaceDiagonals} \\
    &3(a+c+f+h)^2=8(a^2+c^2+f^2+h^2) \label{CubeQuadratic} 
\end{align}
proven in \cite{Stange}. We refer to the linear relations \eqref{CubeBodyDiagonals} as body diagonal relations and the relation \eqref{CubeFaceDiagonals} as a face diagonal relation. Together, these imply similar face diagonal relations for all faces: $a+f=b+e$, $b+g=c+f,$ $c+h=d+g$, $a+h=d+e$, $e+g=f+h$. There is also another quadratic relation $3(b+d+e+g)^2=8(b^2+d^2+e^2+g^2)$ implied by \eqref{CubeQuadratic} and the linear relations. 

Set $w=a+g=b+h=c+e=d+f$, a positive integer.  Using the relations \eqref{CubeBodyDiagonals} and \eqref{CubeFaceDiagonals}, we can rewrite Equation \eqref{CubeQuadratic} as $w^2-4(a+c)w+2(a^2+b^2+c^2+d^2)=0$, or $w=2(a+c) \pm 2\sqrt{ac+bd}$. Thus $w$ is even and $ac+bd=m^2$ for some nonnegative integer $m$.  

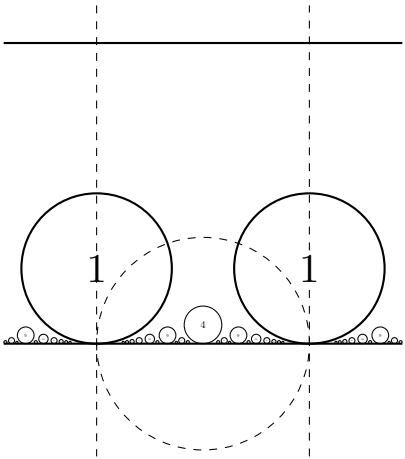
\begin{figure}[h] \label{CubeFordFig}
    \begin{center}
    \begin{tikzpicture}
        \draw[black, dashed] (0.0,0.0) circle (1.4142135623730951);
        \draw[black, dashed] (1.4142135623730951,-1.5) -- (1.4142135623730951,4.5);
        \draw[black, dashed] (-1.4142135623730951,4.5) -- (-1.4142135623730951,-1.5);
        \draw[black, thick] (-2.65,-0.0) -- (2.65,0.0);
        \draw[black, thin] (0.0,0.25000001102611014) circle (0.25000001102611014)node[scale=0.35355340618654824]{4};
        \draw[black, thick] (-1.4142135623730951,0.9999996905511437) circle (0.9999996905511437)node[scale=1.4142131247463254]{1};
        \draw[black, thick] (1.4142135623730951,0.9999996905511437) circle (0.9999996905511437)node[scale=1.4142131247463254]{1};
        \draw[black, thin] (-0.47140455782806895,0.11111111164674763) circle (0.11111111164674763)node[scale=0.15713484102118164]{9};
        \draw[black, thin] (0.47140455782806895,0.11111111164674763) circle (0.11111111164674763)node[scale=0.15713484102118164]{9};
        \draw[black, thin] (-0.28284272047461906,0.04000000006712127) circle (0.04000000006712127)node[scale=0.0565685425898476]{25};
        \draw[black, thin] (0.28284272047461906,0.04000000006712127) circle (0.04000000006712127)node[scale=0.0565685425898476]{25};
        \draw[black, thin] (-0.20203050891044216,0.020408163134628805) circle (0.020408163134628805)node[scale=0.02886150108811467]{49};
        \draw[black, thin] (0.20203050891044216,0.020408163134628805) circle (0.020408163134628805)node[scale=0.02886150108811467]{49};
        \draw[black, thin] (-2.545584404271571,0.04000000006712127) circle (0.04000000006712127)node[scale=0.0565685425898476]{25};
        \draw[black, thin] (2.545584404271571,0.04000000006712127) circle (0.04000000006712127)node[scale=0.0565685425898476]{25};
        \draw[black, thin] (-1.8182745597858163,0.020408163134628805) circle (0.020408163134628805)node[scale=0.02886150108811467]{49};
        \draw[black, thin] (1.8182745597858163,0.020408163134628805) circle (0.020408163134628805)node[scale=0.02886150108811467]{49};
        \draw[black, thin] (-0.848528121423857,0.04000000006712127) circle (0.04000000006712127)node[scale=0.0565685425898476]{25};
        \draw[black, thin] (0.848528121423857,0.04000000006712127) circle (0.04000000006712127)node[scale=0.0565685425898476]{25};
        \draw[black, thin] (-0.6060915267313265,0.020408163134628805) circle (0.020408163134628805)node[scale=0.02886150108811467]{49};
        \draw[black, thin] (0.6060915267313265,0.020408163134628805) circle (0.020408163134628805)node[scale=0.02886150108811467]{49};
        \draw[black, thin] (-1.7677669529663689,0.015624999998597885) circle (0.015624999998597885)node[scale=0.02209708691009672]{64};
        \draw[black, thin] (1.7677669529663689,0.015624999998597885) circle (0.015624999998597885)node[scale=0.02209708691009672]{64};
        \draw[black, thin] (-2.2223355776067004,0.020408163134628805) circle (0.020408163134628805)node[scale=0.02886150108811467]{49};
        \draw[black, thin] (2.2223355776067004,0.020408163134628805) circle (0.020408163134628805)node[scale=0.02886150108811467]{49};
        \draw[black, thin] (-0.9428090323228042,0.027777777911686907) circle (0.027777777911686907)node[scale=0.03928371025529541]{36};
        \draw[black, thin] (0.9428090323228042,0.027777777911686907) circle (0.027777777911686907)node[scale=0.03928371025529541]{36};
        \draw[black, thin] (-2.626396595427585,0.020408163134628805) circle (0.020408163134628805)node[scale=0.02886150108811467]{49};
        \draw[black, thin] (2.626396595427585,0.020408163134628805) circle (0.020408163134628805)node[scale=0.02886150108811467]{49};
        \draw[black, thin] (-2.4748737341529163,0.015624999998597885) circle (0.015624999998597885)node[scale=0.02209708691009672]{64};
        \draw[black, thin] (2.4748737341529163,0.015624999998597885) circle (0.015624999998597885)node[scale=0.02209708691009672]{64};
        \draw[black, thin] (-2.3570225669181215,0.11111111164674763) circle (0.11111111164674763)node[scale=0.15713484102118164]{9};
        \draw[black, thin] (2.3570225669181215,0.11111111164674763) circle (0.11111111164674763)node[scale=0.15713484102118164]{9};
        \draw[black, thin] (-1.0101525445522108,0.020408163134628805) circle (0.020408163134628805)node[scale=0.02886150108811467]{49};
        \draw[black, thin] (1.0101525445522108,0.020408163134628805) circle (0.020408163134628805)node[scale=0.02886150108811467]{49};
        \draw[black, thick] (2.65,3.9999998235822463) -- (-2.65,3.9999998235822463);
        \draw[black, thin] (-0.7071067499365475,0.06249999999439154) circle (0.06249999999439154)node[scale=0.08838834764038687]{16};
        \draw[black, thin] (0.7071067499365475,0.06249999999439154) circle (0.06249999999439154)node[scale=0.08838834764038687]{16};
        \draw[black, thin] (-0.3535533905932738,0.015624999998597885) circle (0.015624999998597885)node[scale=0.02209708691009672]{64};
        \draw[black, thin] (0.3535533905932738,0.015624999998597885) circle (0.015624999998597885)node[scale=0.02209708691009672]{64};
        \draw[black, thin] (-2.1213203123096425,0.06249999999439154) circle (0.06249999999439154)node[scale=0.08838834764038687]{16};
        \draw[black, thin] (2.1213203123096425,0.06249999999439154) circle (0.06249999999439154)node[scale=0.08838834764038687]{16};
        \draw[black, thin] (-1.0606601717798212,0.015624999998597885) circle (0.015624999998597885)node[scale=0.02209708691009672]{64};
        \draw[black, thin] (1.0606601717798212,0.015624999998597885) circle (0.015624999998597885)node[scale=0.02209708691009672]{64};
        \draw[black, thin] (-1.9798990033223334,0.04000000006712127) circle (0.04000000006712127)node[scale=0.0565685425898476]{25};
        \draw[black, thin] (1.9798990033223334,0.04000000006712127) circle (0.04000000006712127)node[scale=0.0565685425898476]{25};
        \draw[black, thin] (-1.8856180924233863,0.027777777911686907) circle (0.027777777911686907)node[scale=0.03928371025529541]{36};
        \draw[black, thin] (1.8856180924233863,0.027777777911686907) circle (0.027777777911686907)node[scale=0.03928371025529541]{36};
    \end{tikzpicture}
    \end{center}
    \caption{Cubic Ford Circles and Duals}
\end{figure}

\begin{figure}[h]
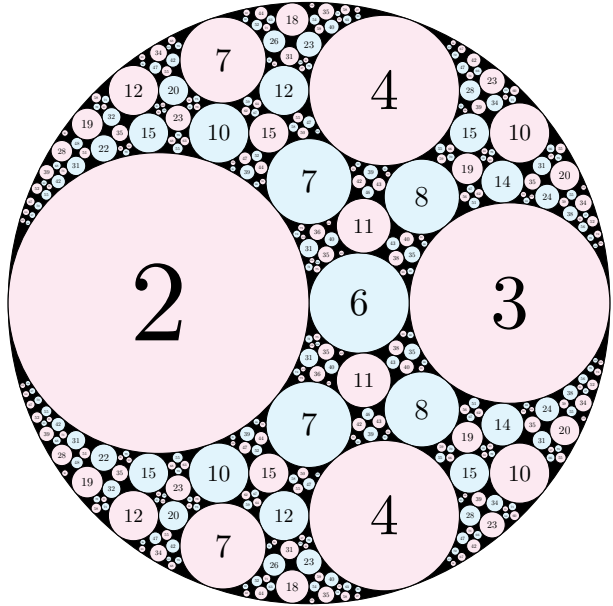
 \label{CubePartialFig}
    \begin{center}
    \include{cubepartial}
    \vspace{-.25in}
    \caption{Cubic Packing with Partial $\chi_2$ Invariant: No Blue Curvatures $n^2$, No Blue Curvatures $2n^2$ for $n$ Odd, and No Red Curvatures $8n^2$}
    \end{center}
\end{figure}

Stange conjectures that the modulus for congruence restrictions in this packing is 4 \cite{Stange}. Based on \cite[Theorem 8.1]{FuchsStangeZhang}, the discussion at the end of \cite[Section 9.2]{FuchsStangeZhang}, and the fact that the Apollonian group contains a congruence subgroup of level 2 (see Figure \ref{CubeFordFig}), we can conclude that the modulus divides 8. In Proposition \ref{CubeMod} we give a complete list of possible primitive integral cubic circle configurations mod 8. In Proposition \ref{CubeGroupMod}, we show that the Apollonian group acts transitively on the set of configurations mod 8 which restrict to a given configuration mod 4. Thus the modulus for congruence restrictions is 4, verifying Stange’s prediction.

\begin{proposition} \label{CubeMod} Any primitive cubic octuple of circles is congruent to one of the following mod $8$:

\begin{center} \, \hfill
\begin{tikzpicture}[scale=.35]
    \node (A) at (-1,1)    {0};
    \node (B) at (1,1)    {1};
    \node (C) at (1,-1)    {1};
    \node (D) at (-1,-1)    {0};
    \node (E) at (-3,3)    {1};
    \node (F) at (3,3)    {2};
    \node (G) at (3,-3)    {2};
    \node (H) at (-3,-3)    {1};
    
    \draw (A) -- (B);
    \draw (B) -- (C);
    \draw (C) -- (D);
    \draw (D) -- (A);
    \draw (A) -- (E);
    \draw (B) -- (F);
    \draw (C) -- (G);
    \draw (D) -- (H);
    \draw (E) -- (F);
    \draw (F) -- (G);
    \draw (G) -- (H);
    \draw (H) -- (E);
\end{tikzpicture} \hfill
\begin{tikzpicture}[scale=.35]
    \node (A) at (-1,1)    {0};
    \node (B) at (1,1)    {1};
    \node (C) at (1,-1)    {5};
    \node (D) at (-1,-1)    {4};
    \node (E) at (-3,3)    {1};
    \node (F) at (3,3)    {2};
    \node (G) at (3,-3)    {6};
    \node (H) at (-3,-3)    {5};
    
    \draw (A) -- (B);
    \draw (B) -- (C);
    \draw (C) -- (D);
    \draw (D) -- (A);
    \draw (A) -- (E);
    \draw (B) -- (F);
    \draw (C) -- (G);
    \draw (D) -- (H);
    \draw (E) -- (F);
    \draw (F) -- (G);
    \draw (G) -- (H);
    \draw (H) -- (E);
\end{tikzpicture} \hfill
\begin{tikzpicture}[scale=.35]
    \node (A) at (-1,1)    {4};
    \node (B) at (1,1)    {5};
    \node (C) at (1,-1)    {5};
    \node (D) at (-1,-1)    {4};
    \node (E) at (-3,3)    {5};
    \node (F) at (3,3)    {6};
    \node (G) at (3,-3)    {6};
    \node (H) at (-3,-3)    {5};
    
    \draw (A) -- (B);
    \draw (B) -- (C);
    \draw (C) -- (D);
    \draw (D) -- (A);
    \draw (A) -- (E);
    \draw (B) -- (F);
    \draw (C) -- (G);
    \draw (D) -- (H);
    \draw (E) -- (F);
    \draw (F) -- (G);
    \draw (G) -- (H);
    \draw (H) -- (E);
\end{tikzpicture} \hfill \,

\, \hfill \begin{tikzpicture}[scale=.35]
    \node (A) at (-1,1)    {0};
    \node (B) at (1,1)    {3};
    \node (C) at (1,-1)    {3};
    \node (D) at (-1,-1)    {0};
    \node (E) at (-3,3)    {3};
    \node (F) at (3,3)    {6};
    \node (G) at (3,-3)    {6};
    \node (H) at (-3,-3)    {3};
    
    \draw (A) -- (B);
    \draw (B) -- (C);
    \draw (C) -- (D);
    \draw (D) -- (A);
    \draw (A) -- (E);
    \draw (B) -- (F);
    \draw (C) -- (G);
    \draw (D) -- (H);
    \draw (E) -- (F);
    \draw (F) -- (G);
    \draw (G) -- (H);
    \draw (H) -- (E);
\end{tikzpicture} \hfill
\begin{tikzpicture}[scale=.35]
    \node (A) at (-1,1)    {0};
    \node (B) at (1,1)    {3};
    \node (C) at (1,-1)    {7};
    \node (D) at (-1,-1)    {4};
    \node (E) at (-3,3)    {3};
    \node (F) at (3,3)    {6};
    \node (G) at (3,-3)    {2};
    \node (H) at (-3,-3)    {7};
    
    \draw (A) -- (B);
    \draw (B) -- (C);
    \draw (C) -- (D);
    \draw (D) -- (A);
    \draw (A) -- (E);
    \draw (B) -- (F);
    \draw (C) -- (G);
    \draw (D) -- (H);
    \draw (E) -- (F);
    \draw (F) -- (G);
    \draw (G) -- (H);
    \draw (H) -- (E);
\end{tikzpicture} \hfill
\begin{tikzpicture}[scale=.35]
    \node (A) at (-1,1)    {4};
    \node (B) at (1,1)    {7};
    \node (C) at (1,-1)    {7};
    \node (D) at (-1,-1)    {4};
    \node (E) at (-3,3)    {7};
    \node (F) at (3,3)    {2};
    \node (G) at (3,-3)    {2};
    \node (H) at (-3,-3)    {7};
    
    \draw (A) -- (B);
    \draw (B) -- (C);
    \draw (C) -- (D);
    \draw (D) -- (A);
    \draw (A) -- (E);
    \draw (B) -- (F);
    \draw (C) -- (G);
    \draw (D) -- (H);
    \draw (E) -- (F);
    \draw (F) -- (G);
    \draw (G) -- (H);
    \draw (H) -- (E);
\end{tikzpicture} \hfill \,
\end{center}
\end{proposition}

\begin{lemma} \label{CubeModLemma1} Let $(a,b,c,d,e,f,g,h)$ be a primitive cubic octuple. Then two of $a, c, f, h$ are even and two are odd. 
\end{lemma}
Note that since $w$ is even, $g$, $e$, $d$, $b$ have the same parities as $a$, $c$, $f$, $h$ respectively. 
\begin{proof}
The quadratic relation \eqref{CubeQuadratic} implies that $a+c+f+h$ is even. If $a, c, f, h$ are all odd, then $a^2+c^2+f^2+h^2 \equiv 4 \bmod 8$, which means that the right side of \eqref{CubeQuadratic} has $2$-adic valuation $5$. However, the left side has an even $2$-adic valuation, a contradiction. If $a, c, f, h$ are all even, then $b, d, e, g$ are all even, and the octuple is imprimitive. We conclude that two of $a, c, f, h$ are even and two are odd. 
\end{proof}

\begin{lemma} \label{CubeModLemma2} Let $(a,b,c,d,e,f,g,h)$ be a primitive cubic octuple. Then $w \equiv 2 \bmod 4$.
\end{lemma}
\begin{proof}
From Lemma \ref{CubeModLemma1}, we have $a^2+c^2+f^2+h^2 \equiv 2 \bmod 4$, so the right side of \eqref{CubeQuadratic} has $2$-adic valuation $4$. Comparing to the left side, we find $a+c+f+h \equiv 4 \bmod 8$. By the relations \eqref{CubeFaceDiagonals} and \eqref{CubeBodyDiagonals}, $a+c+f+h=2w$, so $w \equiv 2 \bmod 4$
\end{proof}

\begin{lemma} \label{CubeModLemma3} Let $(a,b,c,d,e,f,g,h)$ be a primitive cubic octuple. Then the face diagonal sums $a+c, \, b+d, \, a+h, \, d+e, \, a+f, \, b+e, \, b+g, \, c+f, \, c+h, \, d+g, \, e+g, \, f+h$ are nonzero mod $4$. 
\end{lemma}

\begin{proof}
Consider the face $a$, $b$, $c$, $d$. From Lemma \ref{CubeModLemma1}, we can conclude that two of these integers are even and two are odd. If $a+c =b+d \equiv 0 \bmod 4$, we may assume without loss of generality that $a$ and $c$ are odd while $b$ and $d$ are even. Then $ac+bd \equiv 3 \bmod 4$, but we know that $ac+bd$ is a perfect square, a contradiction.
\end{proof}

Consider a primitive integral cubic octuple $(a, b, c, d, e, f, g, h)$. Based on Lemma \ref{CubeModLemma1}, we may assume without loss of generality that $a$, $d$, $f$, $g$ are even while $b$, $c$, $e$, $h$ are odd. From Lemma \ref{CubeModLemma2}, we may assume without loss of generality that $a \equiv 0$, $g \equiv 2$ mod $4$. By Lemma \ref{CubeModLemma3}, $f \equiv 2$ and $d \equiv 0$ mod $4$. And by Lemmas \ref{CubeModLemma2} and \ref{CubeModLemma3}, $b \equiv c \equiv e \equiv h \bmod 4$. This proves that any primitive integral cubic octuple is congruent to one of the following mod 4:

\begin{center} \, \hfill
\begin{tikzpicture}[scale=.35]
    \node (D) at (-1,-1)    {0};
    \node (C) at (1,-1)    {1};
    \node (B) at (1,1)    {1};
    \node (A) at (-1,1)    {0};
    \node (H) at (-3,-3)    {1};
    \node (G) at (3,-3)    {2};
    \node (F) at (3,3)    {2};
    \node (E) at (-3,3)    {1};
    
    \draw (A) -- (B);
    \draw (B) -- (C);
    \draw (C) -- (D);
    \draw (D) -- (A);
    \draw (A) -- (E);
    \draw (B) -- (F);
    \draw (C) -- (G);
    \draw (D) -- (H);
    \draw (E) -- (F);
    \draw (F) -- (G);
    \draw (G) -- (H);
    \draw (H) -- (E);
\end{tikzpicture} \hfill
\begin{tikzpicture}[scale=.35]
    \node (D) at (-1,-1)    {0};
    \node (C) at (1,-1)    {3};
    \node (B) at (1,1)    {3};
    \node (A) at (-1,1)    {0};
    \node (H) at (-3,-3)    {3};
    \node (G) at (3,-3)    {2};
    \node (F) at (3,3)    {2};
    \node (E) at (-3,3)    {3};
    
    \draw (A) -- (B);
    \draw (B) -- (C);
    \draw (C) -- (D);
    \draw (D) -- (A);
    \draw (A) -- (E);
    \draw (B) -- (F);
    \draw (C) -- (G);
    \draw (D) -- (H);
    \draw (E) -- (F);
    \draw (F) -- (G);
    \draw (G) -- (H);
    \draw (H) -- (E);
\end{tikzpicture} \hfill \,
\end{center}

The next lemma is the analogue of \cite[Proposition 3.1]{HKRS} and Lemma \ref{OctModLemma5}. 

\begin{lemma} \label{CubeModLemma4}
If $a$, $b$ are the curvatures of tangent circles in a primitive cubic octuple, then $a+b\not \equiv 5,7 \bmod 8$. 
\end{lemma}

\begin{proof}
We have 
\begin{equation*}
(a+b)(a+d)=a^2+a(b+d)+bd = a^2+a(a+c)+bd=2a^2+ac+bd+2a^2+m^2
\end{equation*}
Suppose that $p$ is an odd prime number which divides $2a^2+m^2$ but does not divide $a$ or $m$. Then $-2$ is a square modulo $p$, so $p$ is congruent to $1$ or $3$ modulo $8$. If a prime $p$ congruent to $5$ or $7$ modulo $8$ divides $2a^2+m^2$, then it must divide both $a$ and $m$, and the multiplicity $v_p(2a^2+m^2)$ must be even. 

Now suppose toward a contradiction that $a+b \equiv 5$ or $7 \bmod 8$. Then $a+b$ must have a prime factor $p$ congruent to $5$ or $7$ modulo $8$ with an odd multiplicity. Because $p$ divides $2a^2+m^2$, we have $p\mid a$, $p\mid m$, and $v_p(2a^2+m^2)$ is even. Thus $p\mid (a+c)$. Since $p$ is a common factor of $a$, $a+b$, and $a+c$, the packing is not primitive.
\end{proof}

Based on this lemma, we see that each mod 4 configuration shown above lifts to a mod 8 configuration in four possible ways, two of which are mirror images of each other. Two of the lifts for each configuration have $w \equiv 2 \bmod 8$ and two have $w \equiv 6 \bmod 8$. This proves Proposition \ref{CubeMod}.

Now we explain the mod 8 behavior across the entire packing. 

\begin{proposition} \label{CubeGroupMod}
The cubic Apollonian group acts on the set of primitive cubic octuples modulo 8 with two orbits, as shown. Each octuple modulo $8$ is fixed by four generators and transformed to a different octuple by the two generators shown.

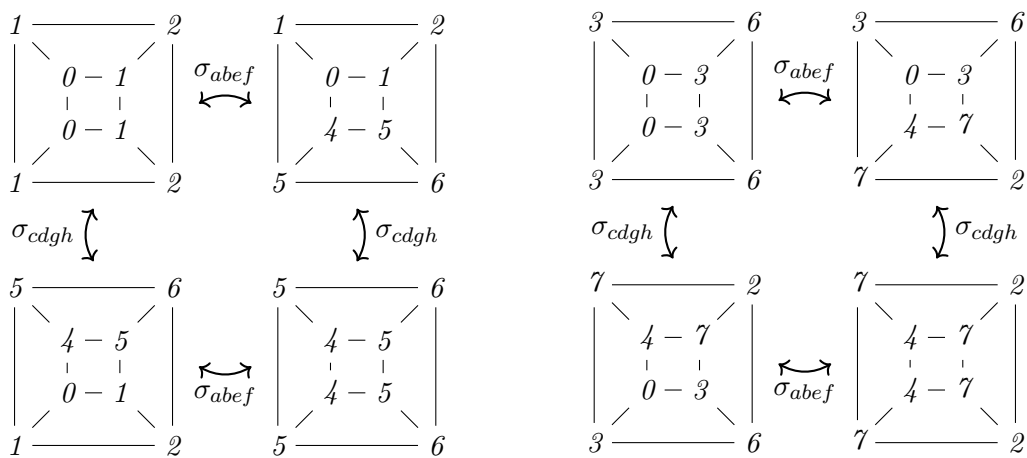
\begin{center}
\, \hfill
\begin{tikzpicture}[scale=.35]
    \node (A) at (-1,1)    {0};
    \node (B) at (1,1)    {1};
    \node (C) at (1,-1)    {1};
    \node (D) at (-1,-1)    {0};
    \node (E) at (-3,3)    {1};
    \node (F) at (3,3)    {2};
    \node (G) at (3,-3)    {2};
    \node (H) at (-3,-3)    {1};
    \draw (A) -- (B);
    \draw (B) -- (C);
    \draw (C) -- (D);
    \draw (D) -- (A);
    \draw (A) -- (E);
    \draw (B) -- (F);
    \draw (C) -- (G);
    \draw (D) -- (H);
    \draw (E) -- (F);
    \draw (F) -- (G);
    \draw (G) -- (H);
    \draw (H) -- (E);
    
    \begin{scope}[shift={(10,0)}]
    \node (A) at (-1,1)    {0};
    \node (B) at (1,1)    {1};
    \node (C) at (1,-1)    {5};
    \node (D) at (-1,-1)    {4};
    \node (E) at (-3,3)    {1};
    \node (F) at (3,3)    {2};
    \node (G) at (3,-3)    {6};
    \node (H) at (-3,-3)    {5};
    \draw (A) -- (B);
    \draw (B) -- (C);
    \draw (C) -- (D);
    \draw (D) -- (A);
    \draw (A) -- (E);
    \draw (B) -- (F);
    \draw (C) -- (G);
    \draw (D) -- (H);
    \draw (E) -- (F);
    \draw (F) -- (G);
    \draw (G) -- (H);
    \draw (H) -- (E);
    \end{scope}
    
    \begin{scope}[shift={(0,-10)}]
    \node (A) at (-1,1)    {4};
    \node (B) at (1,1)    {5};
    \node (C) at (1,-1)    {1};
    \node (D) at (-1,-1)    {0};
    \node (E) at (-3,3)    {5};
    \node (F) at (3,3)    {6};
    \node (G) at (3,-3)    {2};
    \node (H) at (-3,-3)    {1};
    \draw (A) -- (B);
    \draw (B) -- (C);
    \draw (C) -- (D);
    \draw (D) -- (A);
    \draw (A) -- (E);
    \draw (B) -- (F);
    \draw (C) -- (G);
    \draw (D) -- (H);
    \draw (E) -- (F);
    \draw (F) -- (G);
    \draw (G) -- (H);
    \draw (H) -- (E);
    \end{scope}

    \begin{scope}[shift={(10,-10)}]
    \node (A) at (-1,1)    {4};
    \node (B) at (1,1)    {5};
    \node (C) at (1,-1)    {5};
    \node (D) at (-1,-1)    {4};
    \node (E) at (-3,3)    {5};
    \node (F) at (3,3)    {6};
    \node (G) at (3,-3)    {6};
    \node (H) at (-3,-3)    {5};
    \draw (A) -- (B);
    \draw (B) -- (C);
    \draw (C) -- (D);
    \draw (D) -- (A);
    \draw (A) -- (E);
    \draw (B) -- (F);
    \draw (C) -- (G);
    \draw (D) -- (H);
    \draw (E) -- (F);
    \draw (F) -- (G);
    \draw (G) -- (H);
    \draw (H) -- (E);
    \end{scope}
    
    \draw [thick, <->] (4, 0) to[bend left=30] node[above] {$\sigma_{abef}$} (6, 0);
    \draw [thick, <->] (0, -4) to[bend right=30] node[left] {$\sigma_{cdgh}$} (0, -6);
    \draw [thick, <->] (10, -4) to[bend left=30] node[right] {$\sigma_{cdgh}$} (10, -6);
    \draw [thick, <->] (4, -10) to[bend right=30] node[below] {$\sigma_{abef}$} (6, -10);
\end{tikzpicture}
\hfill
\begin{tikzpicture}[scale=.35]
    \node (A) at (-1,1)    {0};
    \node (B) at (1,1)    {3};
    \node (C) at (1,-1)    {3};
    \node (D) at (-1,-1)    {0};
    \node (E) at (-3,3)    {3};
    \node (F) at (3,3)    {6};
    \node (G) at (3,-3)    {6};
    \node (H) at (-3,-3)    {3};
    \draw (A) -- (B);
    \draw (B) -- (C);
    \draw (C) -- (D);
    \draw (D) -- (A);
    \draw (A) -- (E);
    \draw (B) -- (F);
    \draw (C) -- (G);
    \draw (D) -- (H);
    \draw (E) -- (F);
    \draw (F) -- (G);
    \draw (G) -- (H);
    \draw (H) -- (E);
    
    \begin{scope}[shift={(10,0)}]
    \node (A) at (-1,1)    {0};
    \node (B) at (1,1)    {3};
    \node (C) at (1,-1)    {7};
    \node (D) at (-1,-1)    {4};
    \node (E) at (-3,3)    {3};
    \node (F) at (3,3)    {6};
    \node (G) at (3,-3)    {2};
    \node (H) at (-3,-3)    {7};
    \draw (A) -- (B);
    \draw (B) -- (C);
    \draw (C) -- (D);
    \draw (D) -- (A);
    \draw (A) -- (E);
    \draw (B) -- (F);
    \draw (C) -- (G);
    \draw (D) -- (H);
    \draw (E) -- (F);
    \draw (F) -- (G);
    \draw (G) -- (H);
    \draw (H) -- (E);
    \end{scope}
    
    \begin{scope}[shift={(0,-10)}]
    \node (A) at (-1,1)    {4};
    \node (B) at (1,1)    {7};
    \node (C) at (1,-1)    {3};
    \node (D) at (-1,-1)    {0};
    \node (E) at (-3,3)    {7};
    \node (F) at (3,3)    {2};
    \node (G) at (3,-3)    {6};
    \node (H) at (-3,-3)    {3};
    \draw (A) -- (B);
    \draw (B) -- (C);
    \draw (C) -- (D);
    \draw (D) -- (A);
    \draw (A) -- (E);
    \draw (B) -- (F);
    \draw (C) -- (G);
    \draw (D) -- (H);
    \draw (E) -- (F);
    \draw (F) -- (G);
    \draw (G) -- (H);
    \draw (H) -- (E);
    \end{scope}

    \begin{scope}[shift={(10,-10)}]
    \node (A) at (-1,1)    {4};
    \node (B) at (1,1)    {7};
    \node (C) at (1,-1)    {7};
    \node (D) at (-1,-1)    {4};
    \node (E) at (-3,3)    {7};
    \node (F) at (3,3)    {2};
    \node (G) at (3,-3)    {2};
    \node (H) at (-3,-3)    {7};
    \draw (A) -- (B);
    \draw (B) -- (C);
    \draw (C) -- (D);
    \draw (D) -- (A);
    \draw (A) -- (E);
    \draw (B) -- (F);
    \draw (C) -- (G);
    \draw (D) -- (H);
    \draw (E) -- (F);
    \draw (F) -- (G);
    \draw (G) -- (H);
    \draw (H) -- (E);
    \end{scope}
    
    \draw [thick, <->] (4, 0) to[bend left=30] node[above] {$\sigma_{abef}$} (6, 0);
    \draw [thick, <->] (0, -4) to[bend right=30] node[left] {$\sigma_{cdgh}$} (0, -6);
    \draw [thick, <->] (10, -4) to[bend left=30] node[right] {$\sigma_{cdgh}$} (10, -6);
    \draw [thick, <->] (4, -10) to[bend right=30] node[below] {$\sigma_{abef}$} (6, -10);
\end{tikzpicture} \hfill \,
\end{center}
 
\end{proposition}
\begin{proof}
Each generator for the cubic Apollonian group corresponds to a face of the cube. It fixes the four curvatures of the circles on that face of the tangency graph, and changes the value of $w$ and the four curvatures of the circles on the opposite face. For example, the generator $\sigma_{abcd}$ fixes $a, b, c, d$, and transforms $w$ to $4(a+c)-w$. The remaining four curvatures can be computed from the new $w$ value. Since each $w$ is $2$ or $6$ mod $8$, the generator fixes $w \bmod 8$ if the face diagonal sum $a+c$ is odd and changes it if the face diagonal sum is even. We conclude that each possible primitive cubic octuple modulo $8$ is fixed by four generators of the cubic Apollonian group and transformed by two generators, and that the orbits of the group are as illustrated in the statement of the Proposition.
\end{proof}

From this proposition, we see that there are exactly two types of cubic circle packing modulo $8$, and each type contains exactly six integers mod $8$, or three integers mod $4$, as curvatures. We will refer to these as type $(0,1,2)$ and type $(0,2,3)$ based on the curvatures mod $4$. Section \ref{Data} below contains examples of both packing types.

We now seek to  parameterize all circles tangent to a fixed circle within a cubic packing. Figure 1 shows a cubic Ford packing, with a fixed circle of curvature zero along the real axis. The group generated by the reflections in dual circles perpendicular to this axis is a congruence subgroup of level 2 in $\mathrm{SL}(2, \Z)$.

\begin{proposition} \label{CubeFord}
    The generalized Ford circles are parametrized by $x, y \in \Z$ with $y \geq 0$ and $\gcd(x,y)=1$. Their inversive coordinates are as follows:
\begin{equation*}
c(x,y)=(8x^2,y^2,2\sqrt{2}xy, 1)
 \end{equation*}
Each circle is tangent to the real axis is at $\frac{2 \sqrt{2} x}{y}$. 
\end{proposition}

This can be proven similarly to Proposition \ref{OctFord}, by checking Properties (1), (2), and (3) of the parametrization. The computation is simpler because the parametrization is given by a single formula. 



Next we parameterize the curvatures of all circles tangent to a fixed circle in a arbitrary cubic packing.

\begin{proposition} \label{CubeParametrization}
    Suppose a cubic packing contains an octuple of circles with curvatures $(a,b,c,d,e,f, g, h)$. Then all circles tangent to the circle of curvature $a$ are parametrized by $x, y \in \Z$ with $y \geq 0$, $\gcd(x,y)=1$. Their curvatures are as follows:
\begin{equation*}
    Q(x,y)-a=(a+b) x^2 - (a+b+d-e) x y +(a+d) y^2 - a
\end{equation*}   
\end{proposition}

\begin{proof}
We can apply a M\"{o}bius transformation which maps the circles $(0,0,0,-1)$, $(8,0,0, 1)$, $(0,1,0,1)$, $(8, 1, 2\sqrt{2}, 1)$ to circles of curvature $a$, $b$, $d$, $e$ respectively. This transformation is linear in the inversive coordinate system, so when applied to $c(x,y)$, it gives the stated formula.
\end{proof}

Using Equations \eqref{CubeBodyDiagonals}, \eqref{CubeFaceDiagonals}, and \eqref{CubeQuadratic}, we find that $Q$ has discriminant $-8a^2$.

We now have the information necessary to define a  quadratic invariant $\chi_2$ for cubic packings. Our $\chi_2$ will be well-defined across the entire packing in type $(0,1,2)$. We define a partial $\chi_2$ for certain pairs of circles in the packing in type $(0,2,3)$.

Fix a primitive integral cubic octuple $(a,b,c,d,e,f,g,h)$, with $a \neq 0$. The quadratic form $Q(x,y)$ represents the integers $a+b$, $a+d$, and $a+e$. Since $\gcd(a, b, d, e)=1$, for each prime $p$ dividing $a$, $Q$ represents integers coprime to $p$. Thus, by the Chinese remainder theorem, we can find an integer $\rho$ coprime to $a$ which is represented by $Q(x,y)$. We may take $y \geq 0$, $\gcd(x, y)=1$. 

\begin{proposition} \label{CubeNode}
Suppose $(a, b, c, d, e, f, g, h)$ is a primitive cubic octuple. Let $\rho_1$, $\rho_2$ be two integers coprime to $a$ which are represented by $Q$. Let $a'=a/2$ if $a \equiv 2 \bmod 4$ and $a'=a$ otherwise. Then 
\begin{equation*}
\left(\frac{\rho_1}{a'}\right) =\left(\frac{\rho_2}{a'}\right) 
\end{equation*} 
\end{proposition}

This Proposition follows from Lemma \ref{QuadFormKronecker} since the discriminant of $Q$, $-8a^2$, is divisible by $a'$ and a sufficient power of $2$.

If $\C$ is a circle of curvature $a$ in a packing of type $(0, 1, 2)$, we define:
\begin{equation} \label{CubeChi2}
    \chi_2(\mathcal{C})=
    \begin{cases}
        \left(\frac{\rho}{a} \right) & a \equiv 0, 1 \bmod 4 \\
        \left(\frac{-\rho}{a/2}\right) & a \equiv 2 \bmod 4 \\
    \end{cases}
\end{equation}

If $\C$ is a circle of curvature $a$ in a packing of type $(0, 2, 3)$, we define:
\begin{equation} \label{CubePartialChi2}
    \chi_2(\mathcal{C})=
    \begin{cases}
        \left(\frac{\rho}{a} \right) & a \equiv 3 \bmod 4 \\
        -\left(\frac{-\rho}{a} \right) & a \equiv 0 \bmod 4 \\
        -\left(\frac{\rho}{a/2}\right) & a \equiv 2 \bmod 4 \\
    \end{cases}
\end{equation}
Both expressions are well-defined by Proposition \ref{CubeNode}.

Now we study the propagation of $\chi_2$ from circle to circle.

\begin{proposition} \label{CubeEdge}
Suppose that $\mathcal{C}_a$, $\mathcal{C}_b$ are tangent circles with coprime curvatures $a$, $b$ in a primitive integral cubic packing of type $(0,1,2)$. Then
\begin{equation*}
\chi_2(\mathcal{C}_a) = \chi_2(\mathcal{C}_b)
\end{equation*}
\end{proposition}

\begin{proof}
We will use $\rho=a+b$ to compute $\chi_2(\mathcal{C}_a)$ and $\chi_2(\mathcal{C}_b)$. Assume without loss of generality that $a$ is odd, and thus, since we are in type $(0, 1, 2)$, $a \equiv 1 \bmod 4$. By quadratic reciprocity for the Kronecker symbol, $\res{a}{b}=\res{b}{a}$. Then if $b \equiv 0, 1 \bmod 4$, we have
\begin{equation*}
\chi_2(\C_a)=\res{b+a}{a}=\res{b}{a}=\res{a}{b} = \res{b+a}{b} = \chi_2(\C_b)
\end{equation*}
If $b \equiv 2 \bmod 4$, then based on Lemma \ref{CubeModLemma4}, we have $a \equiv 1 \bmod 8$ if and only if $b \equiv 2 \bmod 8$ and $a \equiv 5 \bmod 8$ if and only if $b \equiv 6 \bmod 8$, so $\res{a}{2}=\res{-1}{b/2}$. Then
\begin{equation*}
\chi_2(\C_a)=\res{b+a}{a}=\res{b}{a}=\res{a}{b}=\res{a}{2} \res{a}{b/2} = \res{-a}{b/2} = \res{-a-b}{b/2} = \chi_2(\C_b)
\end{equation*}
\end{proof}

To show that $\chi_2(\mathcal{C}$ is identical for all circles in the packing, we must find a path of coprime curvatures between any two circles. 

\begin{lemma} \label{CubeSimultaneous}
Suppose that $(a,b,c,d,e,f,g,h)$ is an octuple of circles in a cubic packing, so that the circles of curvatures $a,b$ form a square of tangencies with $c,d$, and also with $e,f$. Then all pairs of curvatures which form a square of tangencies with $a,b$ in the packing are parametrized by the following pair of formulas:
\begin{equation*}
    a \left(n^2-n\right)+b \left(n^2-n\right)+d (1-n)+e n, \quad a \left(n^2-n\right)+b \left(n^2-n\right)+c (1-n)+f n
\end{equation*} 
\end{lemma}
\begin{proof}
First we consider the strip packing bounded by the parallel lines with inversive coordinates $(0,0,0,-1)$ and $(8,0,0,1)$. The only pairs of circles which form a square tangency graph with these two lines have radius $1$ and centers at $(2\sqrt{2} n, 1)$ and $(2 \sqrt{2} n, 3)$. The inversive coordinates of these pairs are given by 
$(8n^2, 1, 2\sqrt{2}n, 1)$ and $
(8n^2+8, 1, 2\sqrt{2}n, 3)$ for $n \in \Z$. 

We can apply a M\"obius transformation which maps the four generalized circles $(0,0,0,-1)$, $(8,0,0,1)$, $(0,1,0,1)$, and $(8,1,2\sqrt{2},1)$ to circles of curvatures $a, b, d, e$ respectively. Since the M\"obius transformation is linear in inversive coordinates, it gives the stated formulas when applied to $(8n^2, 1, 2\sqrt{2}n, 1), \,
(8n^2+8, 1, 2\sqrt{2}n, 3)$.
\end{proof}

\begin{lemma} \label{CubeInsert}
Given a pair of tangent circles of curvatures $a$, $b$ in a primitive integral cubic packing, there exist circles of curvatures $k$, $\ell$ forming a square of tangencies with $a$, $b$, such that $\gcd(b, k)=\gcd(\ell,a)=1$ and $\gcd(k,\ell)$ is a power of $2$. 
\end{lemma}

\begin{proof}
We will need the fact that in a primitive cubic octuple $(a,b,c,d,e,f,g,h)$, $\gcd(a, d,e)=1$, $\gcd(b,c,f)=1$ and $\gcd(a-b, c, d, e, f)$ is a power of $2$. If a prime $p$ divides $a, \, d, \, e$, then by the face diagonal relation it divides $h$. By the quadratic formula for $w$, $p$ divides $w$ and thus divides the other four curvatures, contradicting primitivity. The proof that $\gcd(b,c,f)=1$ is similar. If an odd prime $p$ divides $a-b$, $c$, $d$, $e$, $f$, then it divides $w$. By the quadratic formula for $w$, $p$ divides $2a$, so $p$ divides $a$, $b$, $g$ and $h$, contradicting primitivity. 

Since the difference between the two parametrized curvatures in Lemma \ref{CubeSimultaneous} is $d-c+(c-d+e-f)n=d-c=a-b$, the greatest common divisor of these curvatures always divides $a-b$.  

Now we proceed to the main proof. For a prime $p$ dividing $a$, the quadratic expression $a \left(n^2-n\right)+b \left(n^2-n\right)+d (1-n)+e n$ takes the values $d$, $e$, so it takes values not divisible by $p$. Similarly, for a prime $p$ dividing $b$, the quadratic $a \left(n^2-n\right)+b \left(n^2-n\right)+c (1-n)+f n$ takes values not divisible by $p$. For an odd prime $p$ dividing $a-b$, the two quadratic expressions together take values $c$, $d$ ,$e$, $f$, so at least one of them takes a value not divisible by $p$. If $p$ divides $a$, $b$ and $a-b$, then both quadratics take values not divisible by $p$ simultaneously. 

By the Chinese remainder theorem, we can find $n\in\Z$, $k=a \left(n^2-n\right)+b \left(n^2-n\right)+d (1-n)+e n$, $\ell=a \left(n^2-n\right)+b \left(n^2-n\right)+c (1-n)+f n$, such that for all primes $p$ dividing $a$, $p$ does not divide $k$, for all primes $p$ dividing $b$, $p$ does not divide $\ell$, and for all odd primes $p$ dividing $a-b$, $p$ does not divide $k$ or $p$ does not divide $\ell$. This completes the proof. 
\end{proof}

\begin{corollary} \label{CubePath}
Let $\mathcal{C},\mathcal{C}' \in \mathcal{A}$ be two circles in a primitive cubic circle packing. Then there exists a path of circles $\mathcal{C}_1,\mathcal{C}_2,...,\mathcal{C}_k$ such that \begin{enumerate}
    \item $\mathcal{C}_1=\mathcal{C}$ and $\mathcal{C}_k=\mathcal{C}'$;
    \item $\mathcal{C}_i$ is tangent to $\mathcal{C}_{i+1}$ for all $1 \leq i \leq k-1$;
    \item The curvatures of $\mathcal{C}_i$ and $\mathcal{C}_{i+1}$ are coprime for all $1 \leq i \leq k-1$.
\end{enumerate}
\end{corollary}
\begin{proof} Clearly, there exists a path of circles that satisfies the first two requirements. Then, for any consecutive circles $\mathcal{C}_i$ and $\mathcal{C}_{i+1}$ in the path whose curvatures are not coprime, we can insert two circles between them with Lemma \ref{CubeInsert}. The modified path satisfies the third requirement except for possible common factors of $2$ between the two new circles. If either $\C_i$ or $\C_{i+1}$ has even curvature, then Proposition \ref{CubeMod} implies that one of the new circles has odd curvature, so the two new curvatures are coprime. If both $\C_i$ and $\C_{i+1}$ have odd curvature, then the two new curvatures will be even, and have a power of $2$ as a common factor. In this case, we use Lemma \ref{CubeInsert} to insert two additional circles between them, whose curvatures will be coprime. The resulting path satisfies the third requirement. 
\end{proof}

Proposition \ref{CubeEdge} and Corollary \ref{CubePath} imply the following result. 

\begin{proposition} \label{CubeChi}
The value of $\chi_2$ is constant across all circles in a fixed primitive cubic circle packing $\mathcal{A}$ of type $(0,1,2)$.
\end{proposition}

Based on this Proposition, we will refer to $\chi_2(\mathcal{A})$ as a quadratic invariant of the packing. This is used to prove the main result of this section.

\begin{theorem} \label{CubeObstruction}
In primitive integral cubic packing $\mathcal{A}$ of type $(0,1,2)$ with $\chi_2(\mathcal{A})=-1$, no integers of the form $n^2$ or $2n^2$ appear as curvatures. In particular, the Local-Global Conjecture is false for these packings.
\end{theorem}
\begin{proof}
Suppose that some circle in $\mathcal{A}$ has curvature $n^2$. Choose a tangent circle of curvature $a \equiv 0,1 \bmod 4$ with $\gcd(a, n^2)=1$. Then 
\begin{equation*}
\chi_2(\mathcal{A}) = \res{n^2+a}{a} = \res{n^2}{a} =1
\end{equation*}
Similarly, if a circle in $\mathcal{A}$ has curvature $2n^2$, we choose a tangent circle of curvature $a \equiv 1 \bmod 4$ with $\gcd(a, 2n^2)=1$. Then 
\begin{equation*}
\chi_2(\mathcal{A}) = \res{2n^2+a}{a} = \res{2n^2}{a} = \res{2}{a} \res{n^2}{a} = 1
\end{equation*}
using the fact $2n^2$ must be $0$ or $2$ mod $8$, so $a$ must be $1$ mod $8$ by Lemma \ref{CubeModLemma4}. In either case, we contradict $\chi_2(\mathcal{A})=-1$.
\end{proof}

We conclude by describing the situation in type $(0,2,3)$. Instead of Proposition \ref{CubeEdge}, we have the following result, with a similar proof:

\begin{proposition} \label{CubeEdgeAlt}
Suppose that $\mathcal{C}_a$, $\mathcal{C}_b$ are tangent circles with coprime curvatures $a$, $b$ in a primitive integral cubic packing of type $(0,2,3)$. Then
\begin{equation*}
\chi_2(\mathcal{C}_a) = -\chi_2(\mathcal{C}_b)
\end{equation*}
\end{proposition}

There does not seem to be a way of defining $\chi_2(\mathcal{C})$ that is consistent for all circles in the packing--the difficulty is that we have pairs of tangent circles with both curvatures congruent to $3$ mod $4$, and no further congruence information to distinguish them.

There is a unique two-coloring of the vertices of a cube such that adjacent vertices have different colors. This gives rise to a unique two-coloring of all circles in a cubic packing such that tangent circles have different colors. Proposition \ref{CubeEdgeAlt} and Corollary \ref{CubePath} imply the following result. 

\begin{proposition} \label{CubeChiAlt}
Suppose that the circles in a fixed primitive cubic circle packing $\mathcal{A}$ of type $(0,2,3)$ are two-colored blue and red. Then $\chi_2$ has one value for all blue circles, and the opposite value for all red circles. 
\end{proposition}

We obtain partial quadratic obstructions.

\begin{proposition}
Suppose that the circles in a fixed primitive cubic circle packing $\mathcal{A}$ of type $(0,2,3)$ are two-colored blue and red, so that blue circles have $\chi_2(\C)=1$ and red circles have $\chi_2(\C)=-1$. Then no blue circle has curvature $n^2$, no blue circle has curvature $2n^2$ for $n$ odd, and no red circle has curvature $8n^2$. 
\end{proposition}
\begin{proof}
Suppose that some circle in $\mathcal{A}$ has curvature $n^2$. Choose a tangent circle $\C$ of curvature $a \equiv 3 \bmod 4$ with $\gcd(a, n^2)=1$. Then 
\begin{equation*}
\chi_2(\C) = \res{n^2+a}{a} = \res{n^2}{a} =1
\end{equation*}
Similarly, if a circle in $\mathcal{A}$ has curvature $2n^2$ for $n$ odd, we choose a tangent circle of curvature $a \equiv 3 \bmod 4$ with $\gcd(a, 2n^2)=1$. Then 
\begin{equation*}
\chi_2(\C) = \res{2n^2+a}{a} = \res{2n^2}{a} = \res{2}{a} \res{n^2}{a} = 1
\end{equation*}
using the fact $2n^2$ must be $2$ mod $8$, so $a$ must be $7$ mod $8$ by Lemma \ref{CubeModLemma4}. In either case, we conclude that $\C$ is blue, so the circle of curvature $n^2$ or $2n^2$ is red. 

Similarly, if a circle in $\mathcal{A}$ has curvature $8n^2$, we choose a tangent circle of curvature $a \equiv 3 \bmod 4$ with $\gcd(a, 8n^2)=1$. Then 
\begin{equation*}
\chi_2(\mathcal{A}) = \res{2n^2+a}{a} = \res{2n^2}{a} = \res{2}{a} \res{n^2}{a} = -1
\end{equation*}
using the fact $8n^2$ must be $0$ mod $8$, so $a$ must be $3$ mod $8$ by Lemma \ref{CubeModLemma4}. We conclude that $\C$ is red, so the circle of curvature $8n^2$ is blue. 
\end{proof}

We do not obtain quadratic obstructions across all circles in the packing, and computational evidence suggests that such obstructions do not exist. 

\section{Square Packing}
The square packing is built from an infinite collection of circles with a square grid tangency graph as shown.
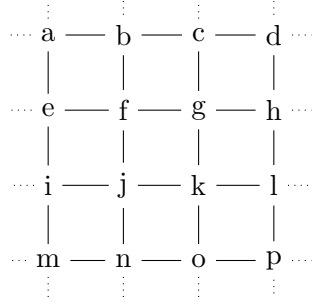
\begin{center}
\begin{tikzpicture}[scale=.5]
    \node (A) at (0,0)    {a};
    \node (B) at (2,0)    {b};
    \node (C) at (4,0)    {c};
    \node (D) at (6,0)    {d};
    \node (E) at (0,-2)    {e};
    \node (F) at (2,-2)    {f};
    \node (G) at (4,-2)    {g};
    \node (H) at (6,-2)    {h};
    \node (I) at (0,-4)    {i};
    \node (J) at (2,-4)    {j};
    \node (K) at (4,-4)    {k};
    \node (L) at (6,-4)    {l};
    \node (M) at (0,-6)    {m};
    \node (N) at (2,-6)    {n};
    \node (O) at (4,-6)    {o};
    \node (P) at (6,-6)    {p};

    \draw (A) -- (B);
    \draw (B) -- (C);
    \draw (C) -- (D);
    \draw (E) -- (F);
    \draw (F) -- (G);
    \draw (G) -- (H);
    \draw (I) -- (J);
    \draw (J) -- (K);
    \draw (K) -- (L);
    \draw (M) -- (N);
    \draw (N) -- (O);
    \draw (O) -- (P);

    \draw (A) -- (E);
    \draw (E) -- (I);
    \draw (I) -- (M);
    \draw (B) -- (F);
    \draw (F) -- (J);
    \draw (J) -- (N);
    \draw (C) -- (G);
    \draw (G) -- (K);
    \draw (K) -- (O);
    \draw (D) -- (H);
    \draw (H) -- (L);
    \draw (L) -- (P);
    
    \draw[dotted] (A) -- (-1,0);
    \draw[dotted] (E) -- (-1,-2);
    \draw[dotted] (I) -- (-1,-4);
    \draw[dotted] (M) -- (-1,-6);
    \draw[dotted] (D) -- (7,0);
    \draw[dotted] (H) -- (7,-2);
    \draw[dotted] (L) -- (7,-4);
    \draw[dotted] (P) -- (7,-6);
    \draw[dotted] (A) -- (0,1);
    \draw[dotted] (B) -- (2,1);
    \draw[dotted] (C) -- (4,1);
    \draw[dotted] (D) -- (6,1);
    \draw[dotted] (M) -- (0,-7);
    \draw[dotted] (N) -- (2,-7);
    \draw[dotted] (O) -- (4,-7);
    \draw[dotted] (P) -- (6,-7);
\end{tikzpicture}
\end{center}
This infinite collection of circles is only part of the packing; we also have infinitely many copies of the configuration obtained by iteratively reflecting it through dual circles. An example of a bounded square packing is shown in Figure \ref{SquarePackingFigure}. 

\begin{figure}[h] 
\begin{center}
\vspace{-.25in}
\include{justsquare}
\vspace{-.25in}
\end{center}
\caption{Square Packing}
\label{SquarePackingFigure}
\end{figure}

The curvatures in a square grid circle configuration are subject to the following linear and quadratic relations
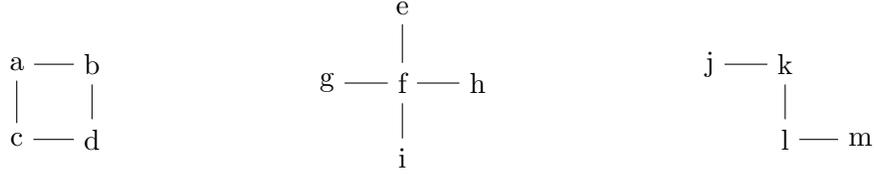
\begin{center}
\, \hfill
\begin{tikzpicture}[scale=.5]
    \node (A) at (0,0)    {a};
    \node (C) at (0,-2)    {c};
    \node (B) at (2,0)    {b};
    \node (D) at (2,-2)    {d};
    \draw (A) -- (B);
    \draw (A) -- (C);
    \draw (D) -- (C);
    \draw (D) -- (B);
    \draw[white] (D) -- (0, -3);
\end{tikzpicture}
\hfill
\begin{tikzpicture}[scale=.5]
    \node (A) at (0,0)    {e};
    \node (C) at (0,-2)    {f};
    \node (B) at (-2,-2)    {g};
    \node (D) at (2,-2)    {h};
    \node (E) at (0,-4)    {i};
    \draw (A) -- (C);
    \draw (B) -- (C);
    \draw (D) -- (C);
    \draw (C) -- (E);
\end{tikzpicture}
\hfill
\begin{tikzpicture}[scale=.5]
    \node (A) at (0,0)    {j};
    \node (B) at (2,0)    {k};
    \node (C) at (2,-2)    {l};
    \node (D) at (4,-2)    {m};
    \draw (A) -- (B);
    \draw (B) -- (C);
    \draw (D) -- (C);
    \draw[white] (D) -- (0, -3);
\end{tikzpicture}
\hfill \,
\end{center}

\begin{align}
    &a+d=b+c\label{SquareFaceDiagonal} \\
    &e+i=g+h \label{SquareT} \\
    &(j-3k)^2+(m-3l)^2=2(j+k)(m+l)
    \label{SquareQuadratic} 
\end{align}
proven in \cite{RWYY}. The relations apply to all configurations of circles arranged in these ways within the grid.

\begin{figure}[h] \label{SquareFordFig}
\begin{center}
\begin{tikzpicture}
\draw[black, thin] (-1.3333333333333333,-1.9444444444444444) circle (0.05555555555555555)node[scale=0.1111111111111111]{18};
\draw[black, thin] (1.3333333333333333,-1.9444444444444444) circle (0.05555555555555555)node[scale=0.1111111111111111]{18};
\draw[black, thin] (-0.4,-1.96) circle (0.04)node[scale=0.08]{25};
\draw[black, thin] (0.4,-1.96) circle (0.04)node[scale=0.08]{25};
\draw[black, thin] (-3.0,-1.875) circle (0.125)node[scale=0.25]{8};
\draw[black, thin] (3.0,-1.875) circle (0.125)node[scale=0.25]{8};
\draw[black, thin] (-2.857142857142857,-1.989795918367347) circle (0.01020408163265306)node[scale=0.02040816326530612]{98};
\draw[black, thin] (2.857142857142857,-1.989795918367347) circle (0.01020408163265306)node[scale=0.02040816326530612]{98};
\draw[black, thin] (-1.2,-1.96) circle (0.04)node[scale=0.08]{25};
\draw[black, thin] (1.2,-1.96) circle (0.04)node[scale=0.08]{25};
\draw[black, thin] (-0.5,-1.96875) circle (0.03125)node[scale=0.0625]{32};
\draw[black, thin] (0.5,-1.96875) circle (0.03125)node[scale=0.0625]{32};
\draw[black, thin] (-2.6666666666666665,-1.9444444444444444) circle (0.05555555555555555)node[scale=0.1111111111111111]{18};
\draw[black, thin] (2.6666666666666665,-1.9444444444444444) circle (0.05555555555555555)node[scale=0.1111111111111111]{18};
\draw[black, thin] (-1.5,-1.96875) circle (0.03125)node[scale=0.0625]{32};
\draw[black, thin] (1.5,-1.96875) circle (0.03125)node[scale=0.0625]{32};
\draw[black, thin] (-0.2857142857142857,-1.9795918367346939) circle (0.02040816326530612)node[scale=0.04081632653061224]{49};
\draw[black, thin] (0.2857142857142857,-1.9795918367346939) circle (0.02040816326530612)node[scale=0.04081632653061224]{49};
\draw[black, thin] (-0.8571428571428571,-1.9795918367346939) circle (0.02040816326530612)node[scale=0.04081632653061224]{49};
\draw[black, thin] (0.8571428571428571,-1.9795918367346939) circle (0.02040816326530612)node[scale=0.04081632653061224]{49};
\draw[black, thin] (-0.8,-1.98) circle (0.02)node[scale=0.04]{50};
\draw[black, thin] (0.8,-1.98) circle (0.02)node[scale=0.04]{50};
\draw[black, thin] (-2.8,-1.96) circle (0.04)node[scale=0.08]{25};
\draw[black, thin] (2.8,-1.96) circle (0.04)node[scale=0.08]{25};
\draw[black, thin] (-1.4285714285714286,-1.9795918367346939) circle (0.02040816326530612)node[scale=0.04081632653061224]{49};
\draw[black, thin] (1.4285714285714286,-1.9795918367346939) circle (0.02040816326530612)node[scale=0.04081632653061224]{49};
\draw[black, thin] (-0.3333333333333333,-1.9861111111111112) circle (0.013888888888888888)node[scale=0.027777777777777776]{72};
\draw[black, thin] (0.3333333333333333,-1.9861111111111112) circle (0.013888888888888888)node[scale=0.027777777777777776]{72};
\draw[black, thin] (-2.5,-1.96875) circle (0.03125)node[scale=0.0625]{32};
\draw[black, thin] (2.5,-1.96875) circle (0.03125)node[scale=0.0625]{32};
\draw[black, thin] (-1.6,-1.98) circle (0.02)node[scale=0.04]{50};
\draw[black, thin] (1.6,-1.98) circle (0.02)node[scale=0.04]{50};
\draw[black, thin] (-0.2222222222222222,-1.9876543209876543) circle (0.012345679012345678)node[scale=0.024691358024691357]{81};
\draw[black, thin] (0.2222222222222222,-1.9876543209876543) circle (0.012345679012345678)node[scale=0.024691358024691357]{81};
\draw[black, thin] (-1.0,-1.875) circle (0.125)node[scale=0.25]{8};
\draw[black, thin] (1.0,-1.875) circle (0.125)node[scale=0.25]{8};
\draw[black, thin] (-0.6666666666666666,-1.8888888888888888) circle (0.1111111111111111)node[scale=0.2222222222222222]{9};
\draw[black, thin] (0.6666666666666666,-1.8888888888888888) circle (0.1111111111111111)node[scale=0.2222222222222222]{9};
\draw[black, thick] (-3.2,-2.0) -- (3.2,-2.0);
\draw[black, thick] (3.2,2.0) -- (-3.2,2.0);
\draw[black, thick] (-2.0,-1.0) circle (1.0)node[scale=2.0]{1};
\draw[black, thick] (2.0,-1.0) circle (1.0)node[scale=2.0]{1};
\draw[black, thin] (0.0,-1.5) circle (0.5)node[scale=1.0]{2};
\draw[black, thin] (-1.1111111111111112,-1.9876543209876543) circle (0.012345679012345678)node[scale=0.024691358024691357]{81};
\draw[black, thin] (1.1111111111111112,-1.9876543209876543) circle (0.012345679012345678)node[scale=0.024691358024691357]{81};
\draw[black, thin] (-0.5714285714285714,-1.989795918367347) circle (0.01020408163265306)node[scale=0.02040816326530612]{98};
\draw[black, thin] (0.5714285714285714,-1.989795918367347) circle (0.01020408163265306)node[scale=0.02040816326530612]{98};
\draw[black, thin] (-2.4,-1.98) circle (0.02)node[scale=0.04]{50};
\draw[black, thin] (2.4,-1.98) circle (0.02)node[scale=0.04]{50};
\draw[black, thin] (-1.6666666666666667,-1.9861111111111112) circle (0.013888888888888888)node[scale=0.027777777777777776]{72};
\draw[black, thin] (1.6666666666666667,-1.9861111111111112) circle (0.013888888888888888)node[scale=0.027777777777777776]{72};
\draw[black, thin] (-2.5714285714285716,-1.9795918367346939) circle (0.02040816326530612)node[scale=0.04081632653061224]{49};
\draw[black, thin] (2.5714285714285716,-1.9795918367346939) circle (0.02040816326530612)node[scale=0.04081632653061224]{49};
\draw[black, thin] (-1.5555555555555556,-1.9876543209876543) circle (0.012345679012345678)node[scale=0.024691358024691357]{81};
\draw[black, thin] (1.5555555555555556,-1.9876543209876543) circle (0.012345679012345678)node[scale=0.024691358024691357]{81};
\draw[black, thin] (-1.1428571428571428,-1.989795918367347) circle (0.01020408163265306)node[scale=0.02040816326530612]{98};
\draw[black, thin] (1.1428571428571428,-1.989795918367347) circle (0.01020408163265306)node[scale=0.02040816326530612]{98};
\draw[black, thin] (-3.142857142857143,-1.9795918367346939) circle (0.02040816326530612)node[scale=0.04081632653061224]{49};
\draw[black, thin] (3.142857142857143,-1.9795918367346939) circle (0.02040816326530612)node[scale=0.04081632653061224]{49};
\draw[black, thin] (-2.3333333333333335,-1.9861111111111112) circle (0.013888888888888888)node[scale=0.027777777777777776]{72};
\draw[black, thin] (2.3333333333333335,-1.9861111111111112) circle (0.013888888888888888)node[scale=0.027777777777777776]{72};
\draw[black, thin] (-1.7142857142857142,-1.989795918367347) circle (0.01020408163265306)node[scale=0.02040816326530612]{98};
\draw[black, thin] (1.7142857142857142,-1.989795918367347) circle (0.01020408163265306)node[scale=0.02040816326530612]{98};
\draw[black, thin] (-2.4444444444444446,-1.9876543209876543) circle (0.012345679012345678)node[scale=0.024691358024691357]{81};
\draw[black, thin] (2.4444444444444446,-1.9876543209876543) circle (0.012345679012345678)node[scale=0.024691358024691357]{81};
\draw[black, thin] (-2.2857142857142856,-1.989795918367347) circle (0.01020408163265306)node[scale=0.02040816326530612]{98};
\draw[black, thin] (2.2857142857142856,-1.989795918367347) circle (0.01020408163265306)node[scale=0.02040816326530612]{98};
\draw[black, thin] (-2.888888888888889,-1.9876543209876543) circle (0.012345679012345678)node[scale=0.024691358024691357]{81};
\draw[black, thin] (2.888888888888889,-1.9876543209876543) circle (0.012345679012345678)node[scale=0.024691358024691357]{81};

\draw[black, dashed] (-2.0,2.5) -- (-2.0,-3.0);
\draw[black, dashed] (2.0,-3.0) -- (2.0,2.5);
\draw[black, dashed] (-1.0,-2.0) circle (1.0);
\draw[black, dashed] (1.0,-2.0) circle (1.0);
\end{tikzpicture}
\end{center}
\caption{Square Ford Circles and Duals}
\end{figure}
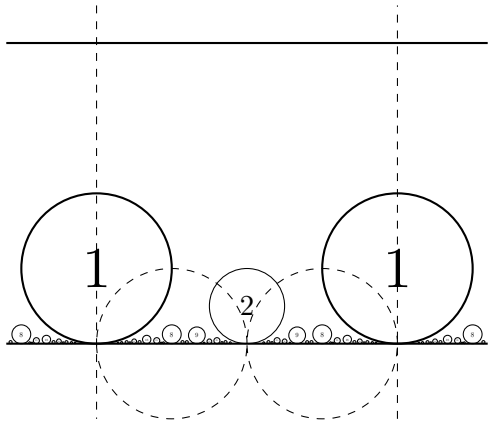

We can use \cite[Theorem 8.1]{FuchsStangeZhang} to determine the modulus for congruence restrictions in this packing. The symmetry group of the packing contains a congruence subgroup of level 2--see Figure \ref{SquareFordFig}. The symmetry group can be generated by the following M\"obius transformations: the translations $T_1=\begin{pmatrix} 1 & 2 \\ 0 & 1 \end{pmatrix}$, $T_2=\begin{pmatrix} 1 & 2i \\ 0 & 1 \end{pmatrix}$, and the inversion $S=\begin{pmatrix} 0 & -1 \\ 1 & 0 \end{pmatrix}$. Setting $H=\begin{pmatrix} 1 & 0 \\ 0 & -1\end{pmatrix}, \, R=\begin{pmatrix} 0 & 1 \\ 0 & 0\end{pmatrix}, \, L=\begin{pmatrix} 0 & 0 \\ 1 & 0\end{pmatrix} \in \mathfrak{sl}(2, \Z)$, we have
\begin{equation*}
\begin{split}
    & SHS^{-1}= -H, \quad S R S^{-1} = -L \quad S L S^{-1}=-R \\
    & T_2 L T_2^{-1} = 2i H + 4R + L, \quad T_2^{-1} H T_2= H+4i R, \quad S T_2 H T_2^{-1} S^{-1} = -H+4i L
\end{split}    
\end{equation*}
Therefore the image of the adjoint action contains $4\mathfrak{sl}(2, \Z[i])$. By \cite[Theorem 8.1]{FuchsStangeZhang}, the modulus for the congruence restrictions for this packing divides $8$. We now find the possible congruence restrictions mod 8.

\begin{proposition} \label{SquareMod}
A primitive integral square grid configuration of curvatures must be congruent to one of the following modulo 8:
\begin{center}
\, \hfill
\begin{tikzpicture}[scale=.5]
    \node (A) at (0,0)    {$1$};
    \node (B) at (0,2)    {$1$};
    \node (C) at (2,0)    {$1$};
    \node (D) at (2,2)    {$1$};
    \draw (A) -- (B);
    \draw (A) -- (C);
    \draw (B) -- (D);
    \draw (C) -- (D);
    \draw[dotted] (A) -- (-1,0);
    \draw[dotted] (A) -- (0,-1);
    \draw[dotted] (B) -- (-1,2);
    \draw[dotted] (B) -- (0,3);
    \draw[dotted] (C) -- (2,-1);
    \draw[dotted] (C) -- (3,0);
    \draw[dotted] (D) -- (3,2);
    \draw[dotted] (D) -- (2,3);
\end{tikzpicture} \hfill
\begin{tikzpicture}[scale=.5]
    \node (A) at (0,0)    {$5$};
    \node (B) at (0,2)    {$5$};
    \node (C) at (2,0)    {$5$};
    \node (D) at (2,2)    {$5$};
    \draw (A) -- (B);
    \draw (A) -- (C);
    \draw (B) -- (D);
    \draw (C) -- (D);
    \draw[dotted] (A) -- (-1,0);
    \draw[dotted] (A) -- (0,-1);
    \draw[dotted] (B) -- (-1,2);
    \draw[dotted] (B) -- (0,3);
    \draw[dotted] (C) -- (2,-1);
    \draw[dotted] (C) -- (3,0);
    \draw[dotted] (D) -- (3,2);
    \draw[dotted] (D) -- (2,3);
\end{tikzpicture} \hfill
\begin{tikzpicture}[scale=.5]
    \node (A) at (0,0)    {$3$};
    \node (B) at (0,2)    {$7$};
    \node (C) at (2,0)    {$7$};
    \node (D) at (2,2)    {$3$};
    \draw (A) -- (B);
    \draw (A) -- (C);
    \draw (B) -- (D);
    \draw (C) -- (D);
    \draw[dotted] (A) -- (-1,0);
    \draw[dotted] (A) -- (0,-1);
    \draw[dotted] (B) -- (-1,2);
    \draw[dotted] (B) -- (0,3);
    \draw[dotted] (C) -- (2,-1);
    \draw[dotted] (C) -- (3,0);
    \draw[dotted] (D) -- (3,2);
    \draw[dotted] (D) -- (2,3);
\end{tikzpicture} 
\hfill \,
\vspace{1cm}

\begin{tikzpicture}[scale=.5]
    \node (A) at (0,0)    {$1$};
    \node (E) at (0,-2)    {$0$};
    \node (I) at (0,-4)    {$1$};
    \node (M) at (0,-6)    {$4$};
    \node (B) at (2,0)    {$3$};
    \node (F) at (2,-2)    {$2$};
    \node (J) at (2,-4)    {$3$};
    \node (N) at (2,-6)    {$6$};
    \node (C) at (4,0)    {$7$};
    \node (G) at (4,-2)    {$6$};
    \node (K) at (4,-4)    {$7$};
    \node (O) at (4,-6)    {$2$};
    \node (D) at (6,0)    {$5$};
    \node (H) at (6,-2)    {$4$};
    \node (L) at (6,-4)    {$5$};
    \node (P) at (6,-6)    {$0$};

    \draw (A) -- (B);
    \draw (B) -- (C);
    \draw (C) -- (D);
    \draw (E) -- (F);
    \draw (F) -- (G);
    \draw (G) -- (H);
    \draw (I) -- (J);
    \draw (J) -- (K);
    \draw (K) -- (L);
    \draw (M) -- (N);
    \draw (N) -- (O);
    \draw (O) -- (P);

    \draw (A) -- (E);
    \draw (E) -- (I);
    \draw (I) -- (M);
    \draw (B) -- (F);
    \draw (F) -- (J);
    \draw (J) -- (N);
    \draw (C) -- (G);
    \draw (G) -- (K);
    \draw (K) -- (O);
    \draw (D) -- (H);
    \draw (H) -- (L);
    \draw (L) -- (P);
    
    \node (AA) at (8,0)    {$5$};
    \node (EE) at (8,-2)    {$4$};
    \node (II) at (8,-4)    {$5$};
    \node (MM) at (8,-6)    {$0$};
    \node (BB) at (10,0)    {$7$};
    \node (FF) at (10,-2)    {$6$};
    \node (JJ) at (10,-4)    {$7$};
    \node (NN) at (10,-6)    {$2$};
    \node (CC) at (12,0)    {$3$};
    \node (GG) at (12,-2)    {$2$};
    \node (KK) at (12,-4)    {$3$};
    \node (OO) at (12,-6)    {$6$};
    \node (DD) at (14,0)    {$1$};
    \node (HH) at (14,-2)    {$0$};
    \node (LL) at (14,-4)    {$1$};
    \node (PP) at (14,-6)    {$4$};

    \draw (AA) -- (BB);
    \draw (BB) -- (CC);
    \draw (CC) -- (DD);
    \draw (EE) -- (FF);
    \draw (FF) -- (GG);
    \draw (GG) -- (HH);
    \draw (II) -- (JJ);
    \draw (JJ) -- (KK);
    \draw (KK) -- (LL);
    \draw (MM) -- (NN);
    \draw (NN) -- (OO);
    \draw (OO) -- (PP);

    \draw (AA) -- (EE);
    \draw (EE) -- (II);
    \draw (II) -- (MM);
    \draw (BB) -- (FF);
    \draw (FF) -- (JJ);
    \draw (JJ) -- (NN);
    \draw (CC) -- (GG);
    \draw (GG) -- (KK);
    \draw (KK) -- (OO);
    \draw (DD) -- (HH);
    \draw (HH) -- (LL);
    \draw (LL) -- (PP);
    
    \draw[dotted] (A) -- (-1,0);
    \draw[dotted] (E) -- (-1,-2);
    \draw[dotted] (I) -- (-1,-4);
    \draw[dotted] (M) -- (-1,-6);
    \draw[dotted] (DD) -- (15,0);
    \draw[dotted] (HH) -- (15,-2);
    \draw[dotted] (LL) -- (15,-4);
    \draw[dotted] (PP) -- (15,-6);
    \draw (D) -- (AA);
    \draw (H) -- (EE);
    \draw (L) -- (II);
    \draw (P) -- (MM);
    \draw[dotted] (A) -- (0,1);
    \draw[dotted] (B) -- (2,1);
    \draw[dotted] (C) -- (4,1);
    \draw[dotted] (D) -- (6,1);
    \draw[dotted] (AA) -- (8,1);
    \draw[dotted] (BB) -- (10,1);
    \draw[dotted] (CC) -- (12,1);
    \draw[dotted] (DD) -- (14,1);
    \draw[dotted] (M) -- (0,-7);
    \draw[dotted] (N) -- (2,-7);
    \draw[dotted] (O) -- (4,-7);
    \draw[dotted] (P) -- (6,-7);
    \draw[dotted] (MM) -- (8,-7);
    \draw[dotted] (NN) -- (10,-7);
    \draw[dotted] (OO) -- (12,-7);
    \draw[dotted] (PP) -- (14,-7);
\end{tikzpicture}
\end{center}

with the pattern repeating periodically. 
\end{proposition}

\begin{lemma} \label{SquareModLemma1}
In a primitive integral square grid configuration of circles, the curvatures mod 2 have one of the following periodic patterns:

\, \hfill
\begin{tikzpicture}[scale=.5]
    \node (A) at (0,0)    {$1$};
    \node (B) at (0,2)    {$1$};
    \node (C) at (2,0)    {$1$};
    \node (D) at (2,2)    {$1$};
    \draw (A) -- (B);
    \draw (A) -- (C);
    \draw (B) -- (D);
    \draw (C) -- (D);
    \draw[dotted] (A) -- (-1,0);
    \draw[dotted] (A) -- (0,-1);
    \draw[dotted] (B) -- (-1,2);
    \draw[dotted] (B) -- (0,3);
    \draw[dotted] (C) -- (2,-1);
    \draw[dotted] (C) -- (3,0);
    \draw[dotted] (D) -- (3,2);
    \draw[dotted] (D) -- (2,3);
\end{tikzpicture} \hfill
\begin{tikzpicture}[scale=.5]
    \node (A) at (0,0)    {$1$};
    \node (B) at (0,2)    {$0$};
    \node (C) at (2,0)    {$1$};
    \node (D) at (2,2)    {$0$};
    \draw (A) -- (B);
    \draw (A) -- (C);
    \draw (B) -- (D);
    \draw (C) -- (D);
    \draw[dotted] (A) -- (-1,0);
    \draw[dotted] (A) -- (0,-1);
    \draw[dotted] (B) -- (-1,2);
    \draw[dotted] (B) -- (0,3);
    \draw[dotted] (C) -- (2,-1);
    \draw[dotted] (C) -- (3,0);
    \draw[dotted] (D) -- (3,2);
    \draw[dotted] (D) -- (2,3);
\end{tikzpicture} \hfill
\begin{tikzpicture}[scale=.5]
    \node (A) at (0,0)    {$0$};
    \node (B) at (0,2)    {$1$};
    \node (C) at (2,0)    {$1$};
    \node (D) at (2,2)    {$0$};
    \draw (A) -- (B);
    \draw (A) -- (C);
    \draw (B) -- (D);
    \draw (C) -- (D);
    \draw[dotted] (A) -- (-1,0);
    \draw[dotted] (A) -- (0,-1);
    \draw[dotted] (B) -- (-1,2);
    \draw[dotted] (B) -- (0,3);
    \draw[dotted] (C) -- (2,-1);
    \draw[dotted] (C) -- (3,0);
    \draw[dotted] (D) -- (3,2);
    \draw[dotted] (D) -- (2,3);
\end{tikzpicture} 
\hfill \,

\end{lemma}

\begin{proof}
    Consider an arbitrary set of four circles in the configuration, with curvatures $j$, $k$, $l$, $m$, arranged as in the quadratic relation \eqref{SquareQuadratic}. Equation \eqref{SquareQuadratic} reduces to $j^2+k^2+l^2+m^2\equiv 0 \bmod 2$ so either all of the curvatures are odd or two are odd and two are even. They cannot all be even because the packing is primitive. 
    
    Moreover, by the linear equations, choosing four circles that satisfy the quadratic relation uniquely determines the configuration. If $j$, $k$, $l$, $m$ are all odd, then all curvatures in the configuration are odd, so we are in the first case. If $j \equiv k \not\equiv l \equiv m$ or $j\equiv m \not \equiv k \equiv l \bmod 2$, then we are in the second case. If $j\equiv l \not \equiv k \equiv m \bmod 2$, then we are in the third case.
\end{proof}

\begin{lemma} \label{SquareModLemma2}
Suppose that $k$ is an odd curvature, and $j$ and $l$ are curvatures of adjacent circles around a square in a primitive integral square grid configuration. Then $v_2((k+j)(k+l))=1$ or $2$. 
\end{lemma}

\begin{proof}
Equation \eqref{SquareQuadratic} can be rewritten as $(m-j-k-3l)^2=8jk+8jl+8kl-8k^2$. Since the left side is a square integer, it follows that $jk+jl+kl-k^2=2n^2$ for some integer $n$. This can be rewritten as $(k+j)(k+l)=2(k^2+n^2)$. Since $k$ is odd, $k^2+n^2 \equiv 1$ or $2 \bmod 4$. Thus $v_2((k+j)(k+l))=v_2(2(k^2+n^2))=1$ or $2$.
\end{proof}

Lemma \ref{SquareModLemma2} rules out the third configuration mod 2 shown in Lemma \ref{SquareModLemma1}. It also implies that the first configuration can only lift to one of the following mod 4:

\, \hfill
\begin{tikzpicture}[scale=.5]
    \node (A) at (0,0)    {$1$};
    \node (B) at (0,2)    {$1$};
    \node (C) at (2,0)    {$1$};
    \node (D) at (2,2)    {$1$};
    \draw (A) -- (B);
    \draw (A) -- (C);
    \draw (B) -- (D);
    \draw (C) -- (D);
    \draw[dotted] (A) -- (-1,0);
    \draw[dotted] (A) -- (0,-1);
    \draw[dotted] (B) -- (-1,2);
    \draw[dotted] (B) -- (0,3);
    \draw[dotted] (C) -- (2,-1);
    \draw[dotted] (C) -- (3,0);
    \draw[dotted] (D) -- (3,2);
    \draw[dotted] (D) -- (2,3);
\end{tikzpicture} \hfill
\begin{tikzpicture}[scale=.5]
    \node (A) at (0,0)    {$3$};
    \node (B) at (0,2)    {$3$};
    \node (C) at (2,0)    {$3$};
    \node (D) at (2,2)    {$3$};
    \draw (A) -- (B);
    \draw (A) -- (C);
    \draw (B) -- (D);
    \draw (C) -- (D);
    \draw[dotted] (A) -- (-1,0);
    \draw[dotted] (A) -- (0,-1);
    \draw[dotted] (B) -- (-1,2);
    \draw[dotted] (B) -- (0,3);
    \draw[dotted] (C) -- (2,-1);
    \draw[dotted] (C) -- (3,0);
    \draw[dotted] (D) -- (3,2);
    \draw[dotted] (D) -- (2,3);
\end{tikzpicture} 
\hfill \,

\begin{lemma} \label{SquareModLemma3}
    Let $j$, $k$ be curvatures of tangent circles in a primitive integral square grid configuration. Then $j+k \not\equiv 3, 6, 7 \bmod 8$. 
\end{lemma}
\begin{proof}
    As in the proof of the previous lemma, we have $(k+j)(k+l)=2(k^2+n^2)$. Assume that $j+k \equiv 3, 6, 7 \bmod 8$, then the odd part of $j+k$ is  $3 \bmod 4$. So, there exists an odd prime number $p \equiv 3 \bmod 4$ such that $p\mid j+k$ and $v_p(j+k)$ is odd. Since $p$ divides twice the sum of two squares, $v_p(k^2+n^2)$ must be even and $p\mid k,n$. Therefore, $p\mid j,k,l$, and by Equations \eqref{SquareFaceDiagonal}, \eqref{SquareT}, \eqref{SquareQuadratic}, all the curvatures in the configuration share a common factor of $p$, which means the configuration is not primitive. 
\end{proof}

Lemma \ref{SquareModLemma3} implies that the first configuration mod 2 shown in Lemma 1 must lift to one of the following mod 8:

\, \hfill
\begin{tikzpicture}[scale=.5]
    \node (A) at (0,0)    {$1$};
    \node (B) at (0,2)    {$1$};
    \node (C) at (2,0)    {$1$};
    \node (D) at (2,2)    {$1$};
    \draw (A) -- (B);
    \draw (A) -- (C);
    \draw (B) -- (D);
    \draw (C) -- (D);
    \draw[dotted] (A) -- (-1,0);
    \draw[dotted] (A) -- (0,-1);
    \draw[dotted] (B) -- (-1,2);
    \draw[dotted] (B) -- (0,3);
    \draw[dotted] (C) -- (2,-1);
    \draw[dotted] (C) -- (3,0);
    \draw[dotted] (D) -- (3,2);
    \draw[dotted] (D) -- (2,3);
\end{tikzpicture} \hfill
\begin{tikzpicture}[scale=.5]
    \node (A) at (0,0)    {$5$};
    \node (B) at (0,2)    {$5$};
    \node (C) at (2,0)    {$5$};
    \node (D) at (2,2)    {$5$};
    \draw (A) -- (B);
    \draw (A) -- (C);
    \draw (B) -- (D);
    \draw (C) -- (D);
    \draw[dotted] (A) -- (-1,0);
    \draw[dotted] (A) -- (0,-1);
    \draw[dotted] (B) -- (-1,2);
    \draw[dotted] (B) -- (0,3);
    \draw[dotted] (C) -- (2,-1);
    \draw[dotted] (C) -- (3,0);
    \draw[dotted] (D) -- (3,2);
    \draw[dotted] (D) -- (2,3);
\end{tikzpicture} \hfill
\begin{tikzpicture}[scale=.5]
    \node (A) at (0,0)    {$7$};
    \node (B) at (0,2)    {$3$};
    \node (C) at (2,0)    {$3$};
    \node (D) at (2,2)    {$7$};
    \draw (A) -- (B);
    \draw (A) -- (C);
    \draw (B) -- (D);
    \draw (C) -- (D);
    \draw[dotted] (A) -- (-1,0);
    \draw[dotted] (A) -- (0,-1);
    \draw[dotted] (B) -- (-1,2);
    \draw[dotted] (B) -- (0,3);
    \draw[dotted] (C) -- (2,-1);
    \draw[dotted] (C) -- (3,0);
    \draw[dotted] (D) -- (3,2);
    \draw[dotted] (D) -- (2,3);
\end{tikzpicture} 
\hfill \,

The second configuration in Lemma \ref{SquareModLemma1} is more complicated, but in fact lifts to mod 8 in a unique way. 

\begin{lemma}
Any primitive integral square grid configuration containing both odd and even curvatures is congruent to the following mod 8:
\begin{center}
\begin{tikzpicture}[scale=.5]
    \node (A) at (0,0)    {$1$};
    \node (E) at (0,-2)    {$0$};
    \node (I) at (0,-4)    {$1$};
    \node (M) at (0,-6)    {$4$};
    \node (B) at (2,0)    {$3$};
    \node (F) at (2,-2)    {$2$};
    \node (J) at (2,-4)    {$3$};
    \node (N) at (2,-6)    {$6$};
    \node (C) at (4,0)    {$7$};
    \node (G) at (4,-2)    {$6$};
    \node (K) at (4,-4)    {$7$};
    \node (O) at (4,-6)    {$2$};
    \node (D) at (6,0)    {$5$};
    \node (H) at (6,-2)    {$4$};
    \node (L) at (6,-4)    {$5$};
    \node (P) at (6,-6)    {$0$};

    \draw (A) -- (B);
    \draw (B) -- (C);
    \draw (C) -- (D);
    \draw (E) -- (F);
    \draw (F) -- (G);
    \draw (G) -- (H);
    \draw (I) -- (J);
    \draw (J) -- (K);
    \draw (K) -- (L);
    \draw (M) -- (N);
    \draw (N) -- (O);
    \draw (O) -- (P);

    \draw (A) -- (E);
    \draw (E) -- (I);
    \draw (I) -- (M);
    \draw (B) -- (F);
    \draw (F) -- (J);
    \draw (J) -- (N);
    \draw (C) -- (G);
    \draw (G) -- (K);
    \draw (K) -- (O);
    \draw (D) -- (H);
    \draw (H) -- (L);
    \draw (L) -- (P);
    
    \node (AA) at (8,0)    {$5$};
    \node (EE) at (8,-2)    {$4$};
    \node (II) at (8,-4)    {$5$};
    \node (MM) at (8,-6)    {$0$};
    \node (BB) at (10,0)    {$7$};
    \node (FF) at (10,-2)    {$6$};
    \node (JJ) at (10,-4)    {$7$};
    \node (NN) at (10,-6)    {$2$};
    \node (CC) at (12,0)    {$3$};
    \node (GG) at (12,-2)    {$2$};
    \node (KK) at (12,-4)    {$3$};
    \node (OO) at (12,-6)    {$6$};
    \node (DD) at (14,0)    {$1$};
    \node (HH) at (14,-2)    {$0$};
    \node (LL) at (14,-4)    {$1$};
    \node (PP) at (14,-6)    {$4$};

    \draw (AA) -- (BB);
    \draw (BB) -- (CC);
    \draw (CC) -- (DD);
    \draw (EE) -- (FF);
    \draw (FF) -- (GG);
    \draw (GG) -- (HH);
    \draw (II) -- (JJ);
    \draw (JJ) -- (KK);
    \draw (KK) -- (LL);
    \draw (MM) -- (NN);
    \draw (NN) -- (OO);
    \draw (OO) -- (PP);

    \draw (AA) -- (EE);
    \draw (EE) -- (II);
    \draw (II) -- (MM);
    \draw (BB) -- (FF);
    \draw (FF) -- (JJ);
    \draw (JJ) -- (NN);
    \draw (CC) -- (GG);
    \draw (GG) -- (KK);
    \draw (KK) -- (OO);
    \draw (DD) -- (HH);
    \draw (HH) -- (LL);
    \draw (LL) -- (PP);
    
    \draw[dotted] (A) -- (-1,0);
    \draw[dotted] (E) -- (-1,-2);
    \draw[dotted] (I) -- (-1,-4);
    \draw[dotted] (M) -- (-1,-6);
    \draw[dotted] (DD) -- (15,0);
    \draw[dotted] (HH) -- (15,-2);
    \draw[dotted] (LL) -- (15,-4);
    \draw[dotted] (PP) -- (15,-6);
    \draw (D) -- (AA);
    \draw (H) -- (EE);
    \draw (L) -- (II);
    \draw (P) -- (MM);
    \draw[dotted] (A) -- (0,1);
    \draw[dotted] (B) -- (2,1);
    \draw[dotted] (C) -- (4,1);
    \draw[dotted] (D) -- (6,1);
    \draw[dotted] (AA) -- (8,1);
    \draw[dotted] (BB) -- (10,1);
    \draw[dotted] (CC) -- (12,1);
    \draw[dotted] (DD) -- (14,1);
    \draw[dotted] (M) -- (0,-7);
    \draw[dotted] (N) -- (2,-7);
    \draw[dotted] (O) -- (4,-7);
    \draw[dotted] (P) -- (6,-7);
    \draw[dotted] (MM) -- (8,-7);
    \draw[dotted] (NN) -- (10,-7);
    \draw[dotted] (OO) -- (12,-7);
    \draw[dotted] (PP) -- (14,-7);
\end{tikzpicture}
\end{center}
with the pattern repeating periodically. 
\end{lemma}
\begin{proof}
Note that such a configuration must be a lifting of the second configuration in Lemma \ref{SquareModLemma1}. Consider an odd curvature $k$ in this configuration. By Lemmas \ref{SquareModLemma2} and \ref{SquareModLemma3}, two adjacent odd curvatures cannot sum to $0$ or $6$ mod 8, so the possible odd curvatures adjacent to $k$ are $2-k$ and $4-k \bmod 8$. The possible even curvatures adjacent to $k$ are $1-k$ and $5-k$. The only way for Equation \eqref{SquareT} to be satisfied is if all four possible curvatures appear adjacent to $k$. These five curvatures then determine all curvatures in the configuration by Equations \eqref{SquareFaceDiagonal} and \eqref{SquareT}. Since the configuration shown above satisfies Equations \eqref{SquareFaceDiagonal}, \eqref{SquareT}, and \eqref{SquareQuadratic} mod 8, and contains every possible odd curvature with four distinct adjacent curvatures, it is the unique way to lift the second configuration in Lemma \ref{SquareModLemma1}.
\end{proof} 

\begin{proposition}
    In a primitive integral square packing, every square grid configuration of circles belongs to the same type modulo 8. 
\end{proposition}
This proposition is proven similarly to Proposition \ref{OctGroupMod}. It suffices to check, in Proposition \ref{SquareMod}, that each configuration is uniquely determined by one square of curvatures.

We will refer to the four types of square grid packings mod 8 as type $(1)$, type $(5)$, type $(3, 7)$, and full type. Section \ref{Data} below contains examples of all four packing types mod 8. 

Next we parametrize the curvatures of all circles tangent to a fixed circle in a square packing, beginning with the Ford circles.

\begin{proposition} \label{SquareFord}
The generalized Ford circles are parametrized by $x, y \in \Z$ with $y \geq 0$ and $\gcd(x,y)=1$. Their inversive coordinates are as follows:
\begin{align*}    
& c_\alpha(x,y)=(4x^2,y^2,2xy,1)\hspace{3pt} \text{ if } xy \text{ is odd} \\
& c_\beta(x,y)=(8x^2, 2y^2, 4xy, 1) \hspace{13pt} \text{ if } xy \text{ is even} 
\end{align*}
Each circle is tangent to the real axis at $\frac{2x}{y}$.
\end{proposition} 

\begin{proposition} \label{SquareParametrization}
    Suppose a square packing contains a circle of curvature $a$, and the four circles around it in a square grid have curvatures $b$, $c$, $d$, $e$ as shown. 

\begin{center}
\begin{tikzpicture}[scale=.3]
    \node (A) at (0,0)    {c};
    \node (C) at (0,-2)    {a};
    \node (B) at (-2,-2)    {b};
    \node (D) at (2,-2)    {d};
    \node (E) at (0,-4)    {e};
    \draw (A) -- (C);
    \draw (B) -- (C);
    \draw (D) -- (C);
    \draw (C) -- (E);
\end{tikzpicture}
\end{center}

Then all circles tangent to the circle of curvature $a$ are parametrized by $x, y \in \Z$ with $y \geq 0$, $\gcd(x,y)=1$. Their curvatures are as follows:
\begin{align*}
& Q_\alpha(x,y)-a=\left(\frac{a+b}{2}\right)x^2 +  \left(\frac{c-e}{2}\right)xy + \left(\frac{a+d}{2}\right)y^2 -a\hspace{16pt} \text{ if } xy \text{ is odd} \\
& Q_\beta(x,y)-a=\left(a+b\right)x^2 +  \left(c-e\right)xy + (a+d)y^2 -a \hspace{3pt} \text{ if } xy \text{ is even} 
\end{align*}
\end{proposition}

Propositions \ref{SquareFord} and \ref{SquareParametrization} can be proved similarly to those in previous sections. Using Equations \eqref{SquareFaceDiagonal}, \eqref{SquareT}, and \eqref{SquareQuadratic}, we find that $Q_\alpha$ has discriminant $-4a^2$ and $Q_\beta$ has discriminant $-16a^2$.

We now define quadratic invariants $\chi_2$ for square packings. The quadratic forms $Q_\alpha$ and $Q_\beta$ both represent the integers $a+b$, $a+c$, $a+d$, and $a+e$ up to factors of 2, and at least one of these forms represents odd integers, so because $\gcd(a,b,c,d,e)=1$, we can conclude that $Q_\alpha(x,y)$ or $Q_\beta(x,y)$ represents an integer $\rho$ coprime to $a$. We may take $y\geq 0$, $\gcd(x,y)=1$, and $xy$ odd if $Q_\alpha$ is used or even if $Q_\beta$ is used. 

\begin{proposition} \label{SquareNode}
Let $a$ be a curvature in a primitive integral square packing, and suppose that all curvatures of circles tangent to $a$ are represented by the quadratic forms $Q_\alpha$ and $Q_\beta$ as in the previous proposition. Let $\rho_1$, $\rho_2$ be two integers coprime to $a$ which are represented by either $Q_{\alpha}$ or $Q_{\beta}$. If $a \not \equiv 3,5 \bmod 8$, then 
\begin{equation*}
\res{\rho_1}{a} = \res{\rho_2}{a}
\end{equation*}
If $a \equiv 3, 5 \bmod 8$, then the equality holds if and only if $\rho_1$, $\rho_2$ are represented by the same quadratic form; if $\rho_1$ is represented by $Q_\alpha$ and $\rho_2$ is represented by $Q_\beta$, then 
\begin{equation*}
\res{\rho_1}{a} = -\res{\rho_2}{a}
\end{equation*}
\end{proposition}

The proof is the same as that of Proposition \ref{OctNode} above, using Lemma \ref{QuadFormKronecker}. If $a$ is odd, then both discriminants $-4a^2$ and $-16a^2$ are divisible by $a$. If $a$ is even we may take $Q_\beta$ to represent the adjacent odd curvatures and the discriminant $-16a^2$ is divisible by $a$ and a sufficient power of $2$. 

We now define the quadratic invariant $\chi_2$ for square packings of type $(1)$ and of full type, where the invariant will be well-defined across the entire packing. For a circle $\C$ of curvature $a$, we set
\begin{align*}
\chi_2(\mathcal{C}) = \begin{cases}
    \left(\frac{\rho}{a}\right) &\text{if } a \equiv 0,1,4 \bmod 8;\\
    -\left(\frac{\rho}{a}\right) &\text{if } a \equiv 2, 6, 7 \bmod 8;\\
    (-1)^{\rho} \left(\frac{\rho}{a}\right) &\text{if } a \equiv 3 \bmod 8;\\
    (-1)^{\rho+1} \left(\frac{\rho}{a}\right) &\text{if } a \equiv 5 \bmod 8;
    \end{cases}
\end{align*}
Where $\rho$ is any curvature represented by $Q_\alpha$ or $Q_\beta$, satisfying $\gcd(\rho, a)=1$. Proposition \ref{SquareNode} guarantees that this expression is independent of the choice of $\rho$--if $a \equiv 3, 5 \bmod 8$ then one of the quadratic forms represents even curvaturess and the other represents odd curvatures, so the factor of $(-1)^\rho$ makes $\chi_2(\C)$ well-defined.  

\begin{proposition} \label{SquareEdge}
Suppose that $\C_a$ and $\C_b$ are tangent circles with coprime curvatures $a$, $b$ in a primitive integral square packing of type $(1)$ or of full type. Then 
\begin{equation*}
\chi_2(\C_a)=\chi_2(\C_b)
\end{equation*}
\end{proposition}
\begin{proof}
We will use $\rho=a+b$ to compute $\chi_2(\C_a)$ and $\chi_2(\C_b)$. Assume without loss of generality that $a$ is odd. If $a\equiv 1 \bmod 8$ then $b\equiv 0, 1, 3, 4 \bmod 8$ and 
\begin{equation*}
\chi_2(\C_a)=\res{a+b}{a} = \res{b}{a} = \res{a}{b} = \res{a+b}{b}=\chi_2(\C_b)
\end{equation*}
If $a\equiv 5 \bmod 8$, then $b \equiv 0, 4, 5, 7 \bmod 8$, and 
\begin{equation*}
\chi_2(\C_a)=(-1)^{a+b+1}\res{a+b}{a} = (-1)^{a+b+1} \res{b}{a} = (-1)^{a+b+1} \res{a}{b} = (-1)^{a+b+1} \res{a+b}{b}=\chi_2(\C_b)
\end{equation*}
For $a \equiv 3$, $b \equiv 7 \bmod 8$,
\begin{equation*}
\chi_2(\C_a)=\res{a+b}{a} = \res{b}{a} = -\res{a}{b} = -\res{a+b}{b}=\chi_2(\C_b)
\end{equation*}
For $a\equiv 3, 7 \bmod 8$, $b \equiv 2 \bmod 8$,
\begin{equation*}
\chi_2(\C_a)=-\res{a+b}{a} = -\res{b}{a} = -\res{a}{b} = -\res{a+b}{b}=\chi_2(\C_b)
\end{equation*}
The equality $\res{a}{b} = \res{a+b}{b}$ is justified because $\res{a}{2}=\res{a+b}{2}$ and $\res{a}{b/2}=\res{a+b}{b/2}$. For $a\equiv 3, 7 \bmod 8$, $b \equiv 6 \bmod 8$,
\begin{equation*}
\chi_2(\C_a)=-\res{a+b}{a} = -\res{b}{a} = \res{a}{b} = -\res{a+b}{b}=\chi_2(\C_b)
\end{equation*}
The equality $\res{a}{b} = -\res{a+b}{b}$ is justified because $\res{a}{2}=-\res{a+b}{2}$ and $\res{a}{b/2}=\res{a+b}{b/2}$. This completes the verification for all possible coprime pairs of adjacent curvatures in type $(1)$ or full type.
\end{proof}

To show that $\chi_2(\mathcal{C})$ is identical for all circles in the packing, we must find a path of coprime curvatures between any two circles. 

\begin{lemma} \label{SquareSimultaneous}
Suppose that $a, b, c, d, e, f$ are curvatures in a square packing, with tangency relations as shown. 
\begin{center}
\begin{tikzpicture}[scale=.3]
    \node (A) at (0,0)    {a};
    \node (C) at (2,0)    {c};
    \node (B) at (0,-2)    {b};
    \node (D) at (2,-2)    {d};
    \node (E) at (-2,0)    {e};
    \node (F) at (-2,-2)    {f};
    \draw (A) -- (B);
    \draw (A) -- (C);
    \draw (A) -- (E);
    \draw (B) -- (D);
    \draw (B) -- (F);
    \draw (D) -- (C);
    \draw (F) -- (E);
\end{tikzpicture}
\end{center}
Then all pairs of curvatures which form a square of tangencies with $a,b$ in the packing are parametrized by the following pair of formulas:
\begin{equation*}
    a \left(2n^2-2n\right)+b \left(2n^2-2n\right)+c (1-n)+e n, \quad a \left(2n^2-2n\right)+b \left(2n^2-2n\right)+d (1-n)+f n
\end{equation*} 
\end{lemma}
The proof is similar to that of Proposition \ref{CubeSimultaneous}. 

\begin{lemma} \label{SquareInsert}
Given a pair of tangent circles of curvatures $a$, $b$ in a primitive integral square packing, there exist circles of curvatures $k$, $\ell$ forming a square of tangencies with $a$, $b$, such that $\gcd(b, k)\mid 4$, $\gcd(\ell,a)\mid 4$ and $\gcd(k,\ell)\mid 8$. 
\end{lemma}

\begin{proof}
The proof is similar to that of Proposition \ref{CubeInsert}, so we will just sketch it briefly. For $a, b, c, d, e, f$ as in the previous lemma, the first quadratic polynomial represents the integers $c, e, 4a+4b+2c-e$. We have $\gcd(a, b, c, e)=1$ because the packing is primitive, so $\gcd(a, c, e, 4a+4b+2c-e)\mid 4$, and the first quadratic polynomial represents integers $\ell$ such that $\gcd(\ell, a)\mid 4$ by the Chinese remainder theorem. Similarly, the second quadratic polynomial represents integers $k$ such that $\gcd(k, b)\mid 4$. Since the difference between the two parametrized curvatures in Lemma \ref{CubeSimultaneous} is always $a-b$, the greatest common divisor of these curvatures always divides $a-b$. Because $\gcd(a-b, c,d,e,f,4a+4b+2c-e)\mid 8$,  we can find $k$, $\ell$ such that $\gcd(k, \ell)\mid 8$.
\end{proof}

\begin{corollary} \label{SquarePath}
Let $\mathcal{C},\mathcal{C}' \in \mathcal{A}$ be two circles in a primitive square circle packing. Then there exists a path of circles $\mathcal{C}_1,\mathcal{C}_2,...,\mathcal{C}_k$ such that \begin{enumerate}
    \item $\mathcal{C}_1=\mathcal{C}$ and $\mathcal{C}_k=\mathcal{C}'$;
    \item $\mathcal{C}_i$ is tangent to $\mathcal{C}_{i+1}$ for all $1 \leq i \leq k-1$;
    \item The curvatures of $\mathcal{C}_i$ and $\mathcal{C}_{i+1}$ are coprime for all $1 \leq i \leq k-1$.
\end{enumerate}
\end{corollary}
The proof is the same as that of Proposition \ref{CubePath}, using Lemma \ref{SquareInsert} to insert pairs of circles between two circles in the path if the curvatures are not coprime. If two tangent circles have even curvatures, Proposition \ref{SquareMod} implies that the two circles inserted between them will have odd curvatures, ensuring coprimality.  

Proposition \ref{SquareEdge} and Corollary \ref{SquarePath} imply the following result. 

\begin{proposition} \label{SquareChi}
The value of $\chi_2$ is constant across all circles in a fixed primitive square circle packing $\mathcal{A}$ of type $(1)$ or full type.
\end{proposition}

Based on this Proposition, we will refer to $\chi_2(\mathcal{A})$ as a quadratic invariant of the packing. This is used to prove the main result of this section.

\begin{theorem} \label{SquareObstruction}
In a primitive integral square packing $\mathcal{A}$ of type $(1)$ or full type with $\chi_2(\mathcal{A})=-1$, no integers of the form $n^2$ appear as curvatures. In particular, the Local-Global Conjecture is false for these packings.
\end{theorem}
\begin{proof}
Suppose that some circle in $\mathcal{A}$ has curvature $n^2$. Choose a tangent circle of curvature $a \equiv 0, 1, 3, 4, 5 \bmod 8$ with $\gcd(a, n^2)=1$. Then 
\begin{equation*}
\chi_2(\mathcal{A}) = \res{n^2+a}{a} = \res{n^2}{a} =1
\end{equation*}
which is a contradiction.
\end{proof}

In full type, we also see partial obstructions on curvatures of the form $2n^2$. If $\chi_2(\A)=-1$, then there are no curvatures of the form $2n^2$ tangent to curvatures $a\equiv 1, 3 \bmod 8$. If $\chi_2(\A)=1$, then there are no curvatures of the form $2n^2$ tangent to curvatures $a\equiv 5,7 \bmod 8$. But these do not give obstructions to the Local-Global Conjecture across all circles in the packing.

\begin{figure}[h]
\begin{center} \, \hfill
\begin{tikzpicture}[scale=.5]
    \node (A) at (0,0)    {5};
    \node (B) at (2,0)    {5};
    \node (C) at (4,0)    {5};
    \node (D) at (6,0)    {5};
    \node (E) at (0,-2)    {5};
    \node (F) at (2,-2)    {5};
    \node (G) at (4,-2)    {5};
    \node (H) at (6,-2)    {5};
    \node (I) at (0,-4)    {5};
    \node (J) at (2,-4)    {5};
    \node (K) at (4,-4)    {5};
    \node (L) at (6,-4)    {5};
    \node (M) at (0,-6)    {5};
    \node (N) at (2,-6)    {5};
    \node (O) at (4,-6)    {5};
    \node (P) at (6,-6)    {5};

    \draw[<->] (A) -- (B);
    \draw[<->] (B) -- (C);
    \draw[<->] (C) -- (D);
    \draw[<->] (E) -- (F);
    \draw[<->] (F) -- (G);
    \draw[<->] (G) -- (H);
    \draw[<->] (I) -- (J);
    \draw[<->] (J) -- (K);
    \draw[<->] (K) -- (L);
    \draw[<->] (M) -- (N);
    \draw[<->] (N) -- (O);
    \draw[<->] (O) -- (P);

    \draw (A) -- (E);
    \draw (E) -- (I);
    \draw (I) -- (M);
    \draw (B) -- (F);
    \draw (F) -- (J);
    \draw (J) -- (N);
    \draw (C) -- (G);
    \draw (G) -- (K);
    \draw (K) -- (O);
    \draw (D) -- (H);
    \draw (H) -- (L);
    \draw (L) -- (P);
    
    \draw[dotted] (A) -- (-1,0);
    \draw[dotted] (E) -- (-1,-2);
    \draw[dotted] (I) -- (-1,-4);
    \draw[dotted] (M) -- (-1,-6);
    \draw[dotted] (D) -- (7,0);
    \draw[dotted] (H) -- (7,-2);
    \draw[dotted] (L) -- (7,-4);
    \draw[dotted] (P) -- (7,-6);
    \draw[dotted] (A) -- (0,1);
    \draw[dotted] (B) -- (2,1);
    \draw[dotted] (C) -- (4,1);
    \draw[dotted] (D) -- (6,1);
    \draw[dotted] (M) -- (0,-7);
    \draw[dotted] (N) -- (2,-7);
    \draw[dotted] (O) -- (4,-7);
    \draw[dotted] (P) -- (6,-7);
\end{tikzpicture} \hfill
\begin{tikzpicture}[scale=.5]
    \node (A) at (0,0)    {3};
    \node (B) at (2,0)    {7};
    \node (C) at (4,0)    {3};
    \node (D) at (6,0)    {7};
    \node (E) at (0,-2)    {7};
    \node (F) at (2,-2)    {3};
    \node (G) at (4,-2)    {7};
    \node (H) at (6,-2)    {3};
    \node (I) at (0,-4)    {3};
    \node (J) at (2,-4)    {7};
    \node (K) at (4,-4)    {3};
    \node (L) at (6,-4)    {7};
    \node (M) at (0,-6)    {7};
    \node (N) at (2,-6)    {3};
    \node (O) at (4,-6)    {7};
    \node (P) at (6,-6)    {3};

    \draw[->] (A) -- (B);
    \draw[<-] (B) -- (C);
    \draw[->] (C) -- (D);
    \draw[->] (E) -- (F);
    \draw[<-] (F) -- (G);
    \draw[->] (G) -- (H);
    \draw[->] (I) -- (J);
    \draw[<-] (J) -- (K);
    \draw[->] (K) -- (L);
    \draw[->] (M) -- (N);
    \draw[<-] (N) -- (O);
    \draw[->] (O) -- (P);

    \draw[<-] (A) -- (E);
    \draw[->] (E) -- (I);
    \draw[<-] (I) -- (M);
    \draw[<-] (B) -- (F);
    \draw[->] (F) -- (J);
    \draw[<-] (J) -- (N);
    \draw[<-] (C) -- (G);
    \draw[->] (G) -- (K);
    \draw[<-] (K) -- (O);
    \draw[<-] (D) -- (H);
    \draw[->] (H) -- (L);
    \draw[<-] (L) -- (P);
    
    \draw[dotted] (A) -- (-1,0);
    \draw[dotted] (E) -- (-1,-2);
    \draw[dotted] (I) -- (-1,-4);
    \draw[dotted] (M) -- (-1,-6);
    \draw[dotted] (D) -- (7,0);
    \draw[dotted] (H) -- (7,-2);
    \draw[dotted] (L) -- (7,-4);
    \draw[dotted] (P) -- (7,-6);
    \draw[dotted] (A) -- (0,1);
    \draw[dotted] (B) -- (2,1);
    \draw[dotted] (C) -- (4,1);
    \draw[dotted] (D) -- (6,1);
    \draw[dotted] (M) -- (0,-7);
    \draw[dotted] (N) -- (2,-7);
    \draw[dotted] (O) -- (4,-7);
    \draw[dotted] (P) -- (6,-7);
\end{tikzpicture} \hfill \,
\end{center}
\caption{Pattern of Kronecker Symbols in Square Packings of Types $(5)$, $(3,7)$}
\label{SquarePartialFigure}
\end{figure}
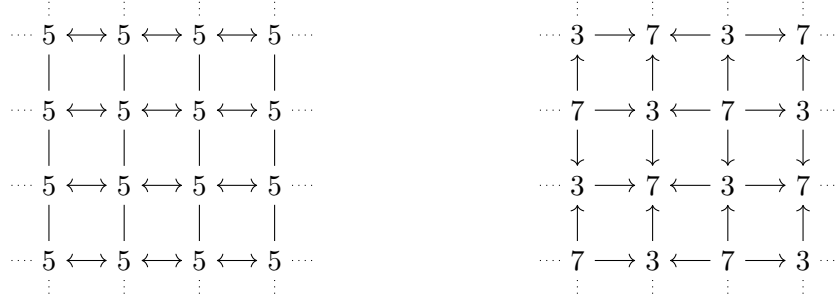

In types $(5)$ and $(3,7)$ we do not find a way of defining $\chi_2$ consistently for all circles in the packing, and we do not find any quadratic obstructions. However, it is still the case that Kronecker symbols between adjacent curvatures have a predictable pattern, as illustrated in Figure \ref{SquarePartialFigure}. In this figure, all pairs of tangent, coprime curvatures $(a,b)$ with an arrow pointing from $a$ to $b$ have the same $\res{a}{b}$ value. All pairs $(a,b)$ without an arrow have the opposite $\res{a}{b}$ value. 

The behavior of Kronecker symbols is governed by Proposition \ref{SquareNode} and quadratic reciprocity. Because each pattern is periodic, with all arrows completely determined by one square of arrows, and the generators of the square Apollonian group all fix four circles arranged in a square, the Apollonian group action is compatible with this pattern. Even though some adjacent pairs of curvatures might not be coprime, Proposition \ref{SquarePath} implies that we can follow the pattern of Kronecker symbols on a path from one to the other.

\section{Triangular Packing}
The triangular packing is built from an infinite collection of circles with a triangular grid tangency graph as shown.
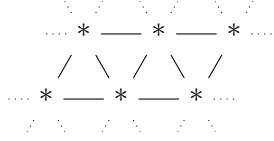
\begin{center}
\begin{tikzpicture}[scale=.5]
    \node (00) at (0,0)    {*};
    \node (10) at (2,0)    {*};  
    \node (20) at (4,0)    {*};     
    \node (01) at (1,0.866*2)    {*};
    \node (11) at (3,0.866*2)    {*};  
    \node (21) at (5,0.866*2)    {*};     
    \draw[dotted] (00) -- (-1,0);
    \draw (00) -- (10);
    \draw (10) -- (20);
    \draw[dotted] (20) -- (5,0);
    \draw[dotted] (01) -- (0,0.866*2);
    \draw (01) -- (11);  
    \draw (11) -- (21);  
    \draw[dotted] (21) -- (6,0.866*2);
    \draw[dotted] (00) -- (-.5,-.866);  
    \draw[dotted] (00) -- (.5,-.866);  
    \draw[dotted] (10) -- (1.5,-.866);  
    \draw[dotted] (10) -- (2.5,-.866);   
    \draw[dotted] (20) -- (3.5,-.866); 
    \draw[dotted] (20) -- (4.5,-.866);   
    \draw (00) -- (01);  
    \draw (10) -- (01);  
    \draw (10) -- (11);   
    \draw (20) -- (11); 
    \draw (20) -- (21); 
    \draw[dotted] (01) -- (.5,3*.866);  
    \draw[dotted] (01) -- (1.5,3*.866);  
    \draw[dotted] (11) -- (2.5,3*.866);  
    \draw[dotted] (11) -- (3.5,3*.866);   
    \draw[dotted] (21) -- (4.5,3*.866); 
    \draw[dotted] (21) -- (5.5,3*.866);   
\end{tikzpicture}
\end{center}

We obtain the triangular packing by iteratively reflecting this confioguration through dual circles. An example of a bounded triangular packing is shown in Figure \ref{TriangularPackingFigure}

\begin{figure}[h] 
\begin{center}
\vspace{-.25in}
\include{justtriangle}
\vspace{-.25in}
\end{center}
\caption{Triangular Packing}
\label{TriangularPackingFigure}
\end{figure}

The curvatures in a triangular grid circle configuration are subject to the following linear and quadratic relations
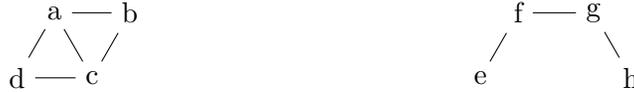
\begin{center}
\, \hfill
\begin{tikzpicture}[scale=.5]
    \node (D) at (0,0)    {d};
    \node (A) at (1,0.866*2)    {a};
    \node (B) at (3,0.866*2)    {b};
    \node (C) at (2,0)    {c};    
    \draw (A) -- (B);
    \draw (A) -- (C);
    \draw (A) -- (D);
    \draw (B) -- (C);  
    \draw (C) -- (D);  
\end{tikzpicture} \hfill
\begin{tikzpicture}[scale=.5]
    \node (E) at (0,0)    {e};
    \node (F) at (1,0.866*2)    {f};
    \node (G) at (3,0.866*2)    {g};
    \node (H) at (4,0)    {h};
    \draw (E) -- (F);
    \draw (F) -- (G);
    \draw (G) -- (H);
\end{tikzpicture}
\hfill \,
\end{center}
\begin{align}
    &(3a-b+3c-d)^2=12ac+4bd \label{DiamondQuadratic} \\
    &e+2g=2f+h\label{Trapezoid}
\end{align}
proven in \cite{RWYY}. The relations apply to all configurations of circles arranged in these ways within the grid.

\begin{figure}[h] \label{TriangleFordFig}
\begin{center}
\begin{tikzpicture}
\draw[black, thin] (-2.5980761875,-0.9375) circle (0.0625)node[scale=0.125]{16};
\draw[black, thin] (2.5980761875,-0.9375) circle (0.0625)node[scale=0.125]{16};
\draw[black, thin] (-1.2371791428571428,-0.9795918367346939) circle (0.02040816326530612)node[scale=0.04081632653061224]{49};
\draw[black, thin] (1.2371791428571428,-0.9795918367346939) circle (0.02040816326530612)node[scale=0.04081632653061224]{49};
\draw[black, thin] (-0.8082903733333333,-0.9866666666666667) circle (0.013333333333333334)node[scale=0.02666666666666667]{75};
\draw[black, thin] (0.8082903733333333,-0.9866666666666667) circle (0.013333333333333334)node[scale=0.02666666666666667]{75};
\draw[black, thin] (-0.69282032,-0.99) circle (0.01)node[scale=0.02]{100};
\draw[black, thin] (0.69282032,-0.99) circle (0.01)node[scale=0.02]{100};
\draw[black, thin] (-2.116951,-0.9629629629629629) circle (0.037037037037037035)node[scale=0.07407407407407407]{27};
\draw[black, thin] (2.116951,-0.9629629629629629) circle (0.037037037037037035)node[scale=0.07407407407407407]{27};
\draw[black, thin] (-1.4433756666666666,-0.9791666666666666) circle (0.020833333333333332)node[scale=0.041666666666666664]{48};
\draw[black, thin] (1.4433756666666666,-0.9791666666666666) circle (0.020833333333333332)node[scale=0.041666666666666664]{48};
\draw[black, thin] (-2.42487112,-0.96) circle (0.04)node[scale=0.08]{25};
\draw[black, thin] (2.42487112,-0.96) circle (0.04)node[scale=0.08]{25};
\draw[black, thin] (-1.299038109375,-0.984375) circle (0.015625)node[scale=0.03125]{64};
\draw[black, thin] (1.299038109375,-0.984375) circle (0.015625)node[scale=0.03125]{64};
\draw[black, thin] (-2.5018511481481482,-0.9629629629629629) circle (0.037037037037037035)node[scale=0.07407407407407407]{27};
\draw[black, thin] (2.5018511481481482,-0.9629629629629629) circle (0.037037037037037035)node[scale=0.07407407407407407]{27};
\draw[black, thin] (-1.2701705866666666,-0.9866666666666667) circle (0.013333333333333334)node[scale=0.02666666666666667]{75};
\draw[black, thin] (1.2701705866666666,-0.9866666666666667) circle (0.013333333333333334)node[scale=0.02666666666666667]{75};
\draw[black, thick] (-3.2,-1.0) -- (3.2,-1.0);
\draw[black, thick] (3.2,1.0) -- (-3.2,1.0);
\draw[black, thick] (-1.732051,0.0) circle (1.0)node[scale=2.0]{1};
\draw[black, thick] (1.732051,0.0) circle (1.0)node[scale=2.0]{1};
\draw[black, thin] (-0.5773503333333333,-0.6666666666666666) circle (0.3333333333333333)node[scale=0.6666666666666666]{3};
\draw[black, thin] (0.5773503333333333,-0.6666666666666666) circle (0.3333333333333333)node[scale=0.6666666666666666]{3};
\draw[black, thin] (0.0,-0.75) circle (0.25)node[scale=0.5]{4};
\draw[black, thin] (-2.0207259375,-0.9791666666666666) circle (0.020833333333333332)node[scale=0.041666666666666664]{48};
\draw[black, thin] (2.0207259375,-0.9791666666666666) circle (0.020833333333333332)node[scale=0.041666666666666664]{48};
\draw[black, thin] (-1.5011106933333334,-0.9866666666666667) circle (0.013333333333333334)node[scale=0.02666666666666667]{75};
\draw[black, thin] (1.5011106933333334,-0.9866666666666667) circle (0.013333333333333334)node[scale=0.02666666666666667]{75};
\draw[black, thin] (-1.1547005,-0.9166666666666666) circle (0.08333333333333333)node[scale=0.16666666666666666]{12};
\draw[black, thin] (1.1547005,-0.9166666666666666) circle (0.08333333333333333)node[scale=0.16666666666666666]{12};
\draw[black, thin] (-0.866025375,-0.9375) circle (0.0625)node[scale=0.125]{16};
\draw[black, thin] (0.866025375,-0.9375) circle (0.0625)node[scale=0.125]{16};
\draw[black, thin] (-0.34641016,-0.96) circle (0.04)node[scale=0.08]{25};
\draw[black, thin] (0.34641016,-0.96) circle (0.04)node[scale=0.08]{25};
\draw[black, thin] (-0.19245007407407405,-0.9629629629629629) circle (0.037037037037037035)node[scale=0.07407407407407407]{27};
\draw[black, thin] (0.19245007407407405,-0.9629629629629629) circle (0.037037037037037035)node[scale=0.07407407407407407]{27};
\draw[black, thin] (-2.8867513333333332,-0.6666666666666666) circle (0.3333333333333333)node[scale=0.6666666666666666]{3};
\draw[black, thin] (2.8867513333333332,-0.6666666666666666) circle (0.3333333333333333)node[scale=0.6666666666666666]{3};
\draw[black, thin] (-3.1176914399999998,-0.96) circle (0.04)node[scale=0.08]{25};
\draw[black, thin] (3.1176914399999998,-0.96) circle (0.04)node[scale=0.08]{25};
\draw[black, thin] (-1.38564065,-0.99) circle (0.01)node[scale=0.02]{100};
\draw[black, thin] (1.38564065,-0.99) circle (0.01)node[scale=0.02]{100};
\draw[black, thin] (-2.2269224693877554,-0.9795918367346939) circle (0.02040816326530612)node[scale=0.04081632653061224]{49};
\draw[black, thin] (2.2269224693877554,-0.9795918367346939) circle (0.02040816326530612)node[scale=0.04081632653061224]{49};
\draw[black, thin] (-2.165063515625,-0.984375) circle (0.015625)node[scale=0.03125]{64};
\draw[black, thin] (2.165063515625,-0.984375) circle (0.015625)node[scale=0.03125]{64};
\draw[black, thin] (-1.9629909200000002,-0.9866666666666667) circle (0.013333333333333334)node[scale=0.02666666666666667]{75};
\draw[black, thin] (1.9629909200000002,-0.9866666666666667) circle (0.013333333333333334)node[scale=0.02666666666666667]{75};
\draw[black, thin] (-2.7217941224489794,-0.9795918367346939) circle (0.02040816326530612)node[scale=0.04081632653061224]{49};
\draw[black, thin] (2.7217941224489794,-0.9795918367346939) circle (0.02040816326530612)node[scale=0.04081632653061224]{49};
\draw[black, thin] (-2.1939310266666667,-0.9866666666666667) circle (0.013333333333333334)node[scale=0.02666666666666667]{75};
\draw[black, thin] (2.1939310266666667,-0.9866666666666667) circle (0.013333333333333334)node[scale=0.02666666666666667]{75};
\draw[black, thin] (-1.0392304799999998,-0.96) circle (0.04)node[scale=0.08]{25};
\draw[black, thin] (1.0392304799999998,-0.96) circle (0.04)node[scale=0.08]{25};
\draw[black, thin] (-0.9622504444444444,-0.9629629629629629) circle (0.037037037037037035)node[scale=0.07407407407407407]{27};
\draw[black, thin] (0.9622504444444444,-0.9629629629629629) circle (0.037037037037037035)node[scale=0.07407407407407407]{27};
\draw[black, thin] (-0.288675125,-0.9791666666666666) circle (0.020833333333333332)node[scale=0.041666666666666664]{48};
\draw[black, thin] (0.288675125,-0.9791666666666666) circle (0.020833333333333332)node[scale=0.041666666666666664]{48};
\draw[black, thin] (-0.24743583673469388,-0.9795918367346939) circle (0.02040816326530612)node[scale=0.04081632653061224]{49};
\draw[black, thin] (0.24743583673469388,-0.9795918367346939) circle (0.02040816326530612)node[scale=0.04081632653061224]{49};
\draw[black, thin] (-2.07846097,-0.99) circle (0.01)node[scale=0.02]{100};
\draw[black, thin] (2.07846097,-0.99) circle (0.01)node[scale=0.02]{100};
\draw[black, thin] (-3.1754264791666666,-0.9791666666666666) circle (0.020833333333333332)node[scale=0.041666666666666664]{48};
\draw[black, thin] (3.1754264791666666,-0.9791666666666666) circle (0.020833333333333332)node[scale=0.041666666666666664]{48};
\draw[black, thin] (-2.6558112400000002,-0.9866666666666667) circle (0.013333333333333334)node[scale=0.02666666666666667]{75};
\draw[black, thin] (2.6558112400000002,-0.9866666666666667) circle (0.013333333333333334)node[scale=0.02666666666666667]{75};
\draw[black, thin] (-3.03108890625,-0.984375) circle (0.015625)node[scale=0.03125]{64};
\draw[black, thin] (3.03108890625,-0.984375) circle (0.015625)node[scale=0.03125]{64};
\draw[black, thin] (-2.3094010833333334,-0.9166666666666666) circle (0.08333333333333333)node[scale=0.16666666666666666]{12};
\draw[black, thin] (2.3094010833333334,-0.9166666666666666) circle (0.08333333333333333)node[scale=0.16666666666666666]{12};
\draw[black, thin] (-1.3471506296296296,-0.9629629629629629) circle (0.037037037037037035)node[scale=0.07407407407407407]{27};
\draw[black, thin] (1.3471506296296296,-0.9629629629629629) circle (0.037037037037037035)node[scale=0.07407407407407407]{27};
\draw[black, thin] (-0.7423074897959183,-0.9795918367346939) circle (0.02040816326530612)node[scale=0.04081632653061224]{49};
\draw[black, thin] (0.7423074897959183,-0.9795918367346939) circle (0.02040816326530612)node[scale=0.04081632653061224]{49};
\draw[black, thin] (-0.433012703125,-0.984375) circle (0.015625)node[scale=0.03125]{64};
\draw[black, thin] (0.433012703125,-0.984375) circle (0.015625)node[scale=0.03125]{64};
\draw[black, thin] (-0.11547005333333334,-0.9866666666666667) circle (0.013333333333333334)node[scale=0.02666666666666667]{75};
\draw[black, thin] (0.11547005333333334,-0.9866666666666667) circle (0.013333333333333334)node[scale=0.02666666666666667]{75};
\draw[black, thin] (-2.77128129,-0.99) circle (0.01)node[scale=0.02]{100};
\draw[black, thin] (2.77128129,-0.99) circle (0.01)node[scale=0.02]{100};

\draw[black, dashed] (-1.732051,1.5) -- (-1.732051,-1.5);
\draw[black, dashed] (1.732051,-1.5) -- (1.732051,1.5);
\draw[black, dashed] (-1.1547004100918505,-1.0) circle (0.5773502050459253);
\draw[black, dashed] (1.1547004100918505,-1.0) circle (0.5773502050459253);
\draw[black, dashed] (-0.28867510252296263,-1.0) circle (0.28867510252296263);
\draw[black, dashed] (0.28867510252296263,-1.0) circle (0.28867510252296263);
\end{tikzpicture}
\end{center}
\caption{Triangular Ford Circles and Duals}
\end{figure}
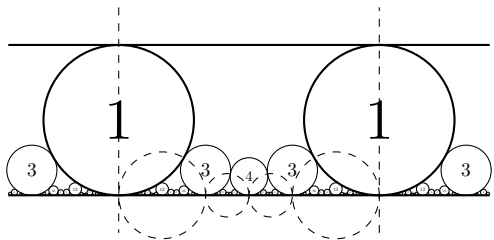

From Figure \ref{TriangleFordFig}, we can observe that the symmetry group for this packing contains a congruence subgroup of level 6. The symmetry group can be generated by the following M\"obius transformations: the translations $T_1=\begin{pmatrix} 1 & 2\sqrt{-3} \\ 0 & 1 \end{pmatrix}$, $T_2=\begin{pmatrix} 1 & 3-\sqrt{-3} \\ 0 & 1 \end{pmatrix}$, and the inversion $S=\begin{pmatrix} 0 & -1 \\ 1 & 0 \end{pmatrix}$. Setting $H=\begin{pmatrix} 1 & 0 \\ 0 & -1\end{pmatrix}, \, R=\begin{pmatrix} 0 & 1 \\ 0 & 0\end{pmatrix}, \, L=\begin{pmatrix} 0 & 0 \\ 1 & 0\end{pmatrix} \in \mathfrak{sl}(2, \Z)$, we have
\begin{equation*}
\begin{split}
    & SHS^{-1}= -H, \quad S R S^{-1} = -L \quad S L S^{-1}=-R \\
    & T_2 L T_2^{-1} = (3-\sqrt{-3}) H + (-6+6\sqrt{-3})R + L, \\
    &T_2^ H T_2^{-1}= H+ (-6+2\sqrt{-3}) R, \quad S T_2 H T_2^{-1} S^{-1} = -H+(6-2\sqrt{-3}) L
\end{split}    
\end{equation*}
Therefore the image of the adjoint action contains $4\mathfrak{sl}(2, \Z[\frac{1+\sqrt{-3}}{2}])$. By \cite[Theorem 8.1]{FuchsStangeZhang}, the modulus for the congruence restrictions for this packing divides $24$. In Proposition \ref{TriangleMod}, we find the possible congruence restrictions mod 12. Lemmas \ref{TriangleModLemma4} and \ref{TriangleModLemma5} imply that the Apollonian group acts transitively on the collection of configurations mod 24 which are lifts of a given configuration mod 12, so 12 is the exact modulus for this packing type.

\begin{proposition} \label{TriangleMod}
A primitive integral triangular grid configuration of curvatures must be congruent to one of the following modulo 12:
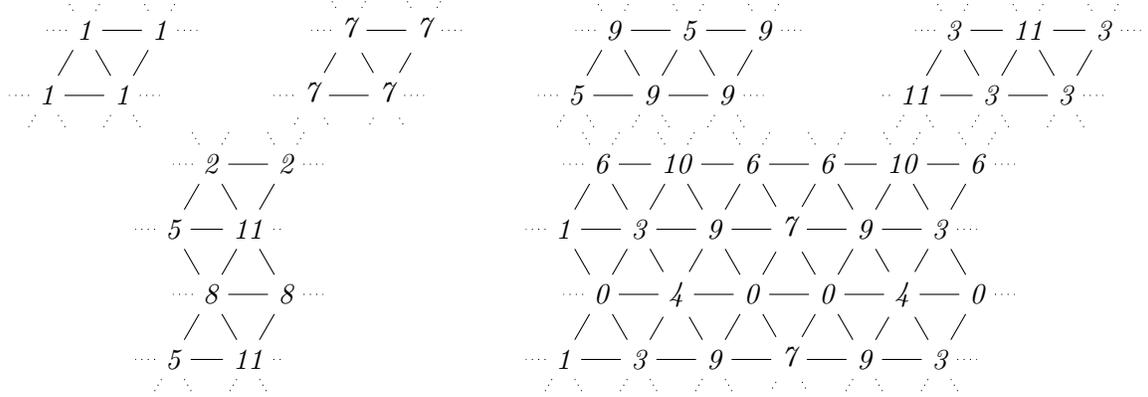
\begin{center}
\, \hfill
\begin{tikzpicture}[scale=.5]
    \node (00) at (0,0)    {1};
    \node (10) at (2,0)    {1};    
    \node (01) at (1,0.866*2)    {1};
    \node (11) at (3,0.866*2)    {1};     
    \draw[dotted] (00) -- (-1,0);
    \draw (00) -- (10);
    \draw[dotted] (10) -- (3,0);
    \draw[dotted] (01) -- (0,0.866*2);
    \draw (01) -- (11);  
    \draw[dotted] (11) -- (4,0.866*2);
    \draw[dotted] (00) -- (-.5,-.866);  
    \draw[dotted] (00) -- (.5,-.866);  
    \draw[dotted] (10) -- (1.5,-.866);  
    \draw[dotted] (10) -- (2.5,-.866);   
    \draw (00) -- (01);  
    \draw (10) -- (01);  
    \draw (10) -- (11);   
    \draw[dotted] (01) -- (.5,3*.866);  
    \draw[dotted] (01) -- (1.5,3*.866);  
    \draw[dotted] (11) -- (2.5,3*.866);  
    \draw[dotted] (11) -- (3.5,3*.866);    
\end{tikzpicture} \hfill
\begin{tikzpicture}[scale=.5]
    \node (00) at (0,0)    {7};
    \node (10) at (2,0)    {7};    
    \node (01) at (1,0.866*2)    {7};
    \node (11) at (3,0.866*2)    {7};     
    \draw[dotted] (00) -- (-1,0);
    \draw (00) -- (10);
    \draw[dotted] (10) -- (3,0);
    \draw[dotted] (01) -- (0,0.866*2);
    \draw (01) -- (11);  
    \draw[dotted] (11) -- (4,0.866*2);
    \draw[dotted] (00) -- (-.5,-.866);  
    \draw[dotted] (00) -- (.5,-.866);  
    \draw[dotted] (10) -- (1.5,-.866);  
    \draw[dotted] (10) -- (2.5,-.866);   
    \draw (00) -- (01);  
    \draw (10) -- (01);  
    \draw (10) -- (11);   
    \draw[dotted] (01) -- (.5,3*.866);  
    \draw[dotted] (01) -- (1.5,3*.866);  
    \draw[dotted] (11) -- (2.5,3*.866);  
    \draw[dotted] (11) -- (3.5,3*.866);    
\end{tikzpicture} \hfill
\begin{tikzpicture}[scale=.5]
    \node (00) at (0,0)    {5};
    \node (10) at (2,0)    {9};  
    \node (20) at (4,0)    {9};     
    \node (01) at (1,0.866*2)    {9};
    \node (11) at (3,0.866*2)    {5};  
    \node (21) at (5,0.866*2)    {9};     
    \draw[dotted] (00) -- (-1,0);
    \draw (00) -- (10);
    \draw (10) -- (20);
    \draw[dotted] (20) -- (5,0);
    \draw[dotted] (01) -- (0,0.866*2);
    \draw (01) -- (11);  
    \draw (11) -- (21);  
    \draw[dotted] (21) -- (6,0.866*2);
    \draw[dotted] (00) -- (-.5,-.866);  
    \draw[dotted] (00) -- (.5,-.866);  
    \draw[dotted] (10) -- (1.5,-.866);  
    \draw[dotted] (10) -- (2.5,-.866);   
    \draw[dotted] (20) -- (3.5,-.866); 
    \draw[dotted] (20) -- (4.5,-.866);   
    \draw (00) -- (01);  
    \draw (10) -- (01);  
    \draw (10) -- (11);   
    \draw (20) -- (11); 
    \draw (20) -- (21); 
    \draw[dotted] (01) -- (.5,3*.866);  
    \draw[dotted] (01) -- (1.5,3*.866);  
    \draw[dotted] (11) -- (2.5,3*.866);  
    \draw[dotted] (11) -- (3.5,3*.866);   
    \draw[dotted] (21) -- (4.5,3*.866); 
    \draw[dotted] (21) -- (5.5,3*.866);   
\end{tikzpicture} \hfill
\begin{tikzpicture}[scale=.5]
    \node (00) at (0,0)    {11};
    \node (10) at (2,0)    {3};  
    \node (20) at (4,0)    {3};     
    \node (01) at (1,0.866*2)    {3};
    \node (11) at (3,0.866*2)    {11};  
    \node (21) at (5,0.866*2)    {3};     
    \draw[dotted] (00) -- (-1,0);
    \draw (00) -- (10);
    \draw (10) -- (20);
    \draw[dotted] (20) -- (5,0);
    \draw[dotted] (01) -- (0,0.866*2);
    \draw (01) -- (11);  
    \draw (11) -- (21);  
    \draw[dotted] (21) -- (6,0.866*2);
    \draw[dotted] (00) -- (-.5,-.866);  
    \draw[dotted] (00) -- (.5,-.866);  
    \draw[dotted] (10) -- (1.5,-.866);  
    \draw[dotted] (10) -- (2.5,-.866);   
    \draw[dotted] (20) -- (3.5,-.866); 
    \draw[dotted] (20) -- (4.5,-.866);   
    \draw (00) -- (01);  
    \draw (10) -- (01);  
    \draw (10) -- (11);   
    \draw (20) -- (11); 
    \draw (20) -- (21); 
    \draw[dotted] (01) -- (.5,3*.866);  
    \draw[dotted] (01) -- (1.5,3*.866);  
    \draw[dotted] (11) -- (2.5,3*.866);  
    \draw[dotted] (11) -- (3.5,3*.866);   
    \draw[dotted] (21) -- (4.5,3*.866); 
    \draw[dotted] (21) -- (5.5,3*.866);   
\end{tikzpicture} \hfill \, \\

\, \hfill \begin{tikzpicture}[scale=.5]
    \node (00) at (0,0)    {5};
    \node (10) at (2,0)    {11};  
    \node (01) at (1,0.866*2)    {8};
    \node (11) at (3,0.866*2)    {8};  
    \node (02) at (0,0.866*4)    {5};
    \node (12) at (2,0.866*4)    {11};   
    \node (03) at (1,0.866*6)    {2};
    \node (13) at (3,0.866*6)    {2};  
    \draw[dotted] (00) -- (-1,0);
    \draw (00) -- (10);
    \draw[dotted] (10) -- (3,0);
    \draw[dotted] (01) -- (0,0.866*2);
    \draw (01) -- (11);  
    \draw[dotted] (11) -- (4,0.866*2);
    \draw[dotted] (02) -- (-1,0.866*4);
    \draw (02) -- (12);
    \draw[dotted] (12) -- (3,0.866*4);
    \draw[dotted] (03) -- (0,0.866*6);
    \draw (03) -- (13);  
    \draw[dotted] (13) -- (4,0.866*6);
    \draw[dotted] (00) -- (-.5,-.866);  
    \draw[dotted] (00) -- (.5,-.866);  
    \draw[dotted] (10) -- (1.5,-.866);  
    \draw[dotted] (10) -- (2.5,-.866);   
    \draw (00) -- (01);  
    \draw (10) -- (01);  
    \draw (10) -- (11);   
    \draw (02) -- (01);  
    \draw (12) -- (01);  
    \draw (12) -- (11);   
    \draw (02) -- (03);  
    \draw (12) -- (03);  
    \draw (12) -- (13);   
    \draw[dotted] (03) -- (.5,7*.866);  
    \draw[dotted] (03) -- (1.5,7*.866);  
    \draw[dotted] (13) -- (2.5,7*.866);  
    \draw[dotted] (13) -- (3.5,7*.866);   
\end{tikzpicture} \hfill
\begin{tikzpicture}[scale=.5]
    \node (00) at (0,0)    {1};
    \node (10) at (2,0)    {3};  
    \node (20) at (4,0)    {9};     
    \node (30) at (6,0)    {7};    
    \node (40) at (8,0)    {9};     
    \node (50) at (10,0)    {3};   
    \node (01) at (1,0.866*2)    {0};
    \node (11) at (3,0.866*2)    {4};  
    \node (21) at (5,0.866*2)    {0};     
    \node (31) at (7,0.866*2)    {0};    
    \node (41) at (9,0.866*2)    {4};     
    \node (51) at (11,0.866*2)    {0}; 
    \node (02) at (0,0.866*4)    {1};
    \node (12) at (2,0.866*4)    {3};  
    \node (22) at (4,0.866*4)    {9};     
    \node (32) at (6,0.866*4)    {7};    
    \node (42) at (8,0.866*4)    {9};     
    \node (52) at (10,0.866*4)    {3};    
    \node (03) at (1,0.866*6)    {6};
    \node (13) at (3,0.866*6)    {10};  
    \node (23) at (5,0.866*6)    {6};     
    \node (33) at (7,0.866*6)    {6};    
    \node (43) at (9,0.866*6)    {10};     
    \node (53) at (11,0.866*6)    {6}; 
    \draw[dotted] (00) -- (-1,0);
    \draw (00) -- (10);
    \draw (10) -- (20);
    \draw (20) -- (30);
    \draw (30) -- (40);  
    \draw (40) -- (50);  
    \draw[dotted] (50) -- (11,0);
    \draw[dotted] (01) -- (0,0.866*2);
    \draw (01) -- (11);  
    \draw (11) -- (21);  
    \draw (21) -- (31);  
    \draw (31) -- (41);  
    \draw (41) -- (51); 
    \draw[dotted] (51) -- (12,0.866*2);
    \draw[dotted] (02) -- (-1,0.866*4);
    \draw (02) -- (12);
    \draw (12) -- (22);
    \draw (22) -- (32);
    \draw (32) -- (42);  
    \draw (42) -- (52);   
    \draw[dotted] (52) -- (11,0.866*4);
    \draw[dotted] (03) -- (0,0.866*6);
    \draw (03) -- (13);  
    \draw (13) -- (23);  
    \draw (23) -- (33);  
    \draw (33) -- (43);  
    \draw (43) -- (53); 
    \draw[dotted] (53) -- (12,0.866*6);
    \draw[dotted] (00) -- (-.5,-.866);  
    \draw[dotted] (00) -- (.5,-.866);  
    \draw[dotted] (10) -- (1.5,-.866);  
    \draw[dotted] (10) -- (2.5,-.866);   
    \draw[dotted] (20) -- (3.5,-.866); 
    \draw[dotted] (20) -- (4.5,-.866); 
    \draw[dotted] (30) -- (5.5,-.866);   
    \draw[dotted] (30) -- (6.5,-.866); 
    \draw[dotted] (40) -- (7.5,-.866);   
    \draw[dotted] (40) -- (8.5,-.866); 
    \draw[dotted] (50) -- (9.5,-.866);   
    \draw[dotted] (50) -- (10.5,-.866);    
    \draw (00) -- (01);  
    \draw (10) -- (01);  
    \draw (10) -- (11);   
    \draw (20) -- (11); 
    \draw (20) -- (21); 
    \draw (30) -- (21);   
    \draw (30) -- (31); 
    \draw (40) -- (31);   
    \draw (40) -- (41); 
    \draw (50) -- (41);   
    \draw (50) -- (51);    
    \draw (02) -- (01);  
    \draw (12) -- (01);  
    \draw (12) -- (11);   
    \draw (22) -- (11); 
    \draw (22) -- (21); 
    \draw (32) -- (21);   
    \draw (32) -- (31); 
    \draw (42) -- (31);   
    \draw (42) -- (41); 
    \draw (52) -- (41);   
    \draw (52) -- (51); 
    \draw (02) -- (03);  
    \draw (12) -- (03);  
    \draw (12) -- (13);   
    \draw (22) -- (13); 
    \draw (22) -- (23); 
    \draw (32) -- (23);   
    \draw (32) -- (33); 
    \draw (42) -- (33);   
    \draw (42) -- (43); 
    \draw (52) -- (43);   
    \draw (52) -- (53); 
    \draw[dotted] (03) -- (.5,7*.866);  
    \draw[dotted] (03) -- (1.5,7*.866);  
    \draw[dotted] (13) -- (2.5,7*.866);  
    \draw[dotted] (13) -- (3.5,7*.866);   
    \draw[dotted] (23) -- (4.5,7*.866); 
    \draw[dotted] (23) -- (5.5,7*.866); 
    \draw[dotted] (33) -- (6.5,7*.866);   
    \draw[dotted] (33) -- (7.5,7*.866); 
    \draw[dotted] (43) -- (8.5,7*.866);   
    \draw[dotted] (43) -- (9.5,7*.866); 
    \draw[dotted] (53) -- (10.5,7*.866);   
    \draw[dotted] (53) -- (11.5,7*.866);   
\end{tikzpicture} \hfill \,
\end{center}
with the pattern repeating horizontally and vertically. 
\end{proposition}

\begin{lemma} \label{TriangleModLemma1}
A primitive integral triangular grid configuration of curvatures must be congruent to one of the following modulo 4:
\end{lemma}
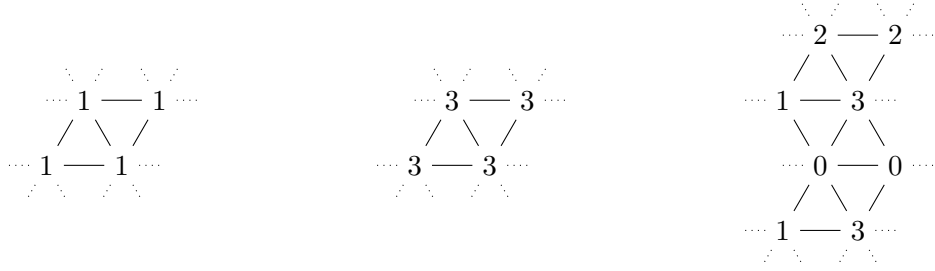
\begin{center}
\, \hfill \begin{tikzpicture}[scale=.5]
    \node (00) at (0,0)    {1};
    \node (10) at (2,0)    {1};    
    \node (01) at (1,0.866*2)    {1};
    \node (11) at (3,0.866*2)    {1};     
    \draw[dotted] (00) -- (-1,0);
    \draw (00) -- (10);
    \draw[dotted] (10) -- (3,0);
    \draw[dotted] (01) -- (0,0.866*2);
    \draw (01) -- (11);  
    \draw[dotted] (11) -- (4,0.866*2);
    \draw[dotted] (00) -- (-.5,-.866);  
    \draw[dotted] (00) -- (.5,-.866);  
    \draw[dotted] (10) -- (1.5,-.866);  
    \draw[dotted] (10) -- (2.5,-.866);   
    \draw (00) -- (01);  
    \draw (10) -- (01);  
    \draw (10) -- (11);   
    \draw[dotted] (01) -- (.5,3*.866);  
    \draw[dotted] (01) -- (1.5,3*.866);  
    \draw[dotted] (11) -- (2.5,3*.866);  
    \draw[dotted] (11) -- (3.5,3*.866); 
    \draw[white] (00) -- (0, -3*.866);
\end{tikzpicture} \hfill
\begin{tikzpicture}[scale=.5]
    \node (00) at (0,0)    {3};
    \node (10) at (2,0)    {3};    
    \node (01) at (1,0.866*2)    {3};
    \node (11) at (3,0.866*2)    {3};     
    \draw[dotted] (00) -- (-1,0);
    \draw (00) -- (10);
    \draw[dotted] (10) -- (3,0);
    \draw[dotted] (01) -- (0,0.866*2);
    \draw (01) -- (11);  
    \draw[dotted] (11) -- (4,0.866*2);
    \draw[dotted] (00) -- (-.5,-.866);  
    \draw[dotted] (00) -- (.5,-.866);  
    \draw[dotted] (10) -- (1.5,-.866);  
    \draw[dotted] (10) -- (2.5,-.866);   
    \draw (00) -- (01);  
    \draw (10) -- (01);  
    \draw (10) -- (11);   
    \draw[dotted] (01) -- (.5,3*.866);  
    \draw[dotted] (01) -- (1.5,3*.866);  
    \draw[dotted] (11) -- (2.5,3*.866);  
    \draw[dotted] (11) -- (3.5,3*.866);  
    \draw[white] (00) -- (0, -3*.866);  
\end{tikzpicture} \hfill
\begin{tikzpicture}[scale=.5]
    \node (00) at (0,0)    {1};
    \node (10) at (2,0)    {3};  
    \node (01) at (1,0.866*2)    {0};
    \node (11) at (3,0.866*2)    {0};  
    \node (02) at (0,0.866*4)    {1};
    \node (12) at (2,0.866*4)    {3};   
    \node (03) at (1,0.866*6)    {2};
    \node (13) at (3,0.866*6)    {2};  
    \draw[dotted] (00) -- (-1,0);
    \draw (00) -- (10);
    \draw[dotted] (10) -- (3,0);
    \draw[dotted] (01) -- (0,0.866*2);
    \draw (01) -- (11);  
    \draw[dotted] (11) -- (4,0.866*2);
    \draw[dotted] (02) -- (-1,0.866*4);
    \draw (02) -- (12);
    \draw[dotted] (12) -- (3,0.866*4);
    \draw[dotted] (03) -- (0,0.866*6);
    \draw (03) -- (13);  
    \draw[dotted] (13) -- (4,0.866*6);
    \draw[dotted] (00) -- (-.5,-.866);  
    \draw[dotted] (00) -- (.5,-.866);  
    \draw[dotted] (10) -- (1.5,-.866);  
    \draw[dotted] (10) -- (2.5,-.866);   
    \draw (00) -- (01);  
    \draw (10) -- (01);  
    \draw (10) -- (11);   
    \draw (02) -- (01);  
    \draw (12) -- (01);  
    \draw (12) -- (11);   
    \draw (02) -- (03);  
    \draw (12) -- (03);  
    \draw (12) -- (13);   
    \draw[dotted] (03) -- (.5,7*.866);  
    \draw[dotted] (03) -- (1.5,7*.866);  
    \draw[dotted] (13) -- (2.5,7*.866);  
    \draw[dotted] (13) -- (3.5,7*.866);   
\end{tikzpicture} \hfill \,
\end{center}
\begin{proof}
Note that four curvatures $a,b,c,d$ satisfying Equation \eqref{DiamondQuadratic} uniquely determine the configuration via the linear relation \eqref{Trapezoid}. Reducing mod 2, we find that $a+b+c+d\equiv 0 \bmod 2$, so either $a,b,c,d$ are all odd or two are odd and two are even. Therefore, either all curvatures in the grid are odd or every diamond in the grid contains two odd and two even curvatures.

Equation \eqref{DiamondQuadratic} gives:
\begin{equation*}
    d=\pm 2 \sqrt{3}\sqrt{ab+bc+ac}+3a+b+3c
\end{equation*}
Since $d$ is integral, we see $\exists m \in \mathbb{Z}$ where $ab+bc+ac=3m^2$. Thus, $(a+b)(a+c)=a^2+3m^2$.

In the case when all curvatures are odd, we have $a^2+3m^2 \equiv 4 \bmod 8$, and thus $a+b, \, a+c \equiv 2 \bmod 4$. It follows that $a\equiv b \equiv c \bmod 4$, and this extends across the entire grid, so we obtain one of the first two configurations. 

In the case where both odd and even curvatures appear, we can find a diamond with $a,b$ odd and $c,d$ even. Then $a^2+3m^2 \equiv 4 \bmod 8$, and thus $a+b \equiv 0 \bmod 4$. Applying the same argument with $c, d$ exchanged for $a,b$ gives $c+d \equiv 0 \bmod 4$. Since the third configuration shown above satisfies the linear relation \eqref{Trapezoid} and contains every possible diamond with $a \equiv -b \bmod 4$ odd and $c\equiv -d \bmod 4$ even, it is the unique configuration with both odd and even curvatures.
\end{proof}

\begin{lemma} \label{TriangleModLemma2}
A primitive integral triangular grid configuration of curvatures must be congruent to one of the following modulo 3:
\end{lemma}
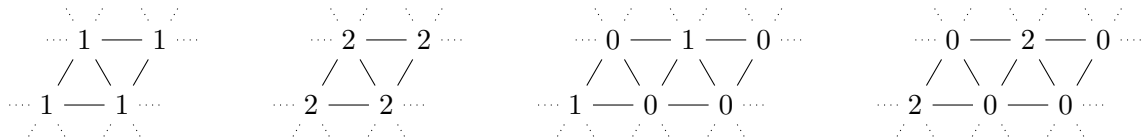
\begin{center}
\, \hfill
\begin{tikzpicture}[scale=.5]
    \node (00) at (0,0)    {1};
    \node (10) at (2,0)    {1};    
    \node (01) at (1,0.866*2)    {1};
    \node (11) at (3,0.866*2)    {1};     
    \draw[dotted] (00) -- (-1,0);
    \draw (00) -- (10);
    \draw[dotted] (10) -- (3,0);
    \draw[dotted] (01) -- (0,0.866*2);
    \draw (01) -- (11);  
    \draw[dotted] (11) -- (4,0.866*2);
    \draw[dotted] (00) -- (-.5,-.866);  
    \draw[dotted] (00) -- (.5,-.866);  
    \draw[dotted] (10) -- (1.5,-.866);  
    \draw[dotted] (10) -- (2.5,-.866);   
    \draw (00) -- (01);  
    \draw (10) -- (01);  
    \draw (10) -- (11);   
    \draw[dotted] (01) -- (.5,3*.866);  
    \draw[dotted] (01) -- (1.5,3*.866);  
    \draw[dotted] (11) -- (2.5,3*.866);  
    \draw[dotted] (11) -- (3.5,3*.866);    
\end{tikzpicture} \hfill
\begin{tikzpicture}[scale=.5]
    \node (00) at (0,0)    {2};
    \node (10) at (2,0)    {2};    
    \node (01) at (1,0.866*2)    {2};
    \node (11) at (3,0.866*2)    {2};     
    \draw[dotted] (00) -- (-1,0);
    \draw (00) -- (10);
    \draw[dotted] (10) -- (3,0);
    \draw[dotted] (01) -- (0,0.866*2);
    \draw (01) -- (11);  
    \draw[dotted] (11) -- (4,0.866*2);
    \draw[dotted] (00) -- (-.5,-.866);  
    \draw[dotted] (00) -- (.5,-.866);  
    \draw[dotted] (10) -- (1.5,-.866);  
    \draw[dotted] (10) -- (2.5,-.866);   
    \draw (00) -- (01);  
    \draw (10) -- (01);  
    \draw (10) -- (11);   
    \draw[dotted] (01) -- (.5,3*.866);  
    \draw[dotted] (01) -- (1.5,3*.866);  
    \draw[dotted] (11) -- (2.5,3*.866);  
    \draw[dotted] (11) -- (3.5,3*.866);    
\end{tikzpicture} \hfill
\begin{tikzpicture}[scale=.5]
    \node (00) at (0,0)    {1};
    \node (10) at (2,0)    {0};  
    \node (20) at (4,0)    {0};     
    \node (01) at (1,0.866*2)    {0};
    \node (11) at (3,0.866*2)    {1};  
    \node (21) at (5,0.866*2)    {0};     
    \draw[dotted] (00) -- (-1,0);
    \draw (00) -- (10);
    \draw (10) -- (20);
    \draw[dotted] (20) -- (5,0);
    \draw[dotted] (01) -- (0,0.866*2);
    \draw (01) -- (11);  
    \draw (11) -- (21);  
    \draw[dotted] (21) -- (6,0.866*2);
    \draw[dotted] (00) -- (-.5,-.866);  
    \draw[dotted] (00) -- (.5,-.866);  
    \draw[dotted] (10) -- (1.5,-.866);  
    \draw[dotted] (10) -- (2.5,-.866);   
    \draw[dotted] (20) -- (3.5,-.866); 
    \draw[dotted] (20) -- (4.5,-.866);   
    \draw (00) -- (01);  
    \draw (10) -- (01);  
    \draw (10) -- (11);   
    \draw (20) -- (11); 
    \draw (20) -- (21); 
    \draw[dotted] (01) -- (.5,3*.866);  
    \draw[dotted] (01) -- (1.5,3*.866);  
    \draw[dotted] (11) -- (2.5,3*.866);  
    \draw[dotted] (11) -- (3.5,3*.866);   
    \draw[dotted] (21) -- (4.5,3*.866); 
    \draw[dotted] (21) -- (5.5,3*.866);   
\end{tikzpicture} \hfill
\begin{tikzpicture}[scale=.5]
    \node (00) at (0,0)    {2};
    \node (10) at (2,0)    {0};  
    \node (20) at (4,0)    {0};     
    \node (01) at (1,0.866*2)    {0};
    \node (11) at (3,0.866*2)    {2};  
    \node (21) at (5,0.866*2)    {0};     
    \draw[dotted] (00) -- (-1,0);
    \draw (00) -- (10);
    \draw (10) -- (20);
    \draw[dotted] (20) -- (5,0);
    \draw[dotted] (01) -- (0,0.866*2);
    \draw (01) -- (11);  
    \draw (11) -- (21);  
    \draw[dotted] (21) -- (6,0.866*2);
    \draw[dotted] (00) -- (-.5,-.866);  
    \draw[dotted] (00) -- (.5,-.866);  
    \draw[dotted] (10) -- (1.5,-.866);  
    \draw[dotted] (10) -- (2.5,-.866);   
    \draw[dotted] (20) -- (3.5,-.866); 
    \draw[dotted] (20) -- (4.5,-.866);   
    \draw (00) -- (01);  
    \draw (10) -- (01);  
    \draw (10) -- (11);   
    \draw (20) -- (11); 
    \draw (20) -- (21); 
    \draw[dotted] (01) -- (.5,3*.866);  
    \draw[dotted] (01) -- (1.5,3*.866);  
    \draw[dotted] (11) -- (2.5,3*.866);  
    \draw[dotted] (11) -- (3.5,3*.866);   
    \draw[dotted] (21) -- (4.5,3*.866); 
    \draw[dotted] (21) -- (5.5,3*.866);   
\end{tikzpicture} \hfill \,
\end{center}
\begin{proof}
Equation \ref{DiamondQuadratic} reduces to $(b-d)^2 \equiv 0 \bmod 3$, so $b \equiv d \bmod 3$. Thus the entire configuration mod 3 is determined by three curvatures $a,b,c$ arranged in a triangle. 

As in the previous proof, we have $ab+bc+ac=3m^2$, so $ab+bc+ac \equiv 0 \bmod 3$. Thus, $a,b,c$ must either all be congruent modulo 3, or at least two of them are 0 modulo 3. Notice that they cannot all be 0 modulo 3, since the packing is primitive. Choosing all the curvatures congruent mod 3 gives the first two configurations. Choosing two curvatures zero and one nonzero gives the other two configurations. 
\end{proof}

\begin{lemma} \label{TriangleModLemma3}
Suppose that $a$, $b$ are the curvatures of two tangent circles in a primitive integral triangular packing. Then, $a+b \not \equiv 5, 10, 11 \bmod 12$.  
\end{lemma}
\begin{proof}
As in Lemma \ref{TriangleModLemma1}, we have $(a+b)(a+c) = a^2 + 3m^2$. Let $p>3$ be a prime number that divides $a^2+3m^2$ and $p \nmid a, m$. Then $-3$ is a square mod $p$, so $p \equiv 1, \, 7 \bmod 12$. And, if a prime $q \equiv 5, \, 11 \bmod 12$ divides $a^2+3m^2$, then $q\mid a, \, m$ and $v_{q}(a^2+3m^2)$ is even. 

Assume that $a+b \equiv 5, 10, 11 \bmod 12$. Then there exists a prime number $q \equiv 5, \, 11 \bmod 12$ such that $q \mid  a+b$ and $v_q(a+b)$ is odd. Then, $q\mid a, \, m$ and $v_q(a^2+3m^2)$ is even. Therefore $q\mid  a+c$, and hence $q\mid  a,b,c$. Moreover, $d = \pm 6m + 3a + b + 3c$ so $q\mid d$. Hence, the packing is not primitive with a common factor of $q$. 
\end{proof}

We can now verify Proposition \ref{TriangleMod}. Using the Chinese remainder theorem, the three possible configurations mod 4 from Lemma \ref{TriangleModLemma1} and the four possible configurations mod 3 from Lemma \ref{TriangleModLemma2} combine into twelve possible configurations mod 12 (the more complicated periodic configurations seem to combine in multiple ways, but they are all the same up to rigid motions). Six of these twelve possibilities are ruled out by Lemma \ref{TriangleModLemma3}, leaving the six shown in Proposition \ref{TriangleMod}.

\begin{proposition} \label{TriangleGroupMod}
In a primitive integral triangular packing, every triangular grid configuration of circles belongs to the same type modulo 12. 
\end{proposition}
This proposition is proven similarly to Proposition \ref{OctGroupMod}. It suffices to check, in Proposition \ref{TriangleMod}, that each configuration is uniquely determined by one triangle of curvatures.

We will refer to the six types of triangular packings mod 12 as type $(1)$, type $(7)$, type $(3, 11)$, type $(5,9)$, type $(2,5,8,11)$, and type $(0,1,3,4,6,7,9,10)$. Section \ref{Data} below contains examples of all six packing types. 

We also need to characterize the lifts of a given configuration mod 12 to mod 24, and to show that the Apollonian group acts transitively on these lifts. This establishes that 12 is the correct modulus for congruence restrictions in these packings. Moreover, lifts to mod 24 will be used to analyze the pattern of Kronecker symbols in primitive integral triangular packings containing both even and odd curvatures. We prove two lemmas describing the lifts from mod 4 to mod 8, first in the case of all odd curvatures, then in the case with odd and even curvatures.

\begin{lemma} \label{TriangleModLemma4} The triangular configuration where all curvatures are 1 mod 4 lifts to mod 8 in five possible ways:
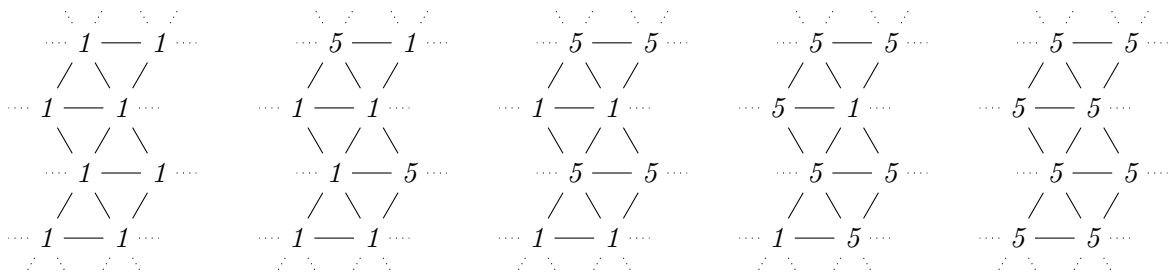
\begin{center} \, \hfill
\begin{tikzpicture}[scale=.5]
    \node (00) at (0,0)    {1};
    \node (10) at (2,0)    {1};  
    \node (01) at (1,0.866*2)    {1};
    \node (11) at (3,0.866*2)    {1};  
    \node (02) at (0,0.866*4)    {1};
    \node (12) at (2,0.866*4)    {1};   
    \node (03) at (1,0.866*6)    {1};
    \node (13) at (3,0.866*6)    {1};  
    \draw[dotted] (00) -- (-1,0);
    \draw (00) -- (10);
    \draw[dotted] (10) -- (3,0);
    \draw[dotted] (01) -- (0,0.866*2);
    \draw (01) -- (11);  
    \draw[dotted] (11) -- (4,0.866*2);
    \draw[dotted] (02) -- (-1,0.866*4);
    \draw (02) -- (12);
    \draw[dotted] (12) -- (3,0.866*4);
    \draw[dotted] (03) -- (0,0.866*6);
    \draw (03) -- (13);  
    \draw[dotted] (13) -- (4,0.866*6);
    \draw[dotted] (00) -- (-.5,-.866);  
    \draw[dotted] (00) -- (.5,-.866);  
    \draw[dotted] (10) -- (1.5,-.866);  
    \draw[dotted] (10) -- (2.5,-.866);   
    \draw (00) -- (01);  
    \draw (10) -- (01);  
    \draw (10) -- (11);   
    \draw (02) -- (01);  
    \draw (12) -- (01);  
    \draw (12) -- (11);   
    \draw (02) -- (03);  
    \draw (12) -- (03);  
    \draw (12) -- (13);   
    \draw[dotted] (03) -- (.5,7*.866);  
    \draw[dotted] (03) -- (1.5,7*.866);  
    \draw[dotted] (13) -- (2.5,7*.866);  
    \draw[dotted] (13) -- (3.5,7*.866);   
\end{tikzpicture} \hfill
\begin{tikzpicture}[scale=.5]
    \node (00) at (0,0)    {1};
    \node (10) at (2,0)    {1};  
    \node (01) at (1,0.866*2)    {1};
    \node (11) at (3,0.866*2)    {5};  
    \node (02) at (0,0.866*4)    {1};
    \node (12) at (2,0.866*4)    {1};   
    \node (03) at (1,0.866*6)    {5};
    \node (13) at (3,0.866*6)    {1};  
    \draw[dotted] (00) -- (-1,0);
    \draw (00) -- (10);
    \draw[dotted] (10) -- (3,0);
    \draw[dotted] (01) -- (0,0.866*2);
    \draw (01) -- (11);  
    \draw[dotted] (11) -- (4,0.866*2);
    \draw[dotted] (02) -- (-1,0.866*4);
    \draw (02) -- (12);
    \draw[dotted] (12) -- (3,0.866*4);
    \draw[dotted] (03) -- (0,0.866*6);
    \draw (03) -- (13);  
    \draw[dotted] (13) -- (4,0.866*6);
    \draw[dotted] (00) -- (-.5,-.866);  
    \draw[dotted] (00) -- (.5,-.866);  
    \draw[dotted] (10) -- (1.5,-.866);  
    \draw[dotted] (10) -- (2.5,-.866);   
    \draw (00) -- (01);  
    \draw (10) -- (01);  
    \draw (10) -- (11);   
    \draw (02) -- (01);  
    \draw (12) -- (01);  
    \draw (12) -- (11);   
    \draw (02) -- (03);  
    \draw (12) -- (03);  
    \draw (12) -- (13);   
    \draw[dotted] (03) -- (.5,7*.866);  
    \draw[dotted] (03) -- (1.5,7*.866);  
    \draw[dotted] (13) -- (2.5,7*.866);  
    \draw[dotted] (13) -- (3.5,7*.866);   
\end{tikzpicture}\hfill
\begin{tikzpicture}[scale=.5]
    \node (00) at (0,0)    {1};
    \node (10) at (2,0)    {1};  
    \node (01) at (1,0.866*2)    {5};
    \node (11) at (3,0.866*2)    {5};  
    \node (02) at (0,0.866*4)    {1};
    \node (12) at (2,0.866*4)    {1};   
    \node (03) at (1,0.866*6)    {5};
    \node (13) at (3,0.866*6)    {5};  
    \draw[dotted] (00) -- (-1,0);
    \draw (00) -- (10);
    \draw[dotted] (10) -- (3,0);
    \draw[dotted] (01) -- (0,0.866*2);
    \draw (01) -- (11);  
    \draw[dotted] (11) -- (4,0.866*2);
    \draw[dotted] (02) -- (-1,0.866*4);
    \draw (02) -- (12);
    \draw[dotted] (12) -- (3,0.866*4);
    \draw[dotted] (03) -- (0,0.866*6);
    \draw (03) -- (13);  
    \draw[dotted] (13) -- (4,0.866*6);
    \draw[dotted] (00) -- (-.5,-.866);  
    \draw[dotted] (00) -- (.5,-.866);  
    \draw[dotted] (10) -- (1.5,-.866);  
    \draw[dotted] (10) -- (2.5,-.866);   
    \draw (00) -- (01);  
    \draw (10) -- (01);  
    \draw (10) -- (11);   
    \draw (02) -- (01);  
    \draw (12) -- (01);  
    \draw (12) -- (11);   
    \draw (02) -- (03);  
    \draw (12) -- (03);  
    \draw (12) -- (13);   
    \draw[dotted] (03) -- (.5,7*.866);  
    \draw[dotted] (03) -- (1.5,7*.866);  
    \draw[dotted] (13) -- (2.5,7*.866);  
    \draw[dotted] (13) -- (3.5,7*.866);   
\end{tikzpicture}\hfill
\begin{tikzpicture}[scale=.5]
    \node (00) at (0,0)    {1};
    \node (10) at (2,0)    {5};  
    \node (01) at (1,0.866*2)    {5};
    \node (11) at (3,0.866*2)    {5};  
    \node (02) at (0,0.866*4)    {5};
    \node (12) at (2,0.866*4)    {1};   
    \node (03) at (1,0.866*6)    {5};
    \node (13) at (3,0.866*6)    {5};  
    \draw[dotted] (00) -- (-1,0);
    \draw (00) -- (10);
    \draw[dotted] (10) -- (3,0);
    \draw[dotted] (01) -- (0,0.866*2);
    \draw (01) -- (11);  
    \draw[dotted] (11) -- (4,0.866*2);
    \draw[dotted] (02) -- (-1,0.866*4);
    \draw (02) -- (12);
    \draw[dotted] (12) -- (3,0.866*4);
    \draw[dotted] (03) -- (0,0.866*6);
    \draw (03) -- (13);  
    \draw[dotted] (13) -- (4,0.866*6);
    \draw[dotted] (00) -- (-.5,-.866);  
    \draw[dotted] (00) -- (.5,-.866);  
    \draw[dotted] (10) -- (1.5,-.866);  
    \draw[dotted] (10) -- (2.5,-.866);   
    \draw (00) -- (01);  
    \draw (10) -- (01);  
    \draw (10) -- (11);   
    \draw (02) -- (01);  
    \draw (12) -- (01);  
    \draw (12) -- (11);   
    \draw (02) -- (03);  
    \draw (12) -- (03);  
    \draw (12) -- (13);   
    \draw[dotted] (03) -- (.5,7*.866);  
    \draw[dotted] (03) -- (1.5,7*.866);  
    \draw[dotted] (13) -- (2.5,7*.866);  
    \draw[dotted] (13) -- (3.5,7*.866);   
\end{tikzpicture}\hfill
\begin{tikzpicture}[scale=.5]
    \node (00) at (0,0)    {5};
    \node (10) at (2,0)    {5};  
    \node (01) at (1,0.866*2)    {5};
    \node (11) at (3,0.866*2)    {5};  
    \node (02) at (0,0.866*4)    {5};
    \node (12) at (2,0.866*4)    {5};   
    \node (03) at (1,0.866*6)    {5};
    \node (13) at (3,0.866*6)    {5};  
    \draw[dotted] (00) -- (-1,0);
    \draw (00) -- (10);
    \draw[dotted] (10) -- (3,0);
    \draw[dotted] (01) -- (0,0.866*2);
    \draw (01) -- (11);  
    \draw[dotted] (11) -- (4,0.866*2);
    \draw[dotted] (02) -- (-1,0.866*4);
    \draw (02) -- (12);
    \draw[dotted] (12) -- (3,0.866*4);
    \draw[dotted] (03) -- (0,0.866*6);
    \draw (03) -- (13);  
    \draw[dotted] (13) -- (4,0.866*6);
    \draw[dotted] (00) -- (-.5,-.866);  
    \draw[dotted] (00) -- (.5,-.866);  
    \draw[dotted] (10) -- (1.5,-.866);  
    \draw[dotted] (10) -- (2.5,-.866);   
    \draw (00) -- (01);  
    \draw (10) -- (01);  
    \draw (10) -- (11);   
    \draw (02) -- (01);  
    \draw (12) -- (01);  
    \draw (12) -- (11);   
    \draw (02) -- (03);  
    \draw (12) -- (03);  
    \draw (12) -- (13);   
    \draw[dotted] (03) -- (.5,7*.866);  
    \draw[dotted] (03) -- (1.5,7*.866);  
    \draw[dotted] (13) -- (2.5,7*.866);  
    \draw[dotted] (13) -- (3.5,7*.866);   
\end{tikzpicture} \hfill \,
\end{center}
Similarly, the triangular configuration where all curvatures are 3 mod 4 lifts to mod 8 in five possible ways containing curvatures of 3 and 7. Moreover, the Apollonian group acts transitively on these lifts. 
\end{lemma}

\begin{proof}
It suffices to describe the lift of four curvatures $(a,b,c,d)$ satisfying Equation \eqref{DiamondQuadratic}; then all curvatures mod 8 are determined by Equation \eqref{Trapezoid}. All possible lifts of $(1,1,1,1)\bmod 4$ satisfy Equation \eqref{DiamondQuadratic} mod 8. The five configurations above correspond to the cases when zero, one, two, three, or four of $a,b,c,d$ lift to $5 \bmod 8$. 

Recall that each Apollonian group generator fixes three circles arranged in a triangle. If circles of curvatures $a$, $b$, $c$, $d$ are arranged in a diamond, as in Equation \eqref{DiamondQuadratic}, then the generator which fixes $a$, $b$, $c$ will transform $d$ to $d'=6a+2b+6c-d$. Assuming that $a\equiv b \equiv c \equiv 1 \bmod 4$, we have $d' \equiv 6-d \bmod 8$ so each generator will transform between $d\equiv 1 \bmod 8$ and $d\equiv 5 \bmod 8$. The group then acts transitively on the five configurations shown above, with each generator mapping a configuration to one of its adjacent configurations.

The argument when $(a,b,c,d)\equiv(3,3,3,3) \bmod 4$ is similar.
\end{proof}

In particular, any primitive integral triangular packing with all odd curvatures contains five different configurations mod 8. 

\begin{lemma} \label{TriangleModLemma5} The unique triangular configuration mod 4 containing even curvatures lifts to mod 8 in two possible ways:

\begin{center} \, \hfill
\begin{tikzpicture}[scale=.45]
    \node (00) at (0,0)    {1};
    \node (10) at (2,0)    {3};  
    \node (20) at (4,0)    {1};     
    \node (30) at (6,0)    {3};
    \node (01) at (1,0.866*2)    {2};
    \node (11) at (3,0.866*2)    {2};  
    \node (21) at (5,0.866*2)    {6};     
    \node (31) at (7,0.866*2)    {6};    
    \node (02) at (0,0.866*4)    {5};
    \node (12) at (2,0.866*4)    {7};  
    \node (22) at (4,0.866*4)    {5};     
    \node (32) at (6,0.866*4)    {7};     
    \node (03) at (1,0.866*6)    {0};
    \node (13) at (3,0.866*6)    {0};  
    \node (23) at (5,0.866*6)    {4};     
    \node (33) at (7,0.866*6)    {4};    
    \node (04) at (0,0.866*8)    {5};
    \node (14) at (2,0.866*8)    {7};  
    \node (24) at (4,0.866*8)    {5};     
    \node (34) at (6,0.866*8)    {7};    
    \node (05) at (1,0.866*10)    {2};
    \node (15) at (3,0.866*10)    {2};  
    \node (25) at (5,0.866*10)    {6};     
    \node (35) at (7,0.866*10)    {6};    
    \node (06) at (0,0.866*12)    {1};
    \node (16) at (2,0.866*12)    {3};  
    \node (26) at (4,0.866*12)    {1};     
    \node (36) at (6,0.866*12)    {3};    
    \node (07) at (1,0.866*14)    {0};
    \node (17) at (3,0.866*14)    {0};  
    \node (27) at (5,0.866*14)    {4};     
    \node (37) at (7,0.866*14)    {4};    
    \draw[dotted] (00) -- (-1,0);
    \draw (00) -- (10);
    \draw (10) -- (20);
    \draw (20) -- (30);
    \draw[dotted] (30) -- (7,0);
    \draw[dotted] (01) -- (0,0.866*2);
    \draw (01) -- (11);  
    \draw (11) -- (21);  
    \draw (21) -- (31);  
    \draw[dotted] (31) -- (8,0.866*2);
    \draw[dotted] (02) -- (-1,0.866*4);
    \draw (02) -- (12);
    \draw (12) -- (22);
    \draw (22) -- (32);  
    \draw[dotted] (32) -- (7,0.866*4);
    \draw[dotted] (03) -- (0,0.866*6);
    \draw (03) -- (13);  
    \draw (13) -- (23);  
    \draw (23) -- (33);    
    \draw[dotted] (33) -- (8,0.866*6);
    \draw[dotted] (04) -- (-1,0.866*8);
    \draw (04) -- (14);  
    \draw (14) -- (24);  
    \draw (24) -- (34);
    \draw[dotted] (34) -- (7,0.866*8);
    \draw[dotted] (05) -- (0,0.866*10);
    \draw (05) -- (15);  
    \draw (15) -- (25);  
    \draw (25) -- (35);  
    \draw[dotted] (35) -- (8,0.866*10);
    \draw[dotted] (06) -- (-1,0.866*12);
    \draw (06) -- (16);  
    \draw (16) -- (26);  
    \draw (26) -- (36);  
    \draw[dotted] (36) -- (7,0.866*12);
    \draw[dotted] (07) -- (0,0.866*14);
    \draw (07) -- (17);  
    \draw (17) -- (27);  
    \draw (27) -- (37);  
    \draw[dotted] (37) -- (8,0.866*14);
    \draw[dotted] (00) -- (-.5,-.866);  
    \draw[dotted] (00) -- (.5,-.866);  
    \draw[dotted] (10) -- (1.5,-.866);  
    \draw[dotted] (10) -- (2.5,-.866);   
    \draw[dotted] (20) -- (3.5,-.866); 
    \draw[dotted] (20) -- (4.5,-.866); 
    \draw[dotted] (30) -- (5.5,-.866);   
    \draw[dotted] (30) -- (6.5,-.866);  
    \draw (00) -- (01);  
    \draw (10) -- (01);  
    \draw (10) -- (11);   
    \draw (20) -- (11); 
    \draw (20) -- (21); 
    \draw (30) -- (21);   
    \draw (30) -- (31); 
    \draw (02) -- (01);  
    \draw (12) -- (01);  
    \draw (12) -- (11);   
    \draw (22) -- (11); 
    \draw (22) -- (21); 
    \draw (32) -- (21);   
    \draw (32) -- (31); 
    \draw (02) -- (03);  
    \draw (12) -- (03);  
    \draw (12) -- (13);   
    \draw (22) -- (13); 
    \draw (22) -- (23); 
    \draw (32) -- (23);   
    \draw (32) -- (33);  
    \draw (04) -- (03);  
    \draw (14) -- (03);  
    \draw (14) -- (13);   
    \draw (24) -- (13); 
    \draw (24) -- (23); 
    \draw (34) -- (23);   
    \draw (34) -- (33); 
    \draw (04) -- (05);  
    \draw (14) -- (05);  
    \draw (14) -- (15);   
    \draw (24) -- (15); 
    \draw (24) -- (25); 
    \draw (34) -- (25);   
    \draw (34) -- (35);  
    \draw (06) -- (05);  
    \draw (16) -- (05);  
    \draw (16) -- (15);   
    \draw (26) -- (15); 
    \draw (26) -- (25); 
    \draw (36) -- (25);   
    \draw (36) -- (35); 
    \draw (06) -- (07);  
    \draw (16) -- (07);  
    \draw (16) -- (17);   
    \draw (26) -- (17); 
    \draw (26) -- (27); 
    \draw (36) -- (27);   
    \draw (36) -- (37); 
    \draw[dotted] (07) -- (.5,15*.866);  
    \draw[dotted] (07) -- (1.5,15*.866);  
    \draw[dotted] (17) -- (2.5,15*.866);  
    \draw[dotted] (17) -- (3.5,15*.866);   
    \draw[dotted] (27) -- (4.5,15*.866); 
    \draw[dotted] (27) -- (5.5,15*.866); 
    \draw[dotted] (37) -- (6.5,15*.866);   
    \draw[dotted] (37) -- (7.5,15*.866);  
\end{tikzpicture}  \hfill
\begin{tikzpicture}[scale=.45]
    \node (00) at (0,0)    {1};
    \node (10) at (2,0)    {3};  
    \node (20) at (4,0)    {1};     
    \node (30) at (6,0)    {3};
    \node (01) at (1,0.866*2)    {2};
    \node (11) at (3,0.866*2)    {6};  
    \node (21) at (5,0.866*2)    {6};     
    \node (31) at (7,0.866*2)    {2};    
    \node (02) at (0,0.866*4)    {5};
    \node (12) at (2,0.866*4)    {7};  
    \node (22) at (4,0.866*4)    {5};     
    \node (32) at (6,0.866*4)    {7};     
    \node (03) at (1,0.866*6)    {4};
    \node (13) at (3,0.866*6)    {0};  
    \node (23) at (5,0.866*6)    {0};     
    \node (33) at (7,0.866*6)    {4};    
    \node (04) at (0,0.866*8)    {5};
    \node (14) at (2,0.866*8)    {7};  
    \node (24) at (4,0.866*8)    {5};     
    \node (34) at (6,0.866*8)    {7};    
    \node (05) at (1,0.866*10)    {2};
    \node (15) at (3,0.866*10)    {6};  
    \node (25) at (5,0.866*10)    {6};     
    \node (35) at (7,0.866*10)    {2};    
    \node (06) at (0,0.866*12)    {1};
    \node (16) at (2,0.866*12)    {3};  
    \node (26) at (4,0.866*12)    {1};     
    \node (36) at (6,0.866*12)    {3};    
    \node (07) at (1,0.866*14)    {4};
    \node (17) at (3,0.866*14)    {0};  
    \node (27) at (5,0.866*14)    {0};     
    \node (37) at (7,0.866*14)    {4};    
    \draw[dotted] (00) -- (-1,0);
    \draw (00) -- (10);
    \draw (10) -- (20);
    \draw (20) -- (30);
    \draw[dotted] (30) -- (7,0);
    \draw[dotted] (01) -- (0,0.866*2);
    \draw (01) -- (11);  
    \draw (11) -- (21);  
    \draw (21) -- (31);  
    \draw[dotted] (31) -- (8,0.866*2);
    \draw[dotted] (02) -- (-1,0.866*4);
    \draw (02) -- (12);
    \draw (12) -- (22);
    \draw (22) -- (32);  
    \draw[dotted] (32) -- (7,0.866*4);
    \draw[dotted] (03) -- (0,0.866*6);
    \draw (03) -- (13);  
    \draw (13) -- (23);  
    \draw (23) -- (33);    
    \draw[dotted] (33) -- (8,0.866*6);
    \draw[dotted] (04) -- (-1,0.866*8);
    \draw (04) -- (14);  
    \draw (14) -- (24);  
    \draw (24) -- (34);
    \draw[dotted] (34) -- (7,0.866*8);
    \draw[dotted] (05) -- (0,0.866*10);
    \draw (05) -- (15);  
    \draw (15) -- (25);  
    \draw (25) -- (35);  
    \draw[dotted] (35) -- (8,0.866*10);
    \draw[dotted] (06) -- (-1,0.866*12);
    \draw (06) -- (16);  
    \draw (16) -- (26);  
    \draw (26) -- (36);  
    \draw[dotted] (36) -- (7,0.866*12);
    \draw[dotted] (07) -- (0,0.866*14);
    \draw (07) -- (17);  
    \draw (17) -- (27);  
    \draw (27) -- (37);  
    \draw[dotted] (37) -- (8,0.866*14);
    \draw[dotted] (00) -- (-.5,-.866);  
    \draw[dotted] (00) -- (.5,-.866);  
    \draw[dotted] (10) -- (1.5,-.866);  
    \draw[dotted] (10) -- (2.5,-.866);   
    \draw[dotted] (20) -- (3.5,-.866); 
    \draw[dotted] (20) -- (4.5,-.866); 
    \draw[dotted] (30) -- (5.5,-.866);   
    \draw[dotted] (30) -- (6.5,-.866);  
    \draw (00) -- (01);  
    \draw (10) -- (01);  
    \draw (10) -- (11);   
    \draw (20) -- (11); 
    \draw (20) -- (21); 
    \draw (30) -- (21);   
    \draw (30) -- (31); 
    \draw (02) -- (01);  
    \draw (12) -- (01);  
    \draw (12) -- (11);   
    \draw (22) -- (11); 
    \draw (22) -- (21); 
    \draw (32) -- (21);   
    \draw (32) -- (31); 
    \draw (02) -- (03);  
    \draw (12) -- (03);  
    \draw (12) -- (13);   
    \draw (22) -- (13); 
    \draw (22) -- (23); 
    \draw (32) -- (23);   
    \draw (32) -- (33);  
    \draw (04) -- (03);  
    \draw (14) -- (03);  
    \draw (14) -- (13);   
    \draw (24) -- (13); 
    \draw (24) -- (23); 
    \draw (34) -- (23);   
    \draw (34) -- (33); 
    \draw (04) -- (05);  
    \draw (14) -- (05);  
    \draw (14) -- (15);   
    \draw (24) -- (15); 
    \draw (24) -- (25); 
    \draw (34) -- (25);   
    \draw (34) -- (35);  
    \draw (06) -- (05);  
    \draw (16) -- (05);  
    \draw (16) -- (15);   
    \draw (26) -- (15); 
    \draw (26) -- (25); 
    \draw (36) -- (25);   
    \draw (36) -- (35); 
    \draw (06) -- (07);  
    \draw (16) -- (07);  
    \draw (16) -- (17);   
    \draw (26) -- (17); 
    \draw (26) -- (27); 
    \draw (36) -- (27);   
    \draw (36) -- (37); 
    \draw[dotted] (07) -- (.5,15*.866);  
    \draw[dotted] (07) -- (1.5,15*.866);  
    \draw[dotted] (17) -- (2.5,15*.866);  
    \draw[dotted] (17) -- (3.5,15*.866);   
    \draw[dotted] (27) -- (4.5,15*.866); 
    \draw[dotted] (27) -- (5.5,15*.866); 
    \draw[dotted] (37) -- (6.5,15*.866);   
    \draw[dotted] (37) -- (7.5,15*.866);  
\end{tikzpicture}  \hfill \,
\end{center}
Moreover, the Apollonian group acts transitively on these lifts.
\end{lemma}

\begin{proof}
In this case, there are four possible liftings which satisfy Equations \eqref{DiamondQuadratic} and \eqref{Trapezoid} mod 8, but two of them are ruled out by the condition that $a+b \equiv 4 \bmod 8$ for tangent odd curvatures $a$, $b$, as verified in the proof of Lemma \ref{TriangleModLemma1}. 

We now check that the Apollonian group maps between configurations of these two types. As in the previous proof, if circles of curvatures $a$, $b$, $c$, $d$ are arranged in a diamond, then the generator which fixes $a$, $b$, $c$ will transform $d$ to $d'=6a+2b+6c-d$. Then the generators corresponding to triangles with one odd and two even curvatures will fix each mod 8 configuration shown above. The generators corresponding to triangles with one even and two odd curvatures will transform one mod 8 configuration to the other. 
\end{proof}

In particular, any primitive integral triangular packing with even and odd curvatures contains both configurations mod 8. 

Via the Chinese remainder theorem, each possible configuration mod 8 listed in Lemmas \ref{TriangleModLemma4} and \ref{TriangleModLemma5} combines with each configuration mod 3 listed in Lemma \ref{TriangleModLemma2} to form a unique mod 24 configuration. Some of these configurations are ruled out by Lemma \ref{TriangleModLemma3}, but this only depends on the curvatures mod 12. For each configuration mod 12 listed in Proposition \ref{TriangleMod}, the Apollonian group acts transitively on the set of lifts to mod 24. See Figure \ref{TrianglePartialEvenFigure} for the possible configurations mod 24 which contain both odd and even curvatures.

Next we parametrize the curvatures of all circles tangent to a fixed circle in a triangular packing, beginning with the Ford circles.

\begin{proposition} \label{TriangleFord}
The generalized Ford circles are parametrized by $x, y \in \Z$ with $y \geq 0$ and $\gcd(x,y)=1$. Their inversive coordinates are as follows:
\begin{align*}    
& c_\alpha(x,y)=(12x^2,y^2,2\sqrt{3}xy, 1) \hspace{3pt} \text{ if } 3 \nmid y \\
& c_\beta(x,y)=(4x^2,\frac{y^2}{3},\frac{2xy}{\sqrt{3}}, 1) \hspace{13pt} \text{ if } 3\mid y 
\end{align*}
Each circle is tangent to the real axis at $\frac{2\sqrt{3}x}{y}$.
\end{proposition} 

\begin{proposition} \label{TriangularParametrization}
    Suppose a triangular packing contains a circle of curvature $a$, and three adjacent circles forming a diamond  have curvatures $b$, $c$, $d$ as in Equation \eqref{DiamondQuadratic}. Then all circles tangent to the circle of curvature $a$ are parametrized by $x, y \in \Z$ with $y \geq 0$, $\gcd(x,y)=1$. Their curvatures are as follows:
\begin{align*}
& Q_\alpha(x,y)-a=3\left(a+b\right)x^2 +  \left(-3a-b-3c+d\right)xy + \left(a+c\right)y^2 -a\hspace{16pt} \text{ if } 3 \nmid y \\
& Q_\beta(x,y)-a= \left(a+b\right)x^2 +  \left(\frac{-3a-b-3c+d}{3}\right)xy + \left(\frac{a+c}{3}\right)y^2 -a \hspace{3pt} \text{ if } 3\mid y 
\end{align*}
\end{proposition}

Propositions \ref{TriangleFord} and \ref{TriangularParametrization} can be proved similarly to those in previous sections. Using Equation \eqref{DiamondQuadratic}, we find that $Q_\alpha$ has discriminant $-12a^2$ and $Q_\beta$ has discriminant $-\frac{4}{3}a^2$.

We now define quadratic invariants $\chi_2$ for triangular packings. The quadratic forms $Q_\alpha$ and $Q_\beta$ both represent the integers $a+b$, $a+c$ and $a+d$, up to factors of 3, and if $3\mid a$ then at least one of these forms represents integers not divisible by 3. Thus because $\gcd(a,b,c,d)=1$, we can conclude that $Q_\alpha(x,y)$ or $Q_\beta(x,y)$ represents an integer $\rho$ coprime to $a$. We may take $y\geq 0$, $\gcd(x,y)=1$, and $3 \nmid y$ if $Q_\alpha$ is used or $3\mid y$ if $Q_\beta$ is used. 

\begin{proposition} \label{TriangleNode}
Let $a$ be a curvature in a primitive integral triangular packing, and suppose that all curvatures of circles tangent to $a$ are represented by the quadratic forms $Q_\alpha$ and $Q_\beta$ as in the previous proposition. Set $a'=a/2$ if $a\equiv 2 \bmod 4$, and $a'=a$ otherwise. Let $\rho_1$, $\rho_2$ be two integers coprime to $a$ which are represented by either $Q_{\alpha}$ or $Q_{\beta}$. If $3\mid a$, or if $a\equiv 1, 2, 11, 13, 22, 23 \bmod 24$, then 
\begin{equation*}
\res{\rho_1}{a'} = \res{\rho_2}{a'}
\end{equation*}
If $a\equiv 5, 7, 10, 14, 17, 19 \bmod 24$ then the equality holds if and only if $\rho_1$, $\rho_2$ are represented by the same quadratic form; if $\rho_1$ is represented by $Q_\alpha$ and $\rho_2$ is represented by $Q_\beta$, then 
\begin{equation*}
\res{\rho_1}{a'} = -\res{\rho_2}{a'}
\end{equation*}
\end{proposition}

Note that when $a\equiv 4, 8, 16, 20 \bmod 24$, we need information about its factorization to determine whether $\res{*}{a}$ takes the same or opposite values for the two quadratic forms. 

\begin{proof}
If $\rho_1$, $\rho_2$ are represented by the same quadratic form, then $\res{\rho_1}{a'}=\res{\rho_2}{a'}$ by Lemma \ref{QuadFormKronecker}, since the discriminants are divisible by $a'$ and a sufficient power of $2$. 

Suppose that $\rho_1$ and $\rho_2$ are represented by $Q_\alpha$ and $Q_\beta$ respectively. This implies that $3 \nmid a$; if $3\mid a$, then one of the forms $Q_\alpha(x,y)$, $Q_\beta(x, 3y)$ only represents multiples of 3. Since $\rho_2$ is represented by $Q_\beta$, $3\rho_2$ is represented by $Q_\alpha$, and 
$$\res{\rho_1}{a'}=\res{3\rho_2}{a'}=\res{3}{a'}\res{\rho_2}{a'}$$
We have $\res{3}{a'}=1$ if $a\equiv 1, 2, 11, 13, 22, 23 \bmod 24$, and $\res{3}{a'}=-1$ if $a\equiv 5, 7, 10, 14, 17, 19 \bmod 24$, which completes the proof.
\end{proof}

We now define the quadratic invariant $\chi_2$ for triangular packings of type $(1)$ and type $(3, 11)$, where the invariant will be well-defined across the entire packing. 

In type $(1)$, for a circle $\C$ of curvature $a$, we set $\chi_2(\C)=\res{\rho}{a}$. In type $(3, 11)$, we set 
$$\chi_2(\C)=\begin{cases}
    \left(\frac{\rho}{a}\right) &\text{if } a \equiv 11 \bmod 12;\\
    -\left(\frac{\rho}{a}\right) &\text{if } a \equiv 3 \bmod 12;
    \end{cases}$$
Where $\rho$ is any curvature represented by $Q_\alpha$ or $Q_\beta$, satisfying $\gcd(\rho, a)=1$. Proposition \ref{TriangleNode} guarantees that this expression is independent of the choice of $\rho$.  

\begin{proposition} \label{TriangleEdge}
Suppose that $\C_a$ and $\C_b$ are tangent circles with coprime curvatures $a$, $b$ in a primitive integral triangular packing of type $(1)$ or $(3,11)$. Then 
\begin{equation*}
\chi_2(\C_a)=\chi_2(\C_b)
\end{equation*}
\end{proposition}
\begin{proof}
We will use $\rho=a+b$ to compute $\chi_2(\C_a)$ and $\chi_2(\C_b)$. If $a\equiv b \equiv 1 \bmod 12$, then
\begin{equation*}
\chi_2(\C_a)=\res{a+b}{a} = \res{b}{a} = \res{a}{b} = \res{a+b}{b}=\chi_2(\C_b)
\end{equation*}
If $a\equiv 3$, $b \equiv 11 \bmod 12$, then
\begin{equation*}
\chi_2(\C_a)=-\res{a+b}{a} = -\res{b}{a} = \res{a}{b} = \res{a+b}{b}=\chi_2(\C_b)
\end{equation*}
These are the only possible cases in types $(1)$ and $(3,11)$.
\end{proof}

To show that $\chi_2(\mathcal{C})$ is identical for all circles in the packing, we must find a path of coprime curvatures between any two circles. 

\begin{lemma} \label{TriangleSimultaneous}
Suppose that $a, b, c, d$ are curvatures in a triangular packing, with tangency relations as in Equation \eqref{DiamondQuadratic}. Then all circles simultaneously tangent to $a,c$ in the packing are parametrized by the following formula:
\begin{equation*}
    a \left(3 n^2-3 n\right)+c \left(3 n^2-3 n\right)+b (1-n)+d n
\end{equation*} 
\end{lemma}
The proof is similar to that of Lemma \ref{OctSimultaneous}. 

\begin{lemma} \label{TriangleInsert}
Given a pair of tangent circles of curvatures $a$, $c$ in a primitive integral square packing, there exists a circle of curvature $k$ tangent to $a$, $c$ such that $\gcd(a, k)\mid 3$ and $\gcd(c, k)\mid 3$. 
\end{lemma}
\begin{proof}
 The quadratic polynomial in Lemma \ref{TriangleSimultaneous} represents the integers $b$, $d$, and $6a+6c-b+2d$. Since the packing is primitive, $\gcd(a,b,c,d)=1$, and thus $\gcd(a, b, d, 6a+6c-b+2d)\mid 6$ and $\gcd(c,b,d, 6a+6c-b+2d)\mid 6$. By the Chinese remainder theorem. the quadratic polynomial represents an integer $k$ whose greatest common divisor with $a$ and $c$ divides $6$. Moreover, if $a$ or $b$ is even, then by Lemma \ref{TriangleModLemma1}, $b$ or $d$ is odd, and we may take $k$ odd. Then $\gcd(a, k)\mid 3$ and $\gcd(c, k)\mid 3$. 
\end{proof}

\begin{corollary} \label{TrianglePath}
Let $\mathcal{C},\mathcal{C}' \in \mathcal{A}$ be two circles in a primitive square circle packing. Then there exists a path of circles $\mathcal{C}_1,\mathcal{C}_2,...,\mathcal{C}_k$ such that \begin{enumerate}
    \item $\mathcal{C}_1=\mathcal{C}$ and $\mathcal{C}_k=\mathcal{C}'$;
    \item $\mathcal{C}_i$ is tangent to $\mathcal{C}_{i+1}$ for all $1 \leq i \leq k-1$;
    \item The curvatures of $\mathcal{C}_i$ and $\mathcal{C}_{i+1}$ are coprime for all $1 \leq i \leq k-1$.
\end{enumerate}
\end{corollary}
The proof is the same as that of Proposition \ref{OctPath}, using Lemma \ref{TriangleInsert} to insert a circle between two circles in the path if the curvatures are not coprime. If two tangent circles have curvatures both divisible by 3, Proposition \ref{TriangleModLemma2} implies that the two circle inserted between them will have curvature not divisible by 3, ensuring coprimality.  

Proposition \ref{TriangleEdge} and Corollary \ref{TrianglePath} imply the following result. 

\begin{proposition} \label{TriangleChi}
The value of $\chi_2$ is constant across all circles in a fixed primitive triangular circle packing $\mathcal{A}$ of type $(1)$ or $(3,11)$.
\end{proposition}

Based on this Proposition, we will refer to $\chi_2(\mathcal{A})$ as a quadratic invariant of the packing. This is used to prove the main result of this section.

\begin{theorem} \label{TriangleObstruction}
In a primitive integral triangular packing $\mathcal{A}$ of type $(1)$ with $\chi_2(\mathcal{A})=-1$, no integers of the form $n^2$ appear as curvatures. In a primitive integral triangular packing $\mathcal{A}$ of type $(3,11)$ with $\chi_2(\mathcal{A})=-1$, no integers of the form $3n^2$ appear as curvatures. In particular, the Local-Global Conjecture is false for these packings.
\end{theorem}
\begin{proof}
Given a circle of curvature $n^2$ in type $(1)$ or $3n^2$ in type $(3,11)$, we can find a tangent circle $\C$ of coprime curvature and compute $\chi_2(\C)=1$, a contradiction.
\end{proof}

\begin{figure}[h]
\begin{center} \, \hfill
\begin{tikzpicture}[scale=.5]
    \node (00) at (0,0)    {7};
    \node (10) at (2,0)    {7};  
    \node (20) at (4,0)    {7};     
    \node (30) at (6,0)    {7};    
    \node (01) at (1,0.866*2)    {7};
    \node (11) at (3,0.866*2)    {7};  
    \node (21) at (5,0.866*2)    {7};     
    \node (31) at (7,0.866*2)    {7};   
    \node (02) at (0,0.866*4)    {7};
    \node (12) at (2,0.866*4)    {7};  
    \node (22) at (4,0.866*4)    {7};     
    \node (32) at (6,0.866*4)    {7};     
    \node (03) at (1,0.866*6)    {7};
    \node (13) at (3,0.866*6)    {7};  
    \node (23) at (5,0.866*6)    {7};     
    \node (33) at (7,0.866*6)    {7};    
    \draw[dotted] (00) -- (-1,0);
    \draw[->] (00) -- (10);
    \draw[->] (10) -- (20);
    \draw[->] (20) -- (30);
    \draw[dotted] (30) -- (7,0);
    \draw[dotted] (01) -- (0,0.866*2);
    \draw[->] (01) -- (11);  
    \draw[->] (11) -- (21);  
    \draw[->] (21) -- (31);  
    \draw[dotted] (31) -- (8,0.866*2);
    \draw[dotted] (02) -- (-1,0.866*4);
    \draw[->] (02) -- (12);
    \draw[->] (12) -- (22);
    \draw[->] (22) -- (32);
    \draw[dotted] (32) -- (7,0.866*4);
    \draw[dotted] (03) -- (0,0.866*6);
    \draw[->] (03) -- (13);  
    \draw[->] (13) -- (23);  
    \draw[->] (23) -- (33);  
    \draw[dotted] (33) -- (8,0.866*6);
    \draw[dotted] (00) -- (-.5,-.866);  
    \draw[dotted] (00) -- (.5,-.866);  
    \draw[dotted] (10) -- (1.5,-.866);  
    \draw[dotted] (10) -- (2.5,-.866);   
    \draw[dotted] (20) -- (3.5,-.866); 
    \draw[dotted] (20) -- (4.5,-.866); 
    \draw[dotted] (30) -- (5.5,-.866);   
    \draw[dotted] (30) -- (6.5,-.866);   
    \draw[<-] (00) -- (01);  
    \draw[->] (10) -- (01);  
    \draw[<-] (10) -- (11);   
    \draw[->] (20) -- (11); 
    \draw[<-] (20) -- (21); 
    \draw[->] (30) -- (21);   
    \draw[<-] (30) -- (31);  
    \draw[<-] (02) -- (01);  
    \draw[->] (12) -- (01);  
    \draw[<-] (12) -- (11);   
    \draw[->] (22) -- (11); 
    \draw[<-] (22) -- (21); 
    \draw[->] (32) -- (21);   
    \draw[<-] (32) -- (31); 
    \draw[<-] (02) -- (03);  
    \draw[->] (12) -- (03);  
    \draw[<-] (12) -- (13);   
    \draw[->] (22) -- (13); 
    \draw[<-] (22) -- (23); 
    \draw[->] (32) -- (23);   
    \draw[<-] (32) -- (33); 
    \draw[dotted] (03) -- (.5,7*.866);  
    \draw[dotted] (03) -- (1.5,7*.866);  
    \draw[dotted] (13) -- (2.5,7*.866);  
    \draw[dotted] (13) -- (3.5,7*.866);   
    \draw[dotted] (23) -- (4.5,7*.866); 
    \draw[dotted] (23) -- (5.5,7*.866); 
    \draw[dotted] (33) -- (6.5,7*.866);   
    \draw[dotted] (33) -- (7.5,7*.866);   
\end{tikzpicture} \hfill 
\begin{tikzpicture}[scale=.5]
    \node (00) at (0,0)    {5};
    \node (10) at (2,0)    {9};  
    \node (20) at (4,0)    {9};     
    \node (30) at (6,0)    {5};    
    \node (01) at (1,0.866*2)    {9};
    \node (11) at (3,0.866*2)    {5};  
    \node (21) at (5,0.866*2)    {9};     
    \node (31) at (7,0.866*2)    {9};    
    \node (02) at (0,0.866*4)    {5};
    \node (12) at (2,0.866*4)    {9};  
    \node (22) at (4,0.866*4)    {9};     
    \node (32) at (6,0.866*4)    {5};     
    \node (03) at (1,0.866*6)    {9};
    \node (13) at (3,0.866*6)    {5};  
    \node (23) at (5,0.866*6)    {9};     
    \node (33) at (7,0.866*6)    {9};    
    \draw[dotted] (00) -- (-1,0);
    \draw[<->] (00) -- (10);
    \draw[dashed] (10) -- (20);
    \draw (20) -- (30);
    \draw[dotted] (30) -- (7,0);
    \draw[dotted] (01) -- (0,0.866*2);
    \draw (01) -- (11);  
    \draw[<->] (11) -- (21);  
    \draw[dashed] (21) -- (31);  
    \draw[dotted] (31) -- (8,0.866*2);
    \draw[dotted] (02) -- (-1,0.866*4);
    \draw[<->] (02) -- (12);
    \draw[dashed] (12) -- (22);
    \draw (22) -- (32); 
    \draw[dotted] (32) -- (7,0.866*4);
    \draw[dotted] (03) -- (0,0.866*6);
    \draw (03) -- (13);  
    \draw[<->] (13) -- (23);  
    \draw[dashed] (23) -- (33);  
    \draw[dotted] (33) -- (8,0.866*6);
    \draw[dotted] (00) -- (-.5,-.866);  
    \draw[dotted] (00) -- (.5,-.866);  
    \draw[dotted] (10) -- (1.5,-.866);  
    \draw[dotted] (10) -- (2.5,-.866);   
    \draw[dotted] (20) -- (3.5,-.866); 
    \draw[dotted] (20) -- (4.5,-.866); 
    \draw[dotted] (30) -- (5.5,-.866);   
    \draw[dotted] (30) -- (6.5,-.866);   
    \draw (00) -- (01);  
    \draw[dashed] (10) -- (01);  
    \draw[<->] (10) -- (11);   
    \draw (20) -- (11); 
    \draw[dashed] (20) -- (21); 
    \draw[<->] (30) -- (21);   
    \draw (30) -- (31); 
    \draw (02) -- (01);  
    \draw[dashed] (12) -- (01);  
    \draw[<->] (12) -- (11);   
    \draw (22) -- (11); 
    \draw[dashed] (22) -- (21); 
    \draw[<->] (32) -- (21);   
    \draw (32) -- (31); 
    \draw (02) -- (03);  
    \draw[dashed] (12) -- (03);  
    \draw[<->] (12) -- (13);   
    \draw (22) -- (13); 
    \draw[dashed] (22) -- (23); 
    \draw[<->] (32) -- (23);   
    \draw (32) -- (33); 
    \draw[dotted] (03) -- (.5,7*.866);  
    \draw[dotted] (03) -- (1.5,7*.866);  
    \draw[dotted] (13) -- (2.5,7*.866);  
    \draw[dotted] (13) -- (3.5,7*.866);   
    \draw[dotted] (23) -- (4.5,7*.866); 
    \draw[dotted] (23) -- (5.5,7*.866); 
    \draw[dotted] (33) -- (6.5,7*.866);   
    \draw[dotted] (33) -- (7.5,7*.866);  
\end{tikzpicture} \hfill \,
\end{center}
\caption{Pattern of Kronecker Symbols in Triangular Packings of Types $(7)$, $(5,9)$}
\label{TrianglePartialFigure}
\end{figure}
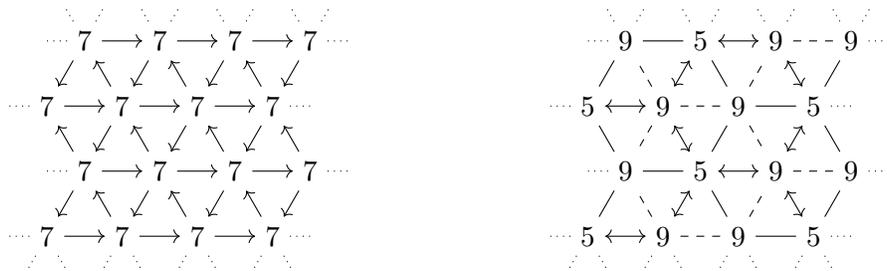

\begin{figure}[H]
\begin{center} \, \hfill \hfill
\begin{tikzpicture}[scale=.45]
    \node (00) at (0,0)    {17};
    \node (10) at (2,0)    {11};  
    \node (20) at (4,0)    {17};     
    \node (30) at (6,0)    {11};
    \node (01) at (1,0.866*2)    {2};
    \node (11) at (3,0.866*2)    {2};  
    \node (21) at (5,0.866*2)    {14};     
    \node (31) at (7,0.866*2)    {14};    
    \node (02) at (0,0.866*4)    {5};
    \node (12) at (2,0.866*4)    {23};  
    \node (22) at (4,0.866*4)    {5};     
    \node (32) at (6,0.866*4)    {23};     
    \node (03) at (1,0.866*6)    {8};
    \node (13) at (3,0.866*6)    {8};  
    \node (23) at (5,0.866*6)    {20};     
    \node (33) at (7,0.866*6)    {20};    
    \node (04) at (0,0.866*8)    {5};
    \node (14) at (2,0.866*8)    {23};  
    \node (24) at (4,0.866*8)    {5};     
    \node (34) at (6,0.866*8)    {23};    
    \node (05) at (1,0.866*10)    {2};
    \node (15) at (3,0.866*10)    {2};  
    \node (25) at (5,0.866*10)    {14};     
    \node (35) at (7,0.866*10)    {14};    
    \node (06) at (0,0.866*12)    {17};
    \node (16) at (2,0.866*12)    {11};  
    \node (26) at (4,0.866*12)    {17};     
    \node (36) at (6,0.866*12)    {11};    
    \node (07) at (1,0.866*14)    {8};
    \node (17) at (3,0.866*14)    {8};  
    \node (27) at (5,0.866*14)    {20};     
    \node (37) at (7,0.866*14)    {20};    
    \draw[dotted] (00) -- (-1,0);
    \draw[<->] (00) -- (10);
    \draw[<->] (10) -- (20);
    \draw (20) -- (30);
    \draw (30) -- (7.5,0);
    \draw[dashed] (01) -- (-.5,0.866*2);
    \draw[dashed] (01) -- (11);  
    \draw[dashed] (11) -- (21);  
    \draw[dashed] (21) -- (31);  
    \draw[dotted] (31) -- (8,0.866*2);
    \draw[dotted] (02) -- (-1,0.866*4);
    \draw (02) -- (12);
    \draw (12) -- (22);
    \draw[<->] (22) -- (32);  
    \draw[<->] (32) -- (7.5,0.866*4);
    \draw[dashed] (03) -- (-.5,0.866*6);
    \draw[dashed] (03) -- (13);  
    \draw[dashed] (13) -- (23);  
    \draw[dashed] (23) -- (33);    
    \draw[dotted] (33) -- (8,0.866*6);
    \draw[dotted] (04) -- (-1,0.866*8);
    \draw (04) -- (14);  
    \draw (14) -- (24);  
    \draw[<->] (24) -- (34);
    \draw[<->] (34) -- (7.5,0.866*8);
    \draw[dashed] (05) -- (-.5,0.866*10);
    \draw[dashed] (05) -- (15);  
    \draw[dashed] (15) -- (25);  
    \draw[dashed] (25) -- (35);  
    \draw[dotted] (35) -- (8,0.866*10);
    \draw[dotted] (06) -- (-1,0.866*12);
    \draw[<->] (06) -- (16);  
    \draw[<->] (16) -- (26);  
    \draw (26) -- (36);  
    \draw (36) -- (7.5,0.866*12);
    \draw[dashed] (07) -- (-.5,0.866*14);
    \draw[dashed] (07) -- (17);  
    \draw[dashed] (17) -- (27);  
    \draw[dashed] (27) -- (37);  
    \draw[dotted] (37) -- (8,0.866*14);
    \draw[dotted] (00) -- (-.5,-.866);  
    \draw (00) -- (.75,-1.5*.866);  
    \draw[{Circle}->] (10) -- (1.25,-1.5*.866);  
    \draw[{Circle}->] (10) -- (2.75,-1.5*.866);   
    \draw (20) -- (3.25,-1.5*.866); 
    \draw[<->] (20) -- (4.75,-1.5*.866); 
    \draw[{Circle}-] (30) -- (5.25,-1.5*.866);
    \draw[{Circle}-] (30) -- (6.75,-1.5*.866); 
    \draw (00) -- (01);  
    \draw[->] (10) -- (01);  
    \draw[->] (10) -- (11);   
    \draw (20) -- (11); 
    \draw[<->] (20) -- (21); 
    \draw (30) -- (21);   
    \draw (30) -- (31); 
    \draw[->] (02) -- (01);  
    \draw (12) -- (01);  
    \draw (12) -- (11);   
    \draw[->] (22) -- (11); 
    \draw[<-] (22) -- (21); 
    \draw[->] (32) -- (21);   
    \draw[->] (32) -- (31); 
    \draw[<->] (02) -- (03);  
    \draw[{Circle}-] (12) -- (03);  
    \draw[{Circle}-] (12) -- (13);   
    \draw[<->] (22) -- (13); 
    \draw (22) -- (23); 
    \draw[{Circle}->] (32) -- (23);   
    \draw[{Circle}->] (32) -- (33);  
    \draw[<->] (04) -- (03);  
    \draw[{Circle}-] (14) -- (03);  
    \draw[{Circle}-] (14) -- (13);   
    \draw[<->] (24) -- (13); 
    \draw (24) -- (23); 
    \draw[{Circle}->] (34) -- (23);   
    \draw[{Circle}->] (34) -- (33); 
    \draw[->] (04) -- (05);  
    \draw (14) -- (05);  
    \draw (14) -- (15);   
    \draw[->] (24) -- (15); 
    \draw[<-] (24) -- (25); 
    \draw[->] (34) -- (25);   
    \draw[->] (34) -- (35);  
    \draw (06) -- (05);  
    \draw[->] (16) -- (05);  
    \draw[->] (16) -- (15);   
    \draw (26) -- (15); 
    \draw[<->] (26) -- (25); 
    \draw (36) -- (25);   
    \draw (36) -- (35); 
    \draw (06) -- (07);  
    \draw[{Circle}->] (16) -- (07);  
    \draw[{Circle}->] (16) -- (17);   
    \draw (26) -- (17); 
    \draw[<->] (26) -- (27); 
    \draw[{Circle}-] (36) -- (27);   
    \draw[{Circle}-] (36) -- (37); 
    \draw[dotted] (07) -- (.5,15*.866);  
    \draw[dotted] (07) -- (1.5,15*.866);  
    \draw[dotted] (17) -- (2.5,15*.866);  
    \draw[dotted] (17) -- (3.5,15*.866);   
    \draw[dotted] (27) -- (4.5,15*.866); 
    \draw[dotted] (27) -- (5.5,15*.866); 
    \draw[dotted] (37) -- (6.5,15*.866);   
    \draw[dotted] (37) -- (7.5,15*.866);  
\end{tikzpicture}  \hfill
\begin{tikzpicture}[scale=.45]
    \node (00) at (0,0)    {17};
    \node (10) at (2,0)    {11};  
    \node (20) at (4,0)    {17};     
    \node (30) at (6,0)    {11};
    \node (01) at (1,0.866*2)    {2};
    \node (11) at (3,0.866*2)    {14};  
    \node (21) at (5,0.866*2)    {14};     
    \node (31) at (7,0.866*2)    {2};    
    \node (02) at (0,0.866*4)    {5};
    \node (12) at (2,0.866*4)    {23};  
    \node (22) at (4,0.866*4)    {5};     
    \node (32) at (6,0.866*4)    {23};     
    \node (03) at (1,0.866*6)    {20};
    \node (13) at (3,0.866*6)    {8};  
    \node (23) at (5,0.866*6)    {8};     
    \node (33) at (7,0.866*6)    {20};    
    \node (04) at (0,0.866*8)    {5};
    \node (14) at (2,0.866*8)    {23};  
    \node (24) at (4,0.866*8)    {5};     
    \node (34) at (6,0.866*8)    {23};    
    \node (05) at (1,0.866*10)    {2};
    \node (15) at (3,0.866*10)    {14};  
    \node (25) at (5,0.866*10)    {14};     
    \node (35) at (7,0.866*10)    {2};    
    \node (06) at (0,0.866*12)    {17};
    \node (16) at (2,0.866*12)    {11};  
    \node (26) at (4,0.866*12)    {17};     
    \node (36) at (6,0.866*12)    {11};    
    \node (07) at (1,0.866*14)    {20};
    \node (17) at (3,0.866*14)    {8};  
    \node (27) at (5,0.866*14)    {8};     
    \node (37) at (7,0.866*14)    {20};    
    \draw[dotted] (00) -- (-1,0);
    \draw[<->] (00) -- (10);
    \draw[<->] (10) -- (20);
    \draw (20) -- (30);
    \draw (30) -- (7.5,0);
    \draw[dashed] (01) -- (-.5,0.866*2);
    \draw[dashed] (01) -- (11);  
    \draw[dashed] (11) -- (21);  
    \draw[dashed] (21) -- (31);  
    \draw[dotted] (31) -- (8,0.866*2);
    \draw[dotted] (02) -- (-1,0.866*4);
    \draw (02) -- (12);
    \draw (12) -- (22);
    \draw[<->] (22) -- (32);  
    \draw[<->] (32) -- (7.5,0.866*4);
    \draw[dashed] (03) -- (-.5,0.866*6);
    \draw[dashed] (03) -- (13);  
    \draw[dashed] (13) -- (23);  
    \draw[dashed] (23) -- (33);    
    \draw[dotted] (33) -- (8,0.866*6);
    \draw[dotted] (04) -- (-1,0.866*8);
    \draw (04) -- (14);  
    \draw (14) -- (24);  
    \draw[<->] (24) -- (34);
    \draw[<->] (34) -- (7.5,0.866*8);
    \draw[dashed] (05) -- (-.5,0.866*10);
    \draw[dashed] (05) -- (15);  
    \draw[dashed] (15) -- (25);  
    \draw[dashed] (25) -- (35);  
    \draw[dotted] (35) -- (8,0.866*10);
    \draw[dotted] (06) -- (-1,0.866*12);
    \draw[<->] (06) -- (16);  
    \draw[<->] (16) -- (26);  
    \draw (26) -- (36);  
    \draw (36) -- (7.5,0.866*12);
    \draw[dashed] (07) -- (-.5,0.866*14);
    \draw[dashed] (07) -- (17);  
    \draw[dashed] (17) -- (27);  
    \draw[dashed] (27) -- (37);  
    \draw[dotted] (37) -- (8,0.866*14);
    \draw[dotted] (00) -- (-.5,-.866);  
    \draw (00) -- (.75,-1.5*.866);  
    \draw[{Circle}->] (10) -- (1.25,-1.5*.866);  
    \draw[{Circle}->] (10) -- (2.75,-1.5*.866);   
    \draw (20) -- (3.25,-1.5*.866); 
    \draw[<->] (20) -- (4.75,-1.5*.866); 
    \draw[{Circle}-] (30) -- (5.25,-1.5*.866); 
    \draw[{Circle}-] (30) -- (6.75,-1.5*.866);
    \draw (00) -- (01);  
    \draw[->] (10) -- (01);  
    \draw[<->] (10) -- (11);   
    \draw (20) -- (11); 
    \draw[<->] (20) -- (21); 
    \draw (30) -- (21);   
    \draw[<-] (30) -- (31); 
    \draw[->] (02) -- (01);  
    \draw (12) -- (01);  
    \draw[<-] (12) -- (11);   
    \draw[->] (22) -- (11); 
    \draw[<-] (22) -- (21); 
    \draw[->] (32) -- (21);   
    \draw[<->] (32) -- (31); 
    \draw[<->] (02) -- (03);  
    \draw[{Circle}-] (12) -- (03);  
    \draw[{Circle}-] (12) -- (13);   
    \draw[<->] (22) -- (13); 
    \draw (22) -- (23); 
    \draw[{Circle}->] (32) -- (23);   
    \draw[{Circle}->] (32) -- (33);  
    \draw[<->] (04) -- (03);  
    \draw[{Circle}-] (14) -- (03);  
    \draw[{Circle}-] (14) -- (13);   
    \draw[<->] (24) -- (13); 
    \draw (24) -- (23); 
    \draw[{Circle}->] (34) -- (23);   
    \draw[{Circle}->] (34) -- (33); 
    \draw[->] (04) -- (05);  
    \draw (14) -- (05);  
    \draw[<-] (14) -- (15);   
    \draw[->] (24) -- (15); 
    \draw[<-] (24) -- (25); 
    \draw[->] (34) -- (25);   
    \draw[<->] (34) -- (35);  
    \draw (06) -- (05);  
    \draw[->] (16) -- (05);  
    \draw[<->] (16) -- (15);   
    \draw (26) -- (15); 
    \draw[<->] (26) -- (25); 
    \draw (36) -- (25);   
    \draw[<-] (36) -- (35); 
    \draw (06) -- (07);  
    \draw[{Circle}->] (16) -- (07);  
    \draw[{Circle}->] (16) -- (17);   
    \draw (26) -- (17); 
    \draw[<->] (26) -- (27); 
    \draw[{Circle}-] (36) -- (27);   
    \draw[{Circle}-] (36) -- (37); 
    \draw[dotted] (07) -- (.5,15*.866);  
    \draw[dotted] (07) -- (1.5,15*.866);  
    \draw[dotted] (17) -- (2.5,15*.866);  
    \draw[dotted] (17) -- (3.5,15*.866);   
    \draw[dotted] (27) -- (4.5,15*.866); 
    \draw[dotted] (27) -- (5.5,15*.866); 
    \draw[dotted] (37) -- (6.5,15*.866);   
    \draw[dotted] (37) -- (7.5,15*.866);  
\end{tikzpicture} \hfill \hfill \,

\, \\

\begin{tikzpicture}[scale=.45]
    \node (00) at (0,0)    {1};
    \node (10) at (2,0)    {3};  
    \node (20) at (4,0)    {9};     
    \node (30) at (6,0)    {19};    
    \node (40) at (8,0)    {9};     
    \node (50) at (10,0)    {3};   
    \node (60) at (12, 0)   {1};   
    \node (70) at (14, 0)   {3};
    \node (80) at (16, 0)   {9};
    \node (90) at (18, 0)   {19};
    \node (100) at (20, 0)   {9};
    \node (110) at (22, 0)   {3};
    \node (01) at (1,0.866*2)    {18};
    \node (11) at (3,0.866*2)    {22};  
    \node (21) at (5,0.866*2)    {6};     
    \node (31) at (7,0.866*2)    {18};    
    \node (41) at (9,0.866*2)    {10};     
    \node (51) at (11,0.866*2)    {6}; 
    \node (61) at (13,0.866*2)    {6};
    \node (71) at (15,0.866*2)    {10};  
    \node (81) at (17,0.866*2)    {18};     
    \node (91) at (19,0.866*2)    {6};    
    \node (101) at (21,0.866*2)    {22};     
    \node (111) at (23,0.866*2)    {18}; 
    \node (02) at (0,0.866*4)    {13};
    \node (12) at (2,0.866*4)    {15};  
    \node (22) at (4,0.866*4)    {21};     
    \node (32) at (6,0.866*4)    {7};    
    \node (42) at (8,0.866*4)    {21};     
    \node (52) at (10,0.866*4)    {15};   
    \node (62) at (12,0.866*4)    {13};
    \node (72) at (14,0.866*4)    {15};  
    \node (82) at (16,0.866*4)    {21};     
    \node (92) at (18,0.866*4)    {7};    
    \node (102) at (20,0.866*4)    {21};     
    \node (112) at (22,0.866*4)    {15};    
    \node (03) at (1,0.866*6)    {12};
    \node (13) at (3,0.866*6)    {16};  
    \node (23) at (5,0.866*6)    {0};     
    \node (33) at (7,0.866*6)    {12};    
    \node (43) at (9,0.866*6)    {4};     
    \node (53) at (11,0.866*6)    {0}; 
    \node (63) at (13,0.866*6)    {0};
    \node (73) at (15,0.866*6)    {4};  
    \node (83) at (17,0.866*6)    {12};     
    \node (93) at (19,0.866*6)    {0};    
    \node (103) at (21,0.866*6)    {16};     
    \node (113) at (23,0.866*6)    {12}; 
    \node (04) at (0,0.866*8)    {13};
    \node (14) at (2,0.866*8)    {15};  
    \node (24) at (4,0.866*8)    {21};     
    \node (34) at (6,0.866*8)    {7};    
    \node (44) at (8,0.866*8)    {21};     
    \node (54) at (10,0.866*8)    {15};   
    \node (64) at (12,0.866*8)    {13};
    \node (74) at (14,0.866*8)    {15};  
    \node (84) at (16,0.866*8)    {21};     
    \node (94) at (18,0.866*8)    {7};    
    \node (104) at (20,0.866*8)    {21};     
    \node (114) at (22,0.866*8)    {15};  
    \node (05) at (1,0.866*10)    {18};
    \node (15) at (3,0.866*10)    {22};  
    \node (25) at (5,0.866*10)    {6};     
    \node (35) at (7,0.866*10)    {18};    
    \node (45) at (9,0.866*10)    {10};     
    \node (55) at (11,0.866*10)    {6}; 
    \node (65) at (13,0.866*10)    {6};
    \node (75) at (15,0.866*10)    {10};  
    \node (85) at (17,0.866*10)    {18};     
    \node (95) at (19,0.866*10)    {6};    
    \node (105) at (21,0.866*10)    {22};     
    \node (115) at (23,0.866*10)    {18}; 
    \node (06) at (0,0.866*12)    {1};
    \node (16) at (2,0.866*12)    {3};  
    \node (26) at (4,0.866*12)    {9};     
    \node (36) at (6,0.866*12)    {19};    
    \node (46) at (8,0.866*12)    {9};     
    \node (56) at (10,0.866*12)    {3};   
    \node (66) at (12,0.866*12)    {1};
    \node (76) at (14,0.866*12)    {3};  
    \node (86) at (16,0.866*12)    {9};     
    \node (96) at (18,0.866*12)    {19};    
    \node (106) at (20,0.866*12)    {9};     
    \node (116) at (22,0.866*12)    {3};  
    \node (07) at (1,0.866*14)    {12};
    \node (17) at (3,0.866*14)    {16};  
    \node (27) at (5,0.866*14)    {0};     
    \node (37) at (7,0.866*14)    {12};    
    \node (47) at (9,0.866*14)    {4};     
    \node (57) at (11,0.866*14)    {0}; 
    \node (67) at (13,0.866*14)    {0};
    \node (77) at (15,0.866*14)    {4};  
    \node (87) at (17,0.866*14)    {12};     
    \node (97) at (19,0.866*14)    {0};    
    \node (107) at (21,0.866*14)    {16};     
    \node (117) at (23,0.866*14)    {12}; 
    \draw[dotted] (00) -- (-1,0);
    \draw[<->] (00) -- (10);
    \draw[dashed] (10) -- (20);
    \draw[<->] (20) -- (30);
    \draw (30) -- (40);  
    \draw[dashed] (40) -- (50);  
    \draw (50) -- (60);
    \draw (60) -- (70);
    \draw[dashed] (70) -- (80);
    \draw (80) -- (90);  
    \draw[<->] (90) -- (100);  
    \draw[dashed] (100) -- (110); 
    \draw[<->] (110) -- (23.5,0);
    \draw[dashed] (01) -- (-.5,0.866*2);
    \draw[dashed] (01) -- (11);  
    \draw[dashed] (11) -- (21);  
    \draw[dashed] (21) -- (31);  
    \draw[dashed] (31) -- (41);  
    \draw[dashed] (41) -- (51); 
    \draw[dashed] (51) -- (61);  
    \draw[dashed] (61) -- (71);  
    \draw[dashed] (71) -- (81);  
    \draw[dashed] (81) -- (91);  
    \draw[dashed] (91) -- (101); 
    \draw[dashed] (101) -- (111); 
    \draw[dotted] (111) -- (24,0.866*2);
    \draw[dotted] (02) -- (-1,0.866*4);
    \draw (02) -- (12);
    \draw[dashed] (12) -- (22);
    \draw (22) -- (32);
    \draw[<->] (32) -- (42);  
    \draw[dashed] (42) -- (52);   
    \draw[<->] (52) -- (62);
    \draw[<->] (62) -- (72);
    \draw[dashed] (72) -- (82);
    \draw[<->] (82) -- (92);  
    \draw (92) -- (102);   
    \draw[dashed] (102) -- (112);   
    \draw (112) -- (23.5,0.866*4);
    \draw[dashed] (03) -- (-.5,0.866*6);
    \draw[dashed] (03) -- (13);  
    \draw[dashed] (13) -- (23);  
    \draw[dashed] (23) -- (33);  
    \draw[dashed] (33) -- (43);  
    \draw[dashed] (43) -- (53); 
    \draw[dashed] (53) -- (63);
    \draw[dashed] (63) -- (73);
    \draw[dashed] (73) -- (83);
    \draw[dashed] (83) -- (93);  
    \draw[dashed] (93) -- (103);   
    \draw[dashed] (103) -- (113);   
    \draw[dotted] (113) -- (24,0.866*6);
    \draw[dotted] (04) -- (-1,0.866*8);
    \draw (04) -- (14);  
    \draw[dashed] (14) -- (24);  
    \draw (24) -- (34);  
    \draw[<->] (34) -- (44);  
    \draw[dashed] (44) -- (54); 
    \draw[<->] (54) -- (64);
    \draw[<->] (64) -- (74);
    \draw[dashed] (74) -- (84);
    \draw[<->] (84) -- (94);  
    \draw (94) -- (104);   
    \draw[dashed] (104) -- (114);   
    \draw (114) -- (23.5,0.866*8);
    \draw[dashed] (05) -- (-.5,0.866*10);
    \draw[dashed] (05) -- (15);  
    \draw[dashed] (15) -- (25);  
    \draw[dashed] (25) -- (35);  
    \draw[dashed] (35) -- (45);  
    \draw[dashed] (45) -- (55); 
    \draw[dashed] (55) -- (65);
    \draw[dashed] (65) -- (75);
    \draw[dashed] (75) -- (85);
    \draw[dashed] (85) -- (95);  
    \draw[dashed] (95) -- (105);   
    \draw[dashed] (105) -- (115);   
    \draw[dotted] (115) -- (24,0.866*10);
    \draw[dotted] (06) -- (-1,0.866*12);
    \draw[<->] (06) -- (16);  
    \draw[dashed] (16) -- (26);  
    \draw[<->] (26) -- (36);  
    \draw (36) -- (46);  
    \draw[dashed] (46) -- (56); 
    \draw (56) -- (66);
    \draw (66) -- (76);
    \draw[dashed] (76) -- (86);
    \draw (86) -- (96);  
    \draw[<->] (96) -- (106);   
    \draw[dashed] (106) -- (116);   
    \draw[<->] (116) -- (23.5,0.866*12);
    \draw[dashed] (07) -- (-.5,0.866*14);
    \draw[dashed] (07) -- (17);  
    \draw[dashed] (17) -- (27);  
    \draw[dashed] (27) -- (37);  
    \draw[dashed] (37) -- (47);  
    \draw[dashed] (47) -- (57); 
    \draw[dashed] (57) -- (67);
    \draw[dashed] (67) -- (77);
    \draw[dashed] (77) -- (87);
    \draw[dashed] (87) -- (97);  
    \draw[dashed] (97) -- (107);   
    \draw[dashed] (107) -- (117);   
    \draw[dotted] (117) -- (24,0.866*14);
    \draw[dotted] (00) -- (-.5,-1.5*.866);  
    \draw[<->] (00) -- (.75,-1.5*.866);  
    \draw[dashed] (10) -- (1.25,-1.5*.866);  
    \draw[{Circle}->] (10) -- (2.75,-1.5*.866);   
    \draw[<->] (20) -- (3.25,-1.5*.866); 
    \draw[dashed] (20) -- (4.75,-1.5*.866); 
    \draw[{Circle}-] (30) -- (5.25,-1.5*.866); 
    \draw[{Circle}->] (30) -- (6.75,-1.5*.866); 
    \draw[dashed] (40) -- (7.25,-1.5*.866);   
    \draw (40) -- (8.75,-1.5*.866); 
    \draw[{Circle}-] (50) -- (9.25,-1.5*.866); 
    \draw[dashed] (50) -- (10.75,-1.5*.866);   
    \draw (60) -- (11.25,-1.5*.866);  
    \draw (60) -- (12.75,-1.5*.866);  
    \draw[dashed] (70) -- (13.25,-1.5*.866);  
    \draw[{Circle}-] (70) -- (14.75,-1.5*.866);   
    \draw (80) -- (15.25,-1.5*.866); 
    \draw[dashed] (80) -- (16.75,-1.5*.866); 
    \draw[{Circle}->] (90) -- (17.25,-1.5*.866);   
    \draw[{Circle}-] (90) -- (18.75,-1.5*.866); 
    \draw[dashed] (100) -- (19.25,-1.5*.866);  
    \draw[<->] (100) -- (20.75,-1.5*.866); 
    \draw[{Circle}->] (110) -- (21.25,-1.5*.866);   
    \draw[dashed] (110) -- (22.75,-1.5*.866);  
    \draw[<->] (00) -- (01);  
    \draw[dashed] (10) -- (01);  
    \draw[<->] (10) -- (11);   
    \draw[<->] (20) -- (11); 
    \draw[dashed] (20) -- (21); 
    \draw (30) -- (21);   
    \draw[->] (30) -- (31); 
    \draw[dashed] (40) -- (31);   
    \draw (40) -- (41); 
    \draw[<-] (50) -- (41);   
    \draw[dashed] (50) -- (51);    
    \draw (60) -- (51);    
    \draw (60) -- (61);  
    \draw[dashed] (70) -- (61);  
    \draw[<-] (70) -- (71);   
    \draw (80) -- (71); 
    \draw[dashed] (80) -- (81); 
    \draw[->] (90) -- (81);   
    \draw (90) -- (91); 
    \draw[dashed] (100) -- (91);   
    \draw[<->] (100) -- (101); 
    \draw[<->] (110) -- (101);   
    \draw[dashed] (110) -- (111);    
    \draw[<-] (02) -- (01);  
    \draw[dashed] (12) -- (01);  
    \draw[<-] (12) -- (11);   
    \draw[<-] (22) -- (11); 
    \draw[dashed] (22) -- (21); 
    \draw[->] (32) -- (21);   
    \draw (32) -- (31); 
    \draw[dashed] (42) -- (31);   
    \draw[->] (42) -- (41); 
    \draw[<->] (52) -- (41);   
    \draw[dashed] (52) -- (51); 
    \draw[->] (62) -- (51);  
    \draw[->] (62) -- (61);  
    \draw[dashed] (72) -- (61);  
    \draw[<->] (72) -- (71);   
    \draw[->] (82) -- (71); 
    \draw[dashed] (82) -- (81); 
    \draw (92) -- (81);   
    \draw[->] (92) -- (91); 
    \draw[dashed] (102) -- (91);   
    \draw[<-] (102) -- (101); 
    \draw[<-] (112) -- (101);   
    \draw[dashed] (112) -- (111);    
    \draw (02) -- (03);  
    \draw[dashed] (12) -- (03);  
    \draw[{Circle}-] (12) -- (13);   
    \draw (22) -- (13); 
    \draw[dashed] (22) -- (23); 
    \draw[{Circle}->] (32) -- (23);   
    \draw[{Circle}-] (32) -- (33); 
    \draw[dashed] (42) -- (33);   
    \draw[<->] (42) -- (43); 
    \draw[{Circle}->] (52) -- (43);   
    \draw[dashed] (52) -- (53); 
    \draw[<->] (62) -- (53);  
    \draw[<->] (62) -- (63);  
    \draw[dashed] (72) -- (63);  
    \draw[{Circle}->] (72) -- (73);   
    \draw[<->] (82) -- (73); 
    \draw[dashed] (82) -- (83); 
    \draw[{Circle}-] (92) -- (83);   
    \draw[{Circle}->] (92) -- (93); 
    \draw[dashed] (102) -- (93);   
    \draw (102) -- (103); 
    \draw[{Circle}-] (112) -- (103);   
    \draw[dashed] (112) -- (113);    
    \draw (04) -- (03);  
    \draw[dashed] (14) -- (03);  
    \draw[{Circle}-] (14) -- (13);   
    \draw (24) -- (13); 
    \draw[dashed] (24) -- (23); 
    \draw[{Circle}->] (34) -- (23);   
    \draw[{Circle}-] (34) -- (33); 
    \draw[dashed] (44) -- (33);   
    \draw[<->] (44) -- (43); 
    \draw[{Circle}->] (54) -- (43);   
    \draw[dashed] (54) -- (53); 
    \draw[<->] (64) -- (53);  
    \draw[<->] (64) -- (63);  
    \draw[dashed] (74) -- (63);  
    \draw[{Circle}->] (74) -- (73);   
    \draw[<->] (84) -- (73); 
    \draw[dashed] (84) -- (83); 
    \draw[{Circle}-] (94) -- (83);   
    \draw[{Circle}->] (94) -- (93); 
    \draw[dashed] (104) -- (93);   
    \draw (104) -- (103); 
    \draw[{Circle}-] (114) -- (103);   
    \draw[dashed] (114) -- (113);  
    \draw[dashed] (110) -- (111);    
    \draw[<-] (04) -- (05);  
    \draw[dashed] (14) -- (05);  
    \draw[<-] (14) -- (15);   
    \draw[<-] (24) -- (15); 
    \draw[dashed] (24) -- (25); 
    \draw[->] (34) -- (25);   
    \draw (34) -- (35); 
    \draw[dashed] (44) -- (35);   
    \draw[->] (44) -- (45); 
    \draw[<->] (54) -- (45);   
    \draw[dashed] (54) -- (55); 
    \draw[->] (64) -- (55);  
    \draw[->] (64) -- (65);  
    \draw[dashed] (74) -- (65);  
    \draw[<->] (74) -- (75);   
    \draw[->] (84) -- (75); 
    \draw[dashed] (84) -- (85); 
    \draw (94) -- (85);   
    \draw[->] (94) -- (95); 
    \draw[dashed] (104) -- (95);   
    \draw[<-] (104) -- (105); 
    \draw[<-] (114) -- (105);   
    \draw[dashed] (114) -- (115); 
    \draw[<->] (06) -- (05);  
    \draw[dashed] (16) -- (05);  
    \draw[<->] (16) -- (15);   
    \draw[<->] (26) -- (15); 
    \draw[dashed] (26) -- (25); 
    \draw (36) -- (25);   
    \draw[->] (36) -- (35); 
    \draw[dashed] (46) -- (35);   
    \draw (46) -- (45); 
    \draw[<-] (56) -- (45);   
    \draw[dashed] (56) -- (55);    
    \draw (66) -- (55);    
    \draw (66) -- (65);  
    \draw[dashed] (76) -- (65);  
    \draw[<-] (76) -- (75);   
    \draw (86) -- (75); 
    \draw[dashed] (86) -- (85); 
    \draw[->] (96) -- (85);   
    \draw (96) -- (95); 
    \draw[dashed] (106) -- (95);   
    \draw[<->] (106) -- (105); 
    \draw[<->] (116) -- (105);   
    \draw[dashed] (116) -- (115); 
    \draw[<->] (06) -- (07);  
    \draw[dashed] (16) -- (07);  
    \draw[{Circle}->] (16) -- (17);   
    \draw[<->] (26) -- (17); 
    \draw[dashed] (26) -- (27); 
    \draw[{Circle}-] (36) -- (27);   
    \draw[{Circle}->] (36) -- (37); 
    \draw[dashed] (46) -- (37);   
    \draw (46) -- (47); 
    \draw[{Circle}-] (56) -- (47);   
    \draw[dashed] (56) -- (57); 
    \draw (66) -- (57);  
    \draw (66) -- (67);  
    \draw[dashed] (76) -- (67);  
    \draw[{Circle}-] (76) -- (77);   
    \draw (86) -- (77); 
    \draw[dashed] (86) -- (87); 
    \draw[{Circle}->] (96) -- (87);   
    \draw[{Circle}-] (96) -- (97); 
    \draw[dashed] (106) -- (97);   
    \draw[<->] (106) -- (107); 
    \draw[{Circle}->] (116) -- (107);   
    \draw[dashed] (116) -- (117);  
    \draw[dotted] (07) -- (.5,15*.866);  
    \draw[dotted] (07) -- (1.5,15*.866);  
    \draw[dotted] (17) -- (2.5,15*.866);  
    \draw[dotted] (17) -- (3.5,15*.866);   
    \draw[dotted] (27) -- (4.5,15*.866); 
    \draw[dotted] (27) -- (5.5,15*.866); 
    \draw[dotted] (37) -- (6.5,15*.866);   
    \draw[dotted] (37) -- (7.5,15*.866); 
    \draw[dotted] (47) -- (8.5,15*.866);   
    \draw[dotted] (47) -- (9.5,15*.866); 
    \draw[dotted] (57) -- (10.5,15*.866);   
    \draw[dotted] (57) -- (11.5,15*.866);   
    \draw[dotted] (67) -- (12.5,15*.866);  
    \draw[dotted] (67) -- (13.5,15*.866);  
    \draw[dotted] (77) -- (14.5,15*.866);  
    \draw[dotted] (77) -- (15.5,15*.866);   
    \draw[dotted] (87) -- (16.5,15*.866); 
    \draw[dotted] (87) -- (17.5,15*.866); 
    \draw[dotted] (97) -- (18.5,15*.866);   
    \draw[dotted] (97) -- (19.5,15*.866); 
    \draw[dotted] (107) -- (20.5,15*.866);   
    \draw[dotted] (107) -- (21.5,15*.866); 
    \draw[dotted] (117) -- (22.5,15*.866);   
    \draw[dotted] (117) -- (23.5,15*.866);  
\end{tikzpicture} 

\, \\

\begin{tikzpicture}[scale=.45]
    \node (00) at (0,0)    {1};
    \node (10) at (2,0)    {3};  
    \node (20) at (4,0)    {9};     
    \node (30) at (6,0)    {19};    
    \node (40) at (8,0)    {9};     
    \node (50) at (10,0)    {3};   
    \node (60) at (12, 0)   {1};   
    \node (70) at (14, 0)   {3};
    \node (80) at (16, 0)   {9};
    \node (90) at (18, 0)   {19};
    \node (100) at (20, 0)   {9};
    \node (110) at (22, 0)   {3};
    \node (01) at (1,0.866*2)    {18};
    \node (11) at (3,0.866*2)    {10};  
    \node (21) at (5,0.866*2)    {6};     
    \node (31) at (7,0.866*2)    {6};    
    \node (41) at (9,0.866*2)    {10};     
    \node (51) at (11,0.866*2)    {18}; 
    \node (61) at (13,0.866*2)    {6};
    \node (71) at (15,0.866*2)    {22};  
    \node (81) at (17,0.866*2)    {18};     
    \node (91) at (19,0.866*2)    {18};    
    \node (101) at (21,0.866*2)    {22};     
    \node (111) at (23,0.866*2)    {6}; 
    \node (02) at (0,0.866*4)    {13};
    \node (12) at (2,0.866*4)    {15};  
    \node (22) at (4,0.866*4)    {21};     
    \node (32) at (6,0.866*4)    {7};    
    \node (42) at (8,0.866*4)    {21};     
    \node (52) at (10,0.866*4)    {15};   
    \node (62) at (12,0.866*4)    {13};
    \node (72) at (14,0.866*4)    {15};  
    \node (82) at (16,0.866*4)    {21};     
    \node (92) at (18,0.866*4)    {7};    
    \node (102) at (20,0.866*4)    {21};     
    \node (112) at (22,0.866*4)    {15};    
    \node (03) at (1,0.866*6)    {0};
    \node (13) at (3,0.866*6)    {4};  
    \node (23) at (5,0.866*6)    {12};     
    \node (33) at (7,0.866*6)    {0};    
    \node (43) at (9,0.866*6)    {16};     
    \node (53) at (11,0.866*6)    {12}; 
    \node (63) at (13,0.866*6)    {12};
    \node (73) at (15,0.866*6)    {16};  
    \node (83) at (17,0.866*6)    {0};     
    \node (93) at (19,0.866*6)    {12};    
    \node (103) at (21,0.866*6)    {4};     
    \node (113) at (23,0.866*6)    {0}; 
    \node (04) at (0,0.866*8)    {13};
    \node (14) at (2,0.866*8)    {15};  
    \node (24) at (4,0.866*8)    {21};     
    \node (34) at (6,0.866*8)    {7};    
    \node (44) at (8,0.866*8)    {21};     
    \node (54) at (10,0.866*8)    {15};   
    \node (64) at (12,0.866*8)    {13};
    \node (74) at (14,0.866*8)    {15};  
    \node (84) at (16,0.866*8)    {21};     
    \node (94) at (18,0.866*8)    {7};    
    \node (104) at (20,0.866*8)    {21};     
    \node (114) at (22,0.866*8)    {15};  
    \node (05) at (1,0.866*10)    {18};
    \node (15) at (3,0.866*10)    {10};  
    \node (25) at (5,0.866*10)    {6};     
    \node (35) at (7,0.866*10)    {6};    
    \node (45) at (9,0.866*10)    {10};     
    \node (55) at (11,0.866*10)    {18}; 
    \node (65) at (13,0.866*10)    {6};
    \node (75) at (15,0.866*10)    {22};  
    \node (85) at (17,0.866*10)    {18};     
    \node (95) at (19,0.866*10)    {18};    
    \node (105) at (21,0.866*10)    {22};     
    \node (115) at (23,0.866*10)    {6}; 
    \node (06) at (0,0.866*12)    {1};
    \node (16) at (2,0.866*12)    {3};  
    \node (26) at (4,0.866*12)    {9};     
    \node (36) at (6,0.866*12)    {19};    
    \node (46) at (8,0.866*12)    {9};     
    \node (56) at (10,0.866*12)    {3};   
    \node (66) at (12,0.866*12)    {1};
    \node (76) at (14,0.866*12)    {3};  
    \node (86) at (16,0.866*12)    {9};     
    \node (96) at (18,0.866*12)    {19};    
    \node (106) at (20,0.866*12)    {9};     
    \node (116) at (22,0.866*12)    {3};  
    \node (07) at (1,0.866*14)    {0};
    \node (17) at (3,0.866*14)    {4};  
    \node (27) at (5,0.866*14)    {12};     
    \node (37) at (7,0.866*14)    {0};    
    \node (47) at (9,0.866*14)    {16};     
    \node (57) at (11,0.866*14)    {12}; 
    \node (67) at (13,0.866*14)    {12};
    \node (77) at (15,0.866*14)    {16};  
    \node (87) at (17,0.866*14)    {0};     
    \node (97) at (19,0.866*14)    {12};    
    \node (107) at (21,0.866*14)    {4};     
    \node (117) at (23,0.866*14)    {0}; 
    \draw[dotted] (00) -- (-1,0);
    \draw[<->] (00) -- (10);
    \draw[dashed] (10) -- (20);
    \draw[<->] (20) -- (30);
    \draw (30) -- (40);  
    \draw[dashed] (40) -- (50);  
    \draw (50) -- (60);
    \draw (60) -- (70);
    \draw[dashed] (70) -- (80);
    \draw (80) -- (90);  
    \draw[<->] (90) -- (100);  
    \draw[dashed] (100) -- (110); 
    \draw[<->] (110) -- (23.5,0);
    \draw[dashed] (01) -- (-.5,0.866*2);
    \draw[dashed] (01) -- (11);  
    \draw[dashed] (11) -- (21);  
    \draw[dashed] (21) -- (31);  
    \draw[dashed] (31) -- (41);  
    \draw[dashed] (41) -- (51); 
    \draw[dashed] (51) -- (61);  
    \draw[dashed] (61) -- (71);  
    \draw[dashed] (71) -- (81);  
    \draw[dashed] (81) -- (91);  
    \draw[dashed] (91) -- (101); 
    \draw[dashed] (101) -- (111); 
    \draw[dotted] (111) -- (24,0.866*2);
    \draw[dotted] (02) -- (-1,0.866*4);
    \draw (02) -- (12);
    \draw[dashed] (12) -- (22);
    \draw (22) -- (32);
    \draw (32)[<->] -- (42);  
    \draw[dashed] (42) -- (52);   
    \draw (52)[<->] -- (62);
    \draw (62)[<->] -- (72);
    \draw[dashed] (72) -- (82);
    \draw (82)[<->] -- (92);  
    \draw (92) -- (102);   
    \draw[dashed] (102) -- (112);   
    \draw (112) -- (23.5,0.866*4);
    \draw[dashed] (03) -- (-.5,0.866*6);
    \draw[dashed] (03) -- (13);  
    \draw[dashed] (13) -- (23);  
    \draw[dashed] (23) -- (33);  
    \draw[dashed] (33) -- (43);  
    \draw[dashed] (43) -- (53); 
    \draw[dashed] (53) -- (63);
    \draw[dashed] (63) -- (73);
    \draw[dashed] (73) -- (83);
    \draw[dashed] (83) -- (93);  
    \draw[dashed] (93) -- (103);   
    \draw[dashed] (103) -- (113);   
    \draw[dotted] (113) -- (24,0.866*6);
    \draw[dotted] (04) -- (-1,0.866*8);
    \draw (04) -- (14);  
    \draw[dashed] (14) -- (24);  
    \draw (24) -- (34);  
    \draw (34)[<->] -- (44);  
    \draw[dashed] (44) -- (54); 
    \draw (54)[<->] -- (64);
    \draw (64)[<->] -- (74);
    \draw[dashed] (74) -- (84);
    \draw (84)[<->] -- (94);  
    \draw (94) -- (104);   
    \draw[dashed] (104) -- (114);   
    \draw (114) -- (23.5,0.866*8);
    \draw[dashed] (05) -- (-.5,0.866*10);
    \draw[dashed] (05) -- (15);  
    \draw[dashed] (15) -- (25);  
    \draw[dashed] (25) -- (35);  
    \draw[dashed] (35) -- (45);  
    \draw[dashed] (45) -- (55); 
    \draw[dashed] (55) -- (65);
    \draw[dashed] (65) -- (75);
    \draw[dashed] (75) -- (85);
    \draw[dashed] (85) -- (95);  
    \draw[dashed] (95) -- (105);   
    \draw[dashed] (105) -- (115);   
    \draw[dotted] (115) -- (24,0.866*10);
    \draw[dotted] (06) -- (-1,0.866*12);
    \draw (06)[<->] -- (16);  
    \draw[dashed] (16) -- (26);  
    \draw (26)[<->] -- (36);  
    \draw (36) -- (46);  
    \draw[dashed] (46) -- (56); 
    \draw (56) -- (66);
    \draw (66) -- (76);
    \draw[dashed] (76) -- (86);
    \draw (86) -- (96);  
    \draw[<->] (96) -- (106);   
    \draw[dashed] (106) -- (116);   
    \draw[<->] (116) -- (23.5,0.866*12);
    \draw[dashed] (07) -- (-.5,0.866*14);
    \draw[dashed] (07) -- (17);  
    \draw[dashed] (17) -- (27);  
    \draw[dashed] (27) -- (37);  
    \draw[dashed] (37) -- (47);  
    \draw[dashed] (47) -- (57); 
    \draw[dashed] (57) -- (67);
    \draw[dashed] (67) -- (77);
    \draw[dashed] (77) -- (87);
    \draw[dashed] (87) -- (97);  
    \draw[dashed] (97) -- (107);   
    \draw[dashed] (107) -- (117);   
    \draw[dotted] (117) -- (24,0.866*14);
    \draw[dotted] (00) -- (-.5,-1.5*.866);  
    \draw[<->] (00) -- (.75,-1.5*.866);  
    \draw[dashed] (10) -- (1.25,-1.5*.866);  
    \draw[{Circle}->] (10) -- (2.75,-1.5*.866);   
    \draw[<->] (20) -- (3.25,-1.5*.866); 
    \draw[dashed] (20) -- (4.75,-1.5*.866); 
    \draw[{Circle}-] (30) -- (5.25,-1.5*.866); 
    \draw[{Circle}->] (30) -- (6.75,-1.5*.866); 
    \draw[dashed] (40) -- (7.25,-1.5*.866);   
    \draw (40) -- (8.75,-1.5*.866); 
    \draw[{Circle}-] (50) -- (9.25,-1.5*.866); 
    \draw[dashed] (50) -- (10.75,-1.5*.866);   
    \draw (60) -- (11.25,-1.5*.866);  
    \draw (60) -- (12.75,-1.5*.866);  
    \draw[dashed] (70) -- (13.25,-1.5*.866);  
    \draw[{Circle}-] (70) -- (14.75,-1.5*.866);   
    \draw (80) -- (15.25,-1.5*.866); 
    \draw[dashed] (80) -- (16.75,-1.5*.866); 
    \draw[{Circle}->] (90) -- (17.25,-1.5*.866);   
    \draw[{Circle}-] (90) -- (18.75,-1.5*.866); 
    \draw[dashed] (100) -- (19.25,-1.5*.866);  
    \draw[<->] (100) -- (20.75,-1.5*.866); 
    \draw[{Circle}->] (110) -- (21.25,-1.5*.866);   
    \draw[dashed] (110) -- (22.75,-1.5*.866); 
    \draw[<->] (00) -- (01);  
    \draw[dashed] (10) -- (01);  
    \draw[->] (10) -- (11);   
    \draw[<->] (20) -- (11); 
    \draw[dashed] (20) -- (21); 
    \draw (30) -- (21);   
    \draw[<->] (30) -- (31); 
    \draw[dashed] (40) -- (31);   
    \draw (40) -- (41); 
    \draw[<-] (50) -- (41);   
    \draw[dashed] (50) -- (51);    
    \draw (60) -- (51);    
    \draw (60) -- (61);  
    \draw[dashed] (70) -- (61);  
    \draw (70) -- (71);   
    \draw (80) -- (71); 
    \draw[dashed] (80) -- (81); 
    \draw[->] (90) -- (81);   
    \draw[<-] (90) -- (91); 
    \draw[dashed] (100) -- (91);   
    \draw[<->] (100) -- (101); 
    \draw[<->] (110) -- (101);   
    \draw[dashed] (110) -- (111);    
    \draw[<-] (02) -- (01);  
    \draw[dashed] (12) -- (01);  
    \draw (12) -- (11);   
    \draw[<-] (22) -- (11); 
    \draw[dashed] (22) -- (21); 
    \draw[->] (32) -- (21);   
    \draw[<-] (32) -- (31); 
    \draw[dashed] (42) -- (31);   
    \draw[->] (42) -- (41); 
    \draw[<->] (52) -- (41);   
    \draw[dashed] (52) -- (51); 
    \draw[->] (62) -- (51);  
    \draw[->] (62) -- (61);  
    \draw[dashed] (72) -- (61);  
    \draw[->] (72) -- (71);   
    \draw[->] (82) -- (71); 
    \draw[dashed] (82) -- (81); 
    \draw (92) -- (81);   
    \draw[<->] (92) -- (91); 
    \draw[dashed] (102) -- (91);   
    \draw[<-] (102) -- (101); 
    \draw[<-] (112) -- (101);   
    \draw[dashed] (112) -- (111);    
    \draw (02) -- (03);  
    \draw[dashed] (12) -- (03);  
    \draw[{Circle}-] (12) -- (13);   
    \draw (22) -- (13); 
    \draw[dashed] (22) -- (23); 
    \draw[{Circle}->] (32) -- (23);   
    \draw[{Circle}-] (32) -- (33); 
    \draw[dashed] (42) -- (33);   
    \draw[->] (42) -- (43); 
    \draw[{Circle}->] (52) -- (43);   
    \draw[dashed] (52) -- (53); 
    \draw[<->] (62) -- (53);  
    \draw[<->] (62) -- (63);  
    \draw[dashed] (72) -- (63);  
    \draw[{Circle}->] (72) -- (73);   
    \draw[<->] (82) -- (73); 
    \draw[dashed] (82) -- (83); 
    \draw[{Circle}-] (92) -- (83);   
    \draw[{Circle}->] (92) -- (93); 
    \draw[dashed] (102) -- (93);   
    \draw (102) -- (103); 
    \draw[{Circle}-] (112) -- (103);   
    \draw[dashed] (112) -- (113);    
    \draw (04) -- (03);  
    \draw[dashed] (14) -- (03);  
    \draw[{Circle}-] (14) -- (13);   
    \draw (24) -- (13); 
    \draw[dashed] (24) -- (23); 
    \draw[{Circle}->] (34) -- (23);   
    \draw[{Circle}-] (34) -- (33); 
    \draw[dashed] (44) -- (33);   
    \draw[->] (44) -- (43); 
    \draw[{Circle}->] (54) -- (43);   
    \draw[dashed] (54) -- (53); 
    \draw[<->] (64) -- (53);  
    \draw[<->] (64) -- (63);  
    \draw[dashed] (74) -- (63);  
    \draw[{Circle}->] (74) -- (73);   
    \draw[<->] (84) -- (73); 
    \draw[dashed] (84) -- (83); 
    \draw[{Circle}-] (94) -- (83);   
    \draw[{Circle}->] (94) -- (93); 
    \draw[dashed] (104) -- (93);   
    \draw (104) -- (103); 
    \draw[{Circle}-] (114) -- (103);   
    \draw[dashed] (114) -- (113);   
    \draw[<-] (04) -- (05);  
    \draw[dashed] (14) -- (05);  
    \draw (14) -- (15);   
    \draw[<-] (24) -- (15); 
    \draw[dashed] (24) -- (25); 
    \draw[->] (34) -- (25);   
    \draw[<-] (34) -- (35); 
    \draw[dashed] (44) -- (35);   
    \draw[->] (44) -- (45); 
    \draw[<->] (54) -- (45);   
    \draw[dashed] (54) -- (55); 
    \draw[->] (64) -- (55);  
    \draw[->] (64) -- (65);  
    \draw[dashed] (74) -- (65);  
    \draw[->] (74) -- (75);   
    \draw[->] (84) -- (75); 
    \draw[dashed] (84) -- (85); 
    \draw (94) -- (85);   
    \draw[<->] (94) -- (95); 
    \draw[dashed] (104) -- (95);   
    \draw[<-] (104) -- (105); 
    \draw[<-] (114) -- (105);   
    \draw[dashed] (114) -- (115);  
    \draw[<->] (06) -- (05);  
    \draw[dashed] (16) -- (05);  
    \draw[->] (16) -- (15);   
    \draw[<->] (26) -- (15); 
    \draw[dashed] (26) -- (25); 
    \draw (36) -- (25);   
    \draw[<->] (36) -- (35); 
    \draw[dashed] (46) -- (35);   
    \draw (46) -- (45); 
    \draw[<-] (56) -- (45);   
    \draw[dashed] (56) -- (55);    
    \draw (66) -- (55);    
    \draw (66) -- (65);  
    \draw[dashed] (76) -- (65);  
    \draw (76) -- (75);   
    \draw (86) -- (75); 
    \draw[dashed] (86) -- (85); 
    \draw[->] (96) -- (85);   
    \draw[<-] (96) -- (95); 
    \draw[dashed] (106) -- (95);   
    \draw[<->] (106) -- (105); 
    \draw[<->] (116) -- (105);   
    \draw[dashed] (116) -- (115);    
    \draw[<->] (06) -- (07);  
    \draw[dashed] (16) -- (07);  
    \draw[{Circle}->] (16) -- (17);   
    \draw[<->] (26) -- (17); 
    \draw[dashed] (26) -- (27); 
    \draw[{Circle}-] (36) -- (27);   
    \draw[{Circle}->] (36) -- (37); 
    \draw[dashed] (46) -- (37);   
    \draw (46) -- (47); 
    \draw[{Circle}-] (56) -- (47);   
    \draw[dashed] (56) -- (57); 
    \draw (66) -- (57);  
    \draw (66) -- (67);  
    \draw[dashed] (76) -- (67);  
    \draw[{Circle}-] (76) -- (77);   
    \draw (86) -- (77); 
    \draw[dashed] (86) -- (87); 
    \draw[{Circle}->] (96) -- (87);   
    \draw[{Circle}-] (96) -- (97); 
    \draw[dashed] (106) -- (97);   
    \draw[<->] (106) -- (107); 
    \draw[{Circle}->] (116) -- (107);   
    \draw[dashed] (116) -- (117);  
    \draw[dotted] (07) -- (.5,15*.866);  
    \draw[dotted] (07) -- (1.5,15*.866);  
    \draw[dotted] (17) -- (2.5,15*.866);  
    \draw[dotted] (17) -- (3.5,15*.866);   
    \draw[dotted] (27) -- (4.5,15*.866); 
    \draw[dotted] (27) -- (5.5,15*.866); 
    \draw[dotted] (37) -- (6.5,15*.866);   
    \draw[dotted] (37) -- (7.5,15*.866); 
    \draw[dotted] (47) -- (8.5,15*.866);   
    \draw[dotted] (47) -- (9.5,15*.866); 
    \draw[dotted] (57) -- (10.5,15*.866);   
    \draw[dotted] (57) -- (11.5,15*.866);   
    \draw[dotted] (67) -- (12.5,15*.866);  
    \draw[dotted] (67) -- (13.5,15*.866);  
    \draw[dotted] (77) -- (14.5,15*.866);  
    \draw[dotted] (77) -- (15.5,15*.866);   
    \draw[dotted] (87) -- (16.5,15*.866); 
    \draw[dotted] (87) -- (17.5,15*.866); 
    \draw[dotted] (97) -- (18.5,15*.866);   
    \draw[dotted] (97) -- (19.5,15*.866); 
    \draw[dotted] (107) -- (20.5,15*.866);   
    \draw[dotted] (107) -- (21.5,15*.866); 
    \draw[dotted] (117) -- (22.5,15*.866);   
    \draw[dotted] (117) -- (23.5,15*.866);  
\end{tikzpicture} 
\end{center}
\caption{Pattern of Kronecker Symbols in Triangular Packings of Type $(2, 5,8,11)$; Pattern in Type $(0,1,3,4,6,7,9,10)$}
\label{TrianglePartialEvenFigure}
\end{figure}
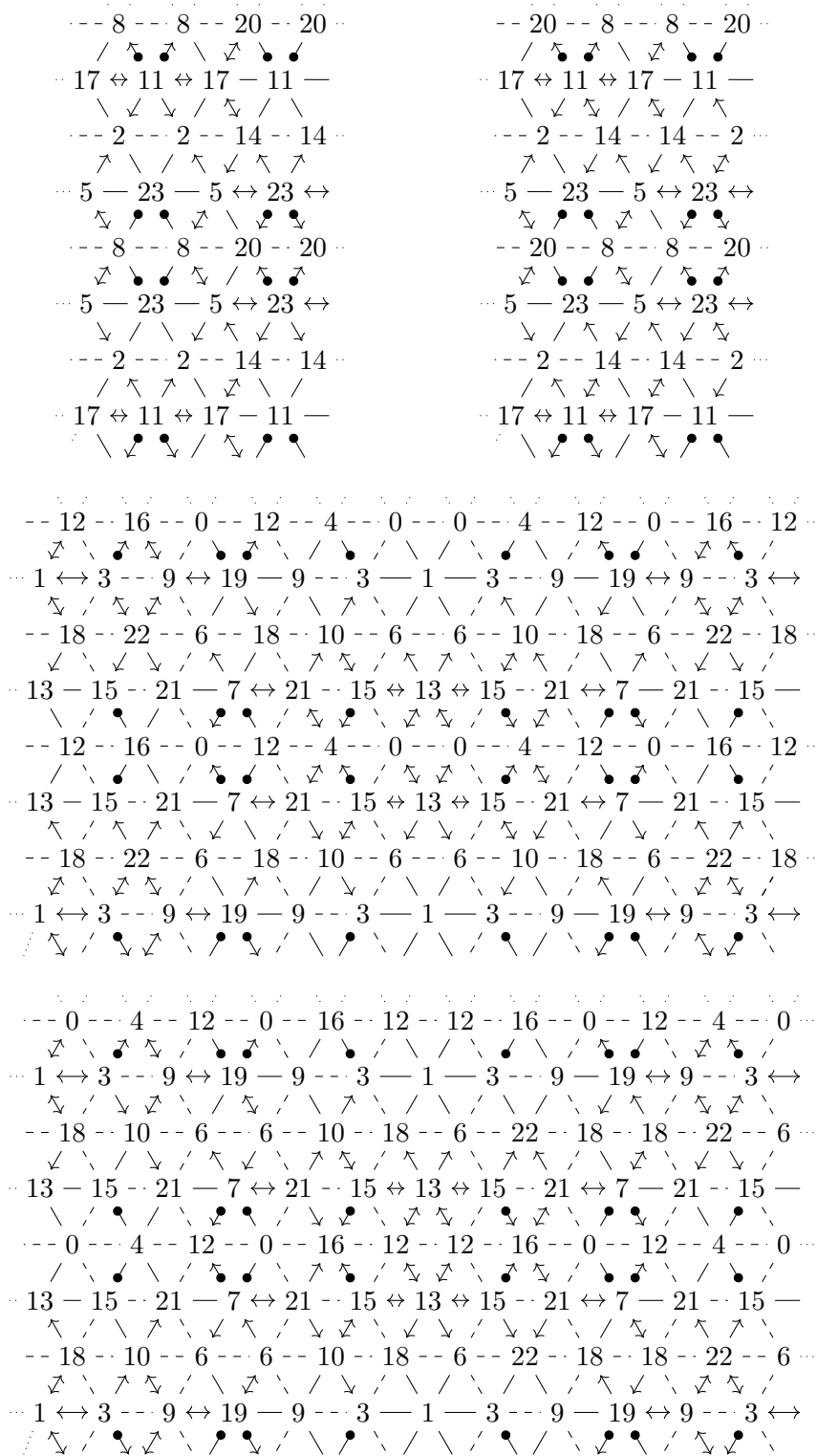

In types $(7)$ and $(5,9)$ we do not find a way of defining $\chi_2$ consistently for all circles in the packing, and we do not find any quadratic obstructions. However, it is still the case that Kronecker symbols between adjacent curvatures have a predictable pattern, as illustrated in Figure \ref{TrianglePartialFigure}. In this figure, all pairs of tangent, coprime curvatures $(a,b)$ with an arrow pointing from $a$ to $b$ have the same $\res{b}{a}$ value. All pairs $(a,b)$ without an arrow have the opposite $\res{b}{a}$ value. All pairs with a dashed line are not coprime.

The behavior of Kronecker symbols is governed by Proposition \ref{TriangleNode} and quadratic reciprocity. Because each pattern is periodic, with all arrows completely determined by one triangle of arrows, and the generators of the triangular Apollonian group all fix three circles arranged in a triangle, the Apollonian group action is compatible with this pattern. Even though some adjacent pairs of curvatures not be coprime, Proposition \ref{TrianglePath} implies that we can follow the pattern of Kronecker symbols on a path from one to the other.

In type $(7)$, we do not find partial obstructions to the Local-Global Conjecture. But in type $(5,9)$ there is a partial obstruction analogous to the partial obstruction for octahedral packings of type $(0,3,6)$ or $(4,5,6)$ mod $8$ in Proposition \ref{OctPartialObstruction}. It is possible to three-color the vertices of the triangular tiling such that adjacent circles have different colors. This gives rise to the unique three-coloring of a triangular packing such that tangent circles have different colors. We may color the circles so that all circles of curvature $5 \bmod 12$ are yellow, and circles of curvature $9 \bmod 12$ are red and blue. Then $\res{a}{b}$ takes one value for all coprime tangent pairs with $a$ red, $b$ yellow, and takes the opposite value for all coprime tangent pairs with $a$ blue, $b$ yellow. Therefore, there are either no blue circles of curvature $n^2$ or no red circles of curvature $n^2$. 

We now describe the situation in types $(2, 5, 8, 11)$ and $(0,1,3,4,6,7,9,10)$. In order to make quadratic reciprocity calculations involving even integers, it is helpful to lift from mod 12 to mod 24. The two possible liftings from mod 4 to mod 8 give rise to two possible liftings from mod 12 to mod 24 for each configuration containing even and odd curvatures. These are illustrated in Figure \ref{TrianglePartialEvenFigure} below. We also illustrate the pattern of Kronecker symbols in each case. All pairs of tangent, coprime curvatures $(a,b)$ with an arrow pointing from $a$ to $b$ have the same $\res{b}{a'}$ value, where $a'=a/2$ if $a \equiv 2 \bmod 4$, and $a'=a$ otherwise. All pairs $(a,b)$ without an arrow have the opposite $\res{b}{a'}$ value. All pairs with a dashed line are not coprime. A circular arrowhead is shown in situations where we need more information about the factorization of $a$ to determine $\res{b}{a'}$. 

The pattern is governed by quadratic reciprocity and Proposition \ref{TriangleNode}. We observe the same rigidity seen in previous cases: if one arrow is known in the configuration, then all other arrows (except those with circular arrowheads) in the configuration are determined. Because the triangular Apollonian group generators fix three circles arranged in a triangle, the arrows across the full packing are determined.

For each mod 12 configuration, the two possible lifts have very different behavior with respect to quadratic invariants and obstructions. In the first lift for each type, the value of $\res{b}{a'}$ is determined by the values of $a, b \bmod 24$, but this is not the case in the second lift. Thus it may be possible to define a quadratic invariant for the first lift, but not the second. In the first lift of type $(2, 5, 8, 11)$, the value of $\res{2b}{a'}$ is the same for all tangent, coprime pairs $a$, $b$ with $b \equiv 2, 8 \bmod 24$. If this value is $-1$, then no circles of curvature $2n^2$ appear in the configuration. However, in the second lift, this does not hold, and we do not find quadratic obstructions. Similarly, in the first lift of type $(0,1,3,4,6,7,9,10)$, the value of $\res{6b}{a'}$ is the same for all tangent, coprime pairs $a$, $b$ with $b \equiv 0, 6 \bmod 24$. If this value is $-1$, then no circles of curvature $6n^2$ appear in the configuration. However, in the second lift, this does not hold, and we do not find quadratic obstructions. 

Because the Apollonian group acts transitively on the lifts, an obstruction within one configuration does not imply an obstruction across the entire packing. In fact, using a sequence of Apollonian group generators which fix a triangle of curvatures while mapping between the two lifts, we can see that every packing contains four distinct types of configurations mod 24: the two lifts shown in Figure \ref{TrianglePartialEvenFigure} and two additional copies of these lifts with all the arrows reversed. We have quadratic obstructions within one of these four configuration types, but not the other three, and not across the packing. 

\pagebreak

\section{Data} \label{Data}

We gathered data for 32 examples of octahedral, cubic, square, and triangular packings, including two examples for each modular restriction type. In the tables below, packings are labeled by the curvatures present in their base configurations, their modular type, their $\chi_2$ value if applicable, and any quadratic obstructions. We computed all curvatures present in the packing up to $N$. We define $S_P(N)$ as the set of sporadic integers, i.e. positive integers not ruled out by modular restrictions or quadratic obstructions, which do not appear in the packing. A large value in the right-hand column means a greater likelihood that all sporadic integers were found. This data confirms the quadratic obstructions that we found, and suggests that there may not be any other obstructions.

\begin{center}
    Octahedral Packings
\end{center}
\begin{adjustbox}{center}
\begin{tabular}[h!]{|p{1.5in}|P{1.0in}|P{0.2in}|P{0.6in}|P{0.6in}|P{0.8in}|P{0.5in}|P{0.7in}|}
\hline
 Base Configuration & Modular Type & $\chi_2$ & Obst. & $N$ & Max($S_P(N)$) & $|S_P(N)|$ & $\frac{N}{\text{Max($S_P(N)$)}}$\\ \hline
 (-7,16,16,18,18,41) & (0,1,2) &1 & & 120000000 & 64037994 & 23123 & 1.87 \\ 
 (-6,10,17,17,24,40) & &-1 & $n^2, 2n^2$&80000000 & 18792216 & 11709 & 4.26 \\
 \hline
 (-2,3,6,8,11,16) & (0,3,6) & & &1000000 & 116523 & 125 & 8.58 \\ 
 (-5,8,14,16,22,35) & & & &4000000 & 1235310 & 785 & 3.24 \\
 \hline
 (-1,2,2,4,4,7) & (2,4,7) &1 & &1000000 & 34716 & 31 & 28.81 \\ 
 (-4,7,10,12,15,26) & &-1 & $n^2, 2n^2$ &40000000 & 314388 & 339 & 127.23 \\
 \hline
 (-2,4,5,5,6,12) & (4,5,6) & & &1000000 & 138798 & 95 & 7.20 \\ 
 (-4,6,13,13,20,30) & & & &80000000 & 6631038 & 2750 & 12.06 \\
 \hline
\end{tabular}
\end{adjustbox}

\begin{center}
    Cubic Packings
\end{center}
\begin{adjustbox}{center}
\begin{tabular}[h!]{|p{1.5in}|P{1in}|P{0.2in}|P{0.6in}|P{0.6in}|P{0.8in}|P{0.5in}|P{0.7in}|}
\hline
 Base Configuration & Modular Type & $\chi_2$ & Obst. & $N$ & Max($S_P(N)$) & $|S_P(N)|$ & $\frac{N}{\text{Max($S_P(N)$)}}$\\ \hline
 \hspace{-.05in}(-7,\hspace{-.01in}16,\hspace{-.01in}18,\hspace{-.01in}25,\hspace{-.01in}41,\hspace{-.01in}48,\hspace{-.01in}50,\hspace{-.01in}73) & (0,1,2) & 1 & &1000000 & 345414 & 1890 & 2.90 \\ 
 (-2,5,5,6,12,13,13,20) & & -1 & $n^2, 2n^2$ &1000000 & 12336 & 118 & 81.06 \\ 
 \hline
 (-1,2,3,4,6,7,8,11) & (0,2,3) & & &1000000 & 312 & 7 & 3205.13 \\
  \hspace{-.05in}(-2,3,10,11,15,16,23,28) & & & &1000000 & 17874 & 168 & 55.95 \\ 
 \hline
\end{tabular}
\end{adjustbox}

\begin{center}
    Square Packings
\end{center}
\begin{adjustbox}{center}
\begin{tabular}[h!]{|p{1.5in}|P{1.0in}|P{0.2in}|P{0.6in}|P{0.6in}|P{0.8in}|P{0.5in}|P{0.7in}|}
\hline
 Base Configuration & Modular Type & $\chi_2$ & Obst. & $N$ & Max($S_P(N)$) & $|S_P(N)|$ & $\frac{N}{\text{Max($S_P(N)$)}}$\\ \hline
 (1,1,1,1) & (1) & 1 & & 100000 & - & 0 & - \\
 (-7,17,17,41) & & -1 & $n^2$ & 100000 & 2665 & 4 & 37.52 \\ 
 \hline
 (-3,5,13,21) & (5) & & & 100000 & - & 0 & - \\
 (-27,37,173,237) & & & & 100000 & 30629 & 201 & 3.26 \\ 
 \hline
 (-3,5,12,20) & full & 1 & & 100000 & 7297 & 390 & 13.70 \\
 (-1,2,3,6) & & -1 & $n^2$ & 100000 & 154 & 10 & 649.35 \\ 
 \hline
 (-1,3,3,7) & (3,7) & & & 100000 & - & 0 & - \\
 (-5,7,31,43) & & & & 100000 & 3827 & 53 & 26.13 \\ 
 \hline
\end{tabular}
\end{adjustbox}

\pagebreak

\begin{center}
    Triangular Packings
\end{center}
\begin{adjustbox}{center}
\begin{tabular}[h!]{|p{1.5in}|P{1.0in}|P{0.2in}|P{0.6in}|P{0.6in}|P{0.8in}|P{0.5in}|P{0.7in}|}
\hline
 Base Configuration & Modular Type & $\chi_2$ & Obst. & $N$ & Max($S_P(N)$) & $|S_P(N)|$ & $\frac{N}{\text{Max($S_P(N)$)}}$\\ \hline
 (1,1,1) & (1) & 1 & & 100000 & - & 0 & -  \\
 (-11,13,73) & & -1 & $n^2$ & 100000 & 27157 & 45 & 3.68 \\ 
 \hline
 (-5,7,19) & (7) & & & 100000 & 175 & 1 & 571.43 \\
 (-17,31,43) & & & & 100000 & 4699 & 17 & 21.28 \\ 
 \hline
 (-1,3,3) & (3,11) & 1 & & 100000 & 1127 & 1 & 88.73 \\
 (-9,15,23) & & -1 & $3n^2$ & 100000 & 15275 & 39 & 6.55 \\ 
 \hline
 (-3,5,9) & (5,9) & & & 100000 & - & 0 & - \\
 (-7,9,33) & & & & 100000 & 67301 & 159 & 1.49 \\ 
 \hline
 (-2,3,6) & (0,1,3,4,6,7,9,10) & & & 100000 & 15106 & 153 & 6.62 \\
 (-3,6,7) & & & & 100000 & 7993 & 97 & 12.51 \\ 
 \hline
 (-1,2,2) & (2,5,8,11) & & & 100000 & - & 0 & - \\
 (-4,5,20) & & & & 100000 & 33923 & 193 & 2.95  \\ 
 \hline
\end{tabular}
\end{adjustbox}

\, \\

\printbibliography
\end{document}